\documentclass[final,3p,times]{elsarticle}

\usepackage{amssymb}
\usepackage{amsmath}
\usepackage{amsthm}
\usepackage{mathrsfs}
\newtheorem{proposition}{Proposition}[section]
\theoremstyle{remark}\newtheorem{remark}{Remark}
\usepackage{mathtools}
\usepackage{hyperref}
\usepackage{graphicx}
\usepackage{subcaption}
\usepackage{float}
\usepackage{lineno}
\usepackage{booktabs}
\usepackage{multirow}  
\usepackage{caption,graphicx,newfloat}
\DeclareCaptionType{InfoBox}
\usepackage{enumitem}
\usepackage{algorithm}
\usepackage{algpseudocode} 

\usepackage{tikz}
\usetikzlibrary{arrows.meta,positioning,fit,calc,shapes.geometric,shapes.misc,backgrounds}
\tikzset{>=Latex}

\begin{document}

\begin{frontmatter}

\title{Neural Spline Flow Diffeomorphisms for Pullback Finite Element Methods on Complex Domains}

\author[label1]{Ioannis Kalogeris}
\author[label1]{Anaksagoras Spitas}   
\author[label2]{Konstantinos Atzarakis} 
\author[label2]{George Stavroulakis} 
\author[label2]{Vissarion Papadopoulos} 

\affiliation[label1]{organization={School of Mechanical Engineering, National Technical University of Athens},
            addressline={Iroon Polytechniou 9}, 
            city={Zografou Campus},
            state={Athens},
            country={Greece}}

\affiliation[label2]{organization={School of Civil Engineering, National Technical University of Athens},
            addressline={Iroon Polytechniou 9}, 
            city={Zografou Campus},
            state={Athens},
            country={Greece}}

\begin{abstract}
The present work proposes a method for solving partial differential equations on arbitrary domains by learning smooth, invertible mappings from a predefined reference domain to the physical domain using Neural Spline Flows (NSF). The learned map is a diffeomorphism by construction, since the NSF architecture guarantees a strictly positive Jacobian determinant through monotone rational-quadratic spline coupling layers, thus no Jacobian-sign penalty is required. Through a diffeomorphic pullback, the weak form of the partial differential equation is transferred from the complex physical domain to the reference domain, where standard structured finite element meshes can be used directly. Specifically, an auto-decoder framework is adopted in which a shared NSF learns a family of deformations from one domain to the other, conditioned on a low-dimensional latent code. For a given learned diffeomorphism, the pullback establishes a rigorous equivalence between the reference domain and the domain represented by the learned map, preserving the mathematical structure of the original problem. Since the learned geometry generally approximates, rather than reproduces exactly, the target physical domain, an additional geometric mapping error is introduced, whose magnitude depends on the accuracy and distortion of the learned transformation. Nevertheless, the numerical results show that the geometric approximation can be made sufficiently accurate to retain reliable finite element solutions. This formulation offers the advantage of eliminating the need to generate unstructured meshes on complex geometries, which often poses a significant bottleneck in the simulation workflow. An additional advantage is that the reference-domain formulation opens the path to low-level computational optimizations, such as memory-aligned data access, efficient matrix assembly, and precomputed quadrature, which are difficult to exploit on irregular meshes. By combining the rigorous mathematical foundations of classical numerical analysis with the representational flexibility of modern machine learning, the proposed approach offers a novel solution framework in which the neural network handles the geometric complexity while the finite element method retains its conventional discretization and convergence properties on the learned domain, with the overall solution accuracy additionally influenced by the quality of the learned geometric mapping.
\end{abstract}



\end{frontmatter}

\section{Introduction}
\label{sec:intro}
 
The numerical solution of partial differential equations (PDEs) on complex geometries is a cornerstone of computational mechanics, underpinning the simulation of physical phenomena ranging from structural deformation and heat transfer to fluid flow and wave propagation. The finite element method (FEM) ~\cite{Hughes2000, Zienkiewicz2013} has been the dominant discretisation framework for over five decades, offering a rigorous mathematical foundation including a well-understood convergence theory and systematic error control, which has made it the method of choice in both academic research and industrial practice. A persistent bottleneck in the FEM workflow, however, is the treatment of geometry. For each new domain, a boundary-conforming mesh must be generated that respects element-quality constraints and adequately resolves the relevant geometric and solution features. Mesh generation for complex two- and three-dimensional geometries can remain a time-consuming process requiring substantial expert intervention~\cite{HOLE1988, MARCO2017,ZOU2024}. Moreover, irregular unstructured meshes make it more difficult to fully exploit high-performance numerical techniques that benefit from regular tensor-product structure, such as sum factorisation, matrix-free operator evaluation, geometric multigrid, vectorisation, and cache-efficient memory access ~\cite{Kronbichler2019,Kronbichler2019b}.

The geometry problem is further compounded in settings that require repeated solves across families of shapes, such as shape optimisation, uncertainty quantification over geometric parameters, and design-space exploration, where generating and managing a distinct mesh for each geometry instance becomes a significant computational burden. The idea of circumventing this burden by solving on a fixed reference domain via a coordinate transformation is classical. In particular, isoparametric elements map a reference element to each physical element, transfinite mapping techniques~\cite{Gordon1973, Provatidis2025} extend this to entire domains, and isogeometric analysis (IGA)~\cite{Hughes2005, Cottrell2009} uses NURBS or T-spline parameterisations to represent the geometry exactly, solving the PDE on a parameter domain that serves as the reference. The present work shares the spirit of IGA, namely, solving on a fixed reference domain with geometry encoded in the mapping, but replaces the CAD-based parameterisation with a learned neural network, removing the requirement for an explicit geometric model and enabling automatic handling of shapes defined only by images or point clouds.
 
The past decade has witnessed an intense effort to integrate machine learning (ML) and data-driven techniques into the simulation of physical systems. Several surrogate-modelling strategies have been proposed, employing, among others, convolutional neural networks~\cite{Nikolopoulos2022, XU2020}, graph neural networks~\cite{Jiang2023}, Gaussian process regression~\cite{OGREN2024,DESHPANDE2025}, and polynomial chaos expansions~\cite{Zhou2020,Sharma2024}, providing computationally inexpensive approximations of otherwise costly PDE-based simulations. In addition, physics-informed neural networks (PINNs)~\cite{Raissi2019} and their variants, such as conservative PINNs~\cite{Jagtap2020}, variational PINNs~\cite{Kharazmi2021}, eXtended PINNs~\cite{Jagtap2020b}, and interval/fuzzy PINNs~\cite{Fuhg2022} to mention a few, solve PDEs by minimising residual losses over collocation points, bypassing the mesh entirely. Neural operators, such as DeepONet~\cite{Lu2021}, Fourier Neural Operators~\cite{Li2021}, Laplace Neural Operators~\cite{Cao2024}, and Variational Physics-informed Neural Operators~\cite{Eshaghi2025}, learn mappings between function spaces and can generalise across forcing terms and boundary conditions. 

These methods have generated considerable excitement as they offer mesh-free evaluation, natural handling of high-dimensional parameter spaces, and the ability to learn from data when the governing equations are partially unknown. However, a practical assessment reveals that many ML-based simulation frameworks do not yet offer the same mature and broadly applicable error-control and convergence guarantees available for classical numerical discretisations. For PINNs, rigorous error and convergence results are available under specific assumptions, but practical a priori error control comparable to that of established numerical discretisations is generally not available, and their convergence is sensitive to training hyperparameters, collocation point placement, and the balance between PDE and boundary losses, while failure modes are difficult to diagnose~\cite{Krishnapriyan2021, Wang2021, Wang2022, Sahli2024}. Standard neural operator architectures do not necessarily preserve the variational structure of the underlying PDE, and their generalisation to geometries or parameter regimes outside the training distribution is not guaranteed~\cite{Lanthaler2022}. More broadly, when a learned surrogate replaces the underlying PDE solver, an additional approximation error is introduced, whose magnitude depends on the training data, model architecture, and generalisation regime and is not, in general, controlled by the classical discretisation-error estimates available for conventional numerical methods.
 
A more principled direction, which has gained momentum in recent years, seeks to combine ML with classical numerical methods in a way that preserves the rigorous guarantees of the latter while leveraging the representational flexibility of the former. Rather than replacing the solver, the idea is to use neural networks to assist specific components of the simulation pipeline. Some examples include learning preconditioners for iterative solvers~\cite{Luz2020, Li2023, Nikolopoulos2024}, accelerating coarse-grid corrections in multigrid~\cite{Greenfeld2019, Katrutsa2020, Huang2023}, generating adaptive mesh refinement strategies~\cite{Zhang2020}, reducing Newton--Raphson iterations in nonlinear solvers~\cite{Kadeethum2022, Papadopoulos2024}, constructing data-driven coarse spaces for domain decomposition~\cite{Heinlein2021}, or building reduced-order bases~\cite{Lee2020, Fresca2022}. In these hybrid approaches the PDE is still solved by a conventional discretisation, and the neural network augments the solver without compromising its mathematical structure. This way the solution inherits the convergence properties of the numerical method, and the ML component merely reduces the computational cost or automates a traditionally manual step. 

In the specific context of domain parameterisation, several recent works have explored learning coordinate transformations for PDE solvers. The Diffeomorphing Mapping Operator Learning (DIMON) framework~\cite{Yin2024} and Diffeomorphic Neural Operators~\cite{Zhao2025} leverage neural operators to map solutions between reference and physical domains. Neural ODEs~\cite{Chen2019} provide a natural framework for constructing diffeomorphisms via continuous-time flows, however, they require numerical ODE integration at every evaluation, which is computationally expensive and provides diffeomorphism guarantees only under regularity assumptions that are difficult to verify in practice. Normalizing flows~\cite{Rezende2015, papamakarios2021} offer an alternative class of invertible neural networks, and among these, Neural Spline Flows (NSF)~\cite{Durkan2019} provide exact invertibility, analytical log-determinants, and guaranteed monotonicity through rational-quadratic spline coupling layers. 
 
The present work follows the hybrid philosophy outlined above. We use a neural network not to approximate the PDE solution, but to learn a geometric mapping in the form of a smooth, invertible coordinate transformation from a reference domain to the physical domain, parameterised as an NSF. The monotone rational-quadratic spline coupling layers of the NSF provide an invertible mapping with a positive Jacobian determinant, avoiding the need for an explicit Jacobian-sign penalty during training. For a given learned map, the governing PDE can then be consistently transferred to the corresponding learned domain through a diffeomorphic pullback. In practice, the learned domain provides an approximation of the target physical geometry, and the resulting geometric mapping error constitutes an additional source of numerical error alongside the finite element discretisation and quadrature errors. To handle families of shapes with a single model, we adopt a two-stage training strategy inspired by the auto-decoder framework of Park et al.~\cite{Park2019}. Specifically, in the first stage, per-shape latent codes are optimised jointly with the flow to learn a smooth latent space. In the second stage, a lightweight CNN encoder is trained to predict the latent code directly from the signed distance function image, supervised by the codes obtained in the first stage. This decomposition enables instant inference on new shapes via a single forward pass. The PDE is subsequently solved on the reference square using a standard Galerkin FEM, with the mapping entering only through the geometric coefficients evaluated at the quadrature points. 

The main contributions are: (i) the use of NSFs as smooth, invertible domain maps for finite element analysis, with monotone spline transformations that maintain a positive Jacobian determinant without requiring an explicit Jacobian-sign penalty; (ii) a two-stage conditioning strategy in which an auto-decoder first learns a smooth latent space of shape deformations, and a CNN encoder is then trained to predict latent codes from signed distance function (SDF) images, enabling instant inference on unseen geometries without iterative fitting; (iii) a pullback finite element formulation on a fixed structured reference mesh, in which the learned NSF mapping enters only through the transformed geometric coefficients; and (iv) validation on benchmark problems demonstrating that the learned mappings provide sufficiently accurate geometric representations to obtain reliable finite element solutions across the considered test cases.
 
The remainder of this paper is organised as follows. Section~\ref{sec:ref_fem} establishes the mathematical framework for reference-domain FEM via diffeomorphic pullback, deriving the transformed weak form, boundary conditions, and the solution equivalence guarantee. Section~\ref{sec:nsf} describes the Neural Spline Flow architecture, the training procedure, and the auto-decoder extension for shape families. Section~\ref{sec:implementation} presents the end-to-end algorithm and discusses efficient evaluation of geometric quantities. Section~\ref{sec:experiments} provides numerical validation on several benchmark problems. Lastly, Section~\ref{sec:conclusions} summarizes the main outcomes of this work and discusses limitations, extensions and connections to broader themes in computational mechanics and scientific machine learning.

\section{Reference-Domain FEM via Diffeomorphic Pullback}
\label{sec:ref_fem}

\subsection{PDE discretization using FEM}
Let $\Omega \subset \mathbb{R}^{d_s}$, $d_s\in\{1,2,3\}$, be a bounded Lipschitz domain with boundary decomposition $\partial\Omega=\overline{\Gamma}_D\cup\overline{\Gamma}_N$ and $\Gamma_D\cap\Gamma_N=\emptyset$. We consider a boundary value problem posed on $\Omega$ for an unknown field $\boldsymbol{u}:\Omega\to\mathbb{R}^{d_f}$, where $d_f\ge 1$ denotes the number of components (e.g., $d_f=1$ for scalar diffusion, or $d_f=3$ for 3D linear elasticity). We write the strong form abstractly as

\begin{equation}
\mathcal{L}\big(\boldsymbol{u} \big) = \boldsymbol{f} \quad \text{in } \Omega,
\qquad
\boldsymbol{u} = \boldsymbol{u}_D \quad \text{on } \Gamma_D,
\qquad
\mathcal{B}_N\big(\boldsymbol{u}\big)= \boldsymbol{t}_N \quad \text{on } \Gamma_N
\label{eq:strong_abstract}
\end{equation}

\noindent where $\mathcal{L}$ is a differential operator, $\mathcal{B}_N$ denotes the natural (Neumann-type) boundary operator (e.g., flux or traction). Throughout this section we focus on stationary boundary value problems to streamline notation. Extensions to time-dependent problems (via time discretization and repeated solves) and to nonlinear problems (via linearization/Newton-type iterations in the weak form) follow standard finite element practice.

We introduce the admissible trial and test spaces
\begin{equation}
\mathcal{U}
=
\left\{
\boldsymbol{u}\in[H^1(\Omega)]^{d_f}:
\boldsymbol{u}|_{\Gamma_D}=\boldsymbol{u}_D
\right\},
\qquad
\mathcal{V}
=
\left\{
\boldsymbol{v}\in[H^1(\Omega)]^{d_f}:
\boldsymbol{v}|_{\Gamma_D}=\boldsymbol{0}
\right\}
\end{equation}
Assuming that the boundary-value problem admits a weak formulation,
we seek $\boldsymbol{u}\in\mathcal{U}$ such that
\begin{equation}
\kappa(\boldsymbol{u},\boldsymbol{v})
=
\ell(\boldsymbol{v}),
\qquad
\forall\,\boldsymbol{v}\in\mathcal{V}
\label{eq:weak_problem}
\end{equation}

\noindent where $\kappa(\cdot,\cdot)$ and $\ell(\cdot)$ encode the PDE operator and forcing/boundary data after integration by parts. For homogeneous Dirichlet data, $\mathcal{U}=\mathcal{V}$, and in the linear coercive case, the standard boundedness and coercivity conditions
\begin{equation}
|\kappa(u,v)|
\leq M\|u\|_{\mathcal V}\|v\|_{\mathcal V},
\qquad
\kappa(v,v)\geq\alpha\|v\|_{\mathcal V}^2,
\qquad \alpha>0
\end{equation}
together with boundedness of $\ell$, guarantee existence and uniqueness
by the Lax--Milgram theorem.

For later use, it is helpful to interpret $\kappa$ and $\ell$ as integrals over $\Omega$ and $\Gamma_N$, respectively. For example, a broad class of second-order problems admits the structure

\begin{equation}\label{eq:generic_integral_forms}
\kappa(\boldsymbol{u},\boldsymbol{v}) = \int_\Omega \mathcal{A}\!\left(\boldsymbol{x}\right)\,
\mathcal{D}\boldsymbol{u} : \mathcal{D}\boldsymbol{v} \,\mathrm{d}\boldsymbol{x}
\;+\;\int_\Omega \mathcal{C}\!\left(\boldsymbol{x}\right)\boldsymbol{u}\cdot \boldsymbol{v}\,\mathrm{d}\boldsymbol{x},
\qquad
\ell(\boldsymbol{v}) = \int_\Omega f\cdot \boldsymbol{v}\,\mathrm{d}\boldsymbol{x} + \int_{\Gamma_N} \boldsymbol{t}_N\cdot \boldsymbol{v}\,\mathrm{d}s
\end{equation}

\noindent where $\mathcal{D}$ denotes the relevant differential operator (e.g., $\nabla$ for diffusion, $\varepsilon(\cdot)$ for elasticity), and $\mathcal{A},\mathcal{C}$ are coefficient tensors. Equation~\eqref{eq:generic_integral_forms} is intentionally generic, since the pullback derivation presented later only requires that the weak form can be written in terms of the field and its derivatives integrated over $\Omega$ (and $\Gamma_N$).

\subsection{Diffeomorphic mapping and reference-square assumptions}
\label{subsec:mapping_assumptions}

Let $\Omega^0\subset\mathbb{R}^{d_s}$ denote a fixed reference domain. In this work we take $\Omega^0$ to be a square in 2D (e.g., $\Omega^0=(-1,1)^2$), or the unit cube in 3D, equipped with a fixed mesh that is reused across geometries. We assume the existence of a $C^1$-diffeomorphism

\begin{equation}
\begin{aligned}
\begin{array}{cccc}
G:& \Omega^0 & \longrightarrow & \Omega \\
&\boldsymbol{\xi} & \longmapsto & \boldsymbol{x}=G(\boldsymbol{\xi})
\end{array}
\qquad
\begin{array}{cccc}
F:=G^{-1}:&\Omega & \longrightarrow & \Omega^0 \\
&\boldsymbol{x} & \longmapsto & \boldsymbol{\xi}=G^{-1}(\boldsymbol{x})
\end{array}
\end{aligned}
\label{eq:Phi_G_def}
\end{equation}
\noindent such that $\det D G(\boldsymbol{\xi}) > 0$ for all $\boldsymbol{\xi}\in \Omega^0$, while $\Omega$ denotes an arbitrary domain.

We denote the Jacobian of the map by

\begin{equation}
\boldsymbol{J}(\boldsymbol{\xi})
:= D G(\boldsymbol{\xi})\in\mathbb{R}^{d_s\times d_s},
\qquad
\text{and} \quad
 \mathcal{J} = \det \boldsymbol{J}(\boldsymbol{\xi})
\label{eq:J_det_def}
\end{equation}

\noindent and we use the standard change-of-variables identities

\begin{equation}\label{eq:cov_change}
\boldsymbol{x} = G(\boldsymbol{\xi}), \qquad
\mathrm{d}\boldsymbol{x} = \mathcal{J} \,\mathrm{d}\boldsymbol{\xi}
\end{equation}

The boundary is mapped accordingly: $\Gamma^0_D := G^{-1}(\Gamma_D)$ and $\Gamma^0_N := G^{-1}(\Gamma_N)$, yielding $\partial\Omega^0 = \overline{\Gamma}^0_D \cup \overline{\Gamma}^0_N$. For Neumann terms we also require the surface measure transformation (\textit{Nanson's} relation)

\begin{equation}\label{eq:ds_transform}
\mathrm{d}s_{\boldsymbol{x}} = \mathcal{J}  \, \|\boldsymbol{J}(\boldsymbol{\xi})^{-T}\,\boldsymbol{n}_{\boldsymbol{\xi}}\|\,\mathrm{d}s_{\boldsymbol{\xi}}
\end{equation}

\noindent where $\boldsymbol{n}_{\boldsymbol{\xi}}$ is the outward unit normal on $\partial\Omega^0$.
In the present work, the map $G:= G_{\boldsymbol{\theta}}$ is realised by a Neural Spline Flow with parameters $\boldsymbol{\theta}$. However, to ease notation in the subsequent sections, the explicit dependence of $G$ on $\boldsymbol{\theta}$ will be suppressed, unless where needed to be explicitly stated.

\begin{remark}[Corners of the reference square]\label{rem:corners}
The square and a smoothly bounded domain such as an ellipse are diffeomorphic as \emph{open} sets, but no $C^1$-diffeomorphism between their closures exists, since an interior angle cannot be mapped onto a smooth boundary by a $C^1$ map. Accordingly, Eq.~\eqref{eq:Phi_G_def} is to be understood on the open sets, with $G$ extending to a homeomorphism of the closures. 

In practice this means that an exact map would have to flatten each $90^{\circ}$ corner of $\partial\Omega^0$ into a straight boundary, driving its Jacobian determinant to zero there and the geometric coefficients $\boldsymbol{A}_{geo}$ to large values. The NSF map, as presented later in section~\ref{sec:nsf}, cannot do this, since its Jacobian determinant is strictly positive by construction, so it smooths the corner instead. Thus, $\boldsymbol{A}_{geo}$ remains bounded and the discrepancy appears as a small loss of geometric fidelity in four neighbourhoods, confined to the corner elements of the reference mesh. Despite this limitation, however, the square is retained deliberately, as a smooth reference domain would remove the corner issue but forfeit the structured tensor-product mesh on which the entire computational rationale of the method rests.
\end{remark}

\subsection{Change of variables and pulled-back weak form on the square}
\label{subsec:pullback_weak}

Define the pullback of a physical-domain field $\boldsymbol{u}$ to the reference domain by

\begin{equation}\label{eq:pullback_def}
\hat{\boldsymbol{u}}(\boldsymbol{\xi}) := \boldsymbol{u}\!\left(G(\boldsymbol{\xi})\right),
\qquad
\left(\text{similarly}\quad \hat{\boldsymbol{v}}(\boldsymbol{\xi}) := \boldsymbol{v}\!\left(G(\boldsymbol{\xi})\right) \right)
\end{equation}

\noindent By the chain rule, for vector-valued fields one has

\begin{equation}\label{eq:grad_transform_vec}
\nabla_{\boldsymbol{x}} \boldsymbol{u}(\boldsymbol{x})
= \nabla_{\boldsymbol{\xi}}\hat{\boldsymbol{u}}(\boldsymbol{\xi})\,J(\boldsymbol{\xi})^{-1},
\qquad
\boldsymbol{x}=G(\boldsymbol{\xi})
\end{equation}

\noindent For scalar fields, Eq.~\eqref{eq:grad_transform_vec} reduces to $\nabla_{\boldsymbol{x}}u = J^{-T}\nabla_{\boldsymbol{\xi}}\hat{u}$.

Substituting Eqs.~\eqref{eq:cov_change}--\eqref{eq:grad_transform_vec} into the weak form produces a reference-domain problem: find $\hat{\boldsymbol{u}}\in\mathcal{U}^0$ such that

\begin{equation}\label{eq:weak_ref}
\hat{\kappa}\!\left(\hat{\boldsymbol{u}},\hat{\boldsymbol{v}}\right)
=
\hat{\ell}\!\left(\hat{\boldsymbol{v}}\right)
\qquad \forall \hat{\boldsymbol{v}}\in\mathcal{V}^0
\end{equation}

\noindent where $\mathcal{U}^0:=\big[H^1_{\Gamma^0_D}(\Omega^0)\big]^{d_f}$ and $\mathcal{V}^0$ is defined analogously.

To make the dependence on the geometric mapping explicit, consider the prototypical diffusion-like term $\int_\Omega (\nabla u)^T \boldsymbol{A}(\boldsymbol{x})\,\nabla v\,\mathrm{d}\boldsymbol{x}$. Using Eq.~\eqref{eq:grad_transform_vec} and $\mathrm{d}\boldsymbol{x}=\mathcal{J}\,\mathrm{d}\boldsymbol{\xi}$ yields

\begin{equation}\label{eq:diffusion_pullback}
\int_\Omega (\nabla_{\boldsymbol{x}} u)^T \boldsymbol{A}(\boldsymbol{x})\,\nabla_{\boldsymbol{x}} v\,\mathrm{d}\boldsymbol{x}
=
\int_{\Omega^0}
(\nabla_{\boldsymbol{\xi}}\hat{u})^T
\underbrace{\Big(\mathcal{J}\,\boldsymbol{J}^{-1}\,\boldsymbol{A}\!\big(G(\boldsymbol{\xi})\big)\,\boldsymbol{J}^{-T}\Big)}_{=:~\boldsymbol{A}_{\mathrm{geo}}(\boldsymbol{\xi})}
\nabla_{\boldsymbol{\xi}}\hat{v}\,\mathrm{d}\boldsymbol{\xi}
\end{equation}

\noindent Thus, even when the physical coefficients $\boldsymbol{A}$ are constant, the pulled-back operator contains the \emph{geometric coefficients} $\boldsymbol{A}_{\mathrm{geo}}(\boldsymbol{\xi})$ induced by $\boldsymbol{J}$ and $\mathcal{J}$.
Analogous expressions arise for vector problems. For instance, in linear elasticity one replaces $\nabla(\cdot)$ by the symmetric gradient $\varepsilon(\cdot)$ and obtains a pulled-back bilinear form whose integrand involves $\boldsymbol{J}^{-1}$ and $\mathcal{J}$ through the transformed strain operator.

Body force and reaction terms transform similarly:

\begin{equation}\label{eq:rhs_pullback}
\int_\Omega \boldsymbol{f}(\boldsymbol{x})\cdot \boldsymbol{v}(\boldsymbol{x})\,\mathrm{d}\boldsymbol{x}
=
\int_{\Omega^0}
\hat{\boldsymbol{f}}(\boldsymbol{\xi})\cdot \hat{\boldsymbol{v}}(\boldsymbol{\xi})\,
\mathcal{J}(\boldsymbol{\xi})\,\mathrm{d}\boldsymbol{\xi},
\qquad \text{where }
\hat{\boldsymbol{f}}(\boldsymbol{\xi}) := \boldsymbol{f}\!\left(G(\boldsymbol{\xi})\right)
\end{equation}

Regarding the transformation of boundary conditions, for Dirichlet data $\boldsymbol{u}=\boldsymbol{u}_D$ on $\Gamma_D$, the pulled-back condition reads

\begin{equation}\label{eq:dirichlet_pullback}
\hat{\boldsymbol{u}}(\boldsymbol{\xi}) = \hat{\boldsymbol{u}}_D(\boldsymbol{\xi})
:= \boldsymbol{u}_D\!\left(G(\boldsymbol{\xi})\right)
\qquad \text{on } \Gamma^0_D
\end{equation}

\noindent In the common homogeneous case $\boldsymbol{u}_D\equiv 0$, one simply enforces $\hat{\boldsymbol{u}}=0$ on $\Gamma^0_D$.

Neumann data enters the weak form trough boundary integrals.
Using the surface transformation of Eq. \eqref{eq:ds_transform}, one obtains

\begin{equation}\label{eq:neumann_pullback}
\int_{\Gamma_N} \boldsymbol{t}_N(\boldsymbol{x})\cdot \boldsymbol{v}(\boldsymbol{x})\,\mathrm{d}s_{\boldsymbol{x}}
=
\int_{\Gamma_N^0}
\underbrace{
\boldsymbol{t}_N\!\left(G(\boldsymbol{\xi})\right)\,
\mathcal{J}\,\|J(\boldsymbol{\xi})^{-T}\,\boldsymbol{n}_{\boldsymbol{\xi}}\|
}_{=:~\hat{\boldsymbol{t}}_N(\boldsymbol{\xi})}
\cdot \hat{\boldsymbol{v}}(\boldsymbol{\xi})\,\mathrm{d}s_{\boldsymbol{\xi}}
\end{equation}

\subsection{Discrete solution transfer back to the physical domain}
\label{subsec:map_back}

Let $\mathcal{U}^0_h\subset\mathcal{U}^0$ and $\mathcal{V}^0_h\subset\mathcal{V}^0$ be conforming finite-dimensional spaces on the fixed mesh $\mathcal{T}^0_h$ on the reference domain, with bases $\{\hat{\boldsymbol{\varphi}}_j\}_{j=1}^{d}$ and $\{\hat{\boldsymbol{\psi}}_i\}_{i=1}^{d}$.
The discrete reference problem reads: find

\begin{equation}\label{eq:discrete_ref}
\hat{\boldsymbol{u}}_h(\boldsymbol{\xi})=\sum_{j=1}^{d}\hat{\boldsymbol{u}}_j\,\hat{\boldsymbol{\varphi}}_j(\boldsymbol{\xi}) \in \mathcal{U}^0_h
\quad \text{such that} \quad
\sum_{j=1}^{d}\hat{\kappa}(\hat{\boldsymbol{\varphi}}_j,\hat{\boldsymbol{\psi}}_i)\,\hat{\boldsymbol{u}}_j
=
\hat{\ell}(\hat{\boldsymbol{\psi}}_i),\;\; i=1,\ldots,d
\end{equation}
leading to the $d\times d$ linear system
\begin{equation}\label{eq:lin_system_ref}
\hat{\mathbf{K}}\,\hat{\mathbf{u}}=\hat{\mathbf{f}}
\end{equation}

Once $\hat{\boldsymbol{u}}_h$ is computed on the reference square, the corresponding approximation on the physical domain is obtained by composition:

\begin{equation}\label{eq:pushforward_solution}
\boldsymbol{u}_h(\boldsymbol{x})
:= \hat{\boldsymbol{u}}_h\!\left(G^{-1}(\boldsymbol{x})\right),
\qquad \boldsymbol{x}\in\Omega
\end{equation}

In particular, evaluating $\boldsymbol{u}_h$ at a physical point $\boldsymbol{x}$ requires (i) computing $\boldsymbol{\xi}=G^{-1}(\boldsymbol{x})$, (ii) locating the element in $\mathcal{T}^0_h$ containing $\boldsymbol{\xi}$, and (iii) evaluating the reference basis there. As mentioned, in this work the reference domain is chosen to be a square since it allows for a perfectly structured mesh with square 4-node elements. As a consequence, step~(ii) reduces to integer arithmetic on regular grids, avoiding the search cost associated with unstructured meshes. Derived quantities (gradients, fluxes, stresses) are recovered via the chain rule. For instance, for vector fields

\begin{equation}\label{eq:grad_mapback}
\nabla_{\boldsymbol{x}} \boldsymbol{u}_h(\boldsymbol{x})
=
\nabla_{\boldsymbol{\xi}}\hat{\boldsymbol{u}}_h(\boldsymbol{\xi})\,
\boldsymbol{J}(\boldsymbol{\xi})^{-1},
\qquad \boldsymbol{\xi}=G^{-1}(\boldsymbol{x})
\end{equation}

\noindent Eqs.~\eqref{eq:pushforward_solution}--\eqref{eq:grad_mapback} complete the transfer of the numerical solution from the reference square back to the original domain.

\subsection{Solution equivalence under diffeomorphisms}
\label{subsec:solution_equivalence}

The preceding pullback derivation is based on a standard change of variables. For completeness, we state explicitly the resulting equivalence between the physical- and reference-domain problems.

\begin{proposition}[Solution equivalence]\label{prop:equivalence}
Let $G:\Omega^0\to\Omega$ be a $C^1$-diffeomorphism with $\mathcal{J}=\det DG>0$ everywhere, and assume that $G$ and $G^{-1}$ are sufficiently regular for the pullback to define an isomorphism between the corresponding $H^1$ spaces. Let the bilinear and linear forms $\kappa,\ell$ satisfy the hypotheses of Lax--Milgram on $\Omega$. Then the pulled-back problem on $\Omega^0$ satisfies the corresponding well-posedness conditions, with constants depending on the mapping $G$, and:
\begin{enumerate}[label=(\roman*),leftmargin=2em]
    \item $\boldsymbol{u}\in\mathcal{U}$ solves the physical-domain problem in Eq.~\eqref{eq:weak_problem} if and only if $\hat{\boldsymbol{u}}=\boldsymbol{u}\circ G\in\mathcal{U}^0$ solves the reference-domain problem in Eq.~\eqref{eq:weak_ref};
    
    \item the correspondence is exact, so that the change of variables itself introduces no modelling error;
    
    \item for a given exact mapping $G$, the remaining numerical errors arise from the finite element discretisation on $\Omega^0$ and from the numerical quadrature used to evaluate the transformed coefficients.
\end{enumerate}
\end{proposition}

\noindent
Statement~(i) follows from the bijectivity of the pullback operator $\boldsymbol{u}\mapsto\boldsymbol{u}\circ G$ between the corresponding function spaces, together with the change-of-variables identities established in Eqs.~\eqref{eq:diffusion_pullback}--\eqref{eq:rhs_pullback}. 
Statement~(ii) emphasizes that, for an exact diffeomorphism, the pullback is not an approximation or model reduction of the governing PDE, but an equivalent reformulation on a different domain.

For a fixed exact map $G$, the numerical approximation therefore introduces the usual finite element discretisation error, together with a possible quadrature error. The latter deserves some attention because the transformed coefficients contain quantities such as $\boldsymbol{J}^{-1}$ and $\mathcal{J}$, which are generally non-polynomial functions of the reference coordinates. Consequently, standard Gauss quadrature is not in general exact for the pulled-back integrands. When required, this error can be reduced by using higher-order quadrature (over-integration), so that it remains small relative to the finite element discretisation error.

An important distinction arises when the diffeomorphism itself is learned, as is the case in this work. In the proposed framework, the Neural Spline Flow provides a map $G_{\boldsymbol{\theta}}:\Omega^0\rightarrow\Omega_{\boldsymbol{\theta}}$, where
\begin{equation}
    \Omega_{\boldsymbol{\theta}}:= G_{\boldsymbol{\theta}}(\Omega^0)
\end{equation}
is the domain represented by the learned transformation. For any fixed $G_\theta$ satisfying the diffeomorphism assumptions above, Proposition~\ref{prop:equivalence} applies exactly between $\Omega^0$ and $\Omega_{\boldsymbol{\theta}}$. In practice, however, $G_{\boldsymbol{\theta}}$ is trained to approximate a target physical domain $\Omega$, and generally $\Omega_{\boldsymbol{\theta}}\neq\Omega$ exactly. The discrepancy between the learned and target geometries therefore introduces an additional \emph{geometric mapping error}, which is distinct from both the finite element discretisation error and the quadrature error. As will be demonstrated in the numerical examples, this error propagates through the transformed geometric coefficients and may ultimately pose a limit on the accuracy of the physical-domain solution.

\section{Learning Diffeomorphisms between reference and physical domain with Neural Spline Flows}
\label{sec:nsf}

This section describes a strategy to construct a smooth, invertible mapping between the fixed reference square $\Omega^0=(-1,1)^2$ and a physical domain $\Omega$. Although the NSF architecture itself naturally extends to higher dimensions, the present work focuses exclusively on two-dimensional domains, while extensions of the complete geometry-generation and FEM pipeline to 3D are left for future work.

The downstream objective is the reference-domain FEM formulation of Section~\ref{sec:ref_fem}, which requires evaluating the map $G:\Omega^0\to\Omega$, its Jacobian $\boldsymbol{J}=DG$, and the determinant $\mathcal{J}=\det \boldsymbol{J}$ at a large number of quadrature points. We therefore seek a mapping representation that (i) guarantees $\mathcal{J}>0$ by construction, (ii) provides the full Jacobian matrix via automatic differentiation, (iii) generalises across a family of shapes, and (iv) is amenable to efficient batch evaluation.

We parameterise the map using a \emph{Neural Spline Flow} (NSF)~\cite{Durkan2019}, an invertible neural network architecture based on monotone rational-quadratic spline coupling layers. The NSF forward pass defines $G_{\boldsymbol{\theta}}:\Omega^0\to\Omega$ directly, where the subscript $\boldsymbol{\theta}$ denotes the network's trainable parameters. Specifically, the network takes a reference-square point as input and returns the corresponding physical-domain point.

To handle multiple shapes with a single model, we adopt an \emph{auto-decoder} framework~\cite{Park2019}: each shape is associated with a low-dimensional latent code $\mathbf{z}\in\mathbb{R}^{d_z}$ that conditions the spline parameters, producing a shape-dependent diffeomorphism $G_{\boldsymbol{\theta}}(\cdot\,;\mathbf{z}):\Omega^0\to\Omega_{\mathbf{z}}$.

\subsection{The Neural Spline Flow map}
\label{subsec:nsf_architecture}

A Neural Spline Flow represents the map $G_{\boldsymbol{\theta}}$ as a composition of $L$ invertible
coupling layers, $G_{\boldsymbol{\theta}} = h^{(L)}\circ\cdots\circ h^{(1)}$, acting on intermediate
coordinate vectors $\boldsymbol{\xi}^{(\ell)}\in\mathbb{R}^{d}$. Following the mapping direction
adopted above, the endpoints of this composition are identified with the reference and
physical coordinates,
\begin{equation}\label{eq:flow_endpoints}
\boldsymbol{\xi}^{(0)} = \boldsymbol{\xi}\ \ (\text{reference input}), \qquad
\boldsymbol{\xi}^{(L)} = \mathbf{x}\ \ (\text{physical output})
\end{equation}
while the intermediate states $\boldsymbol{\xi}^{(\ell)}$, $0<\ell<L$, interpolate between the two
and carry no direct geometric meaning. In each layer, the coordinates are split into a
\emph{passive} group and an \emph{active} group. The passive coordinates are left unchanged
and serve as input to a neural network that predicts the parameters of a monotone
transformation applied to the active coordinates. Successive layers alternate which
coordinates are active, ensuring that all dimensions are transformed.

Concretely, for $d=2$, consider a coupling layer that transforms the active coordinate $\xi_2$
while keeping $\xi_1$ fixed (here $\xi_1,\xi_2$ denote the components of the intermediate state
$\boldsymbol{\xi}^{(\ell)}$). The layer computes
\begin{equation}\label{eq:coupling_layer}
\xi_1' = \xi_1, \qquad \xi_2' = s\!\left(\xi_2\,;\,\mathrm{MLP}(\xi_1,\mathbf{z})\right),
\end{equation}
where $s(\cdot\,;\mathrm{MLP})$ is a monotone scalar function parameterised by
$\mathrm{MLP}$, and $\mathrm{MLP}(\xi_1,\mathbf{z})$ is the output of a small Multilayer Perceptron (MLP) that takes the passive coordinate $\xi_1$ and the conditioning vector $\mathbf{z}$ as input.
The next layer swaps roles: $\xi_1'$ becomes active and $\xi_2'$ becomes passive.

The key innovation of~\cite{Durkan2019} is to choose $s$ as a \emph{rational-quadratic
spline}, where the domain is partitioned into $K$ bins by knot points, and within each bin $s$ is
defined by a rational-quadratic segment
\begin{equation}\label{eq:rq_spline}
s(\zeta;\mathrm{MLP})
= \frac{\alpha(\zeta)\,\zeta'^{(k)} + \beta(\zeta)\,\zeta'^{(k+1)}}{\gamma(\zeta)}
\end{equation}
where $\zeta\in[\zeta^{(k)}, \zeta^{(k+1)}]$ is the input within the $k$-th bin,
$\{\zeta^{(k)}, \zeta'^{(k)}\}_{k=0}^{K}$ are the knot positions (input) and heights
(output), and $\alpha,\beta,\gamma$ are quadratic polynomials in $\zeta$ determined by the
knots and the positive derivatives at each interior knot. The precise
functional form, given in \cite{Durkan2019}, ensures $C^1$ continuity at knot boundaries, while
outside the spline domain, the transformation defaults to the identity.

\begin{figure}[H]
  \centering
  \includegraphics[width=\textwidth]{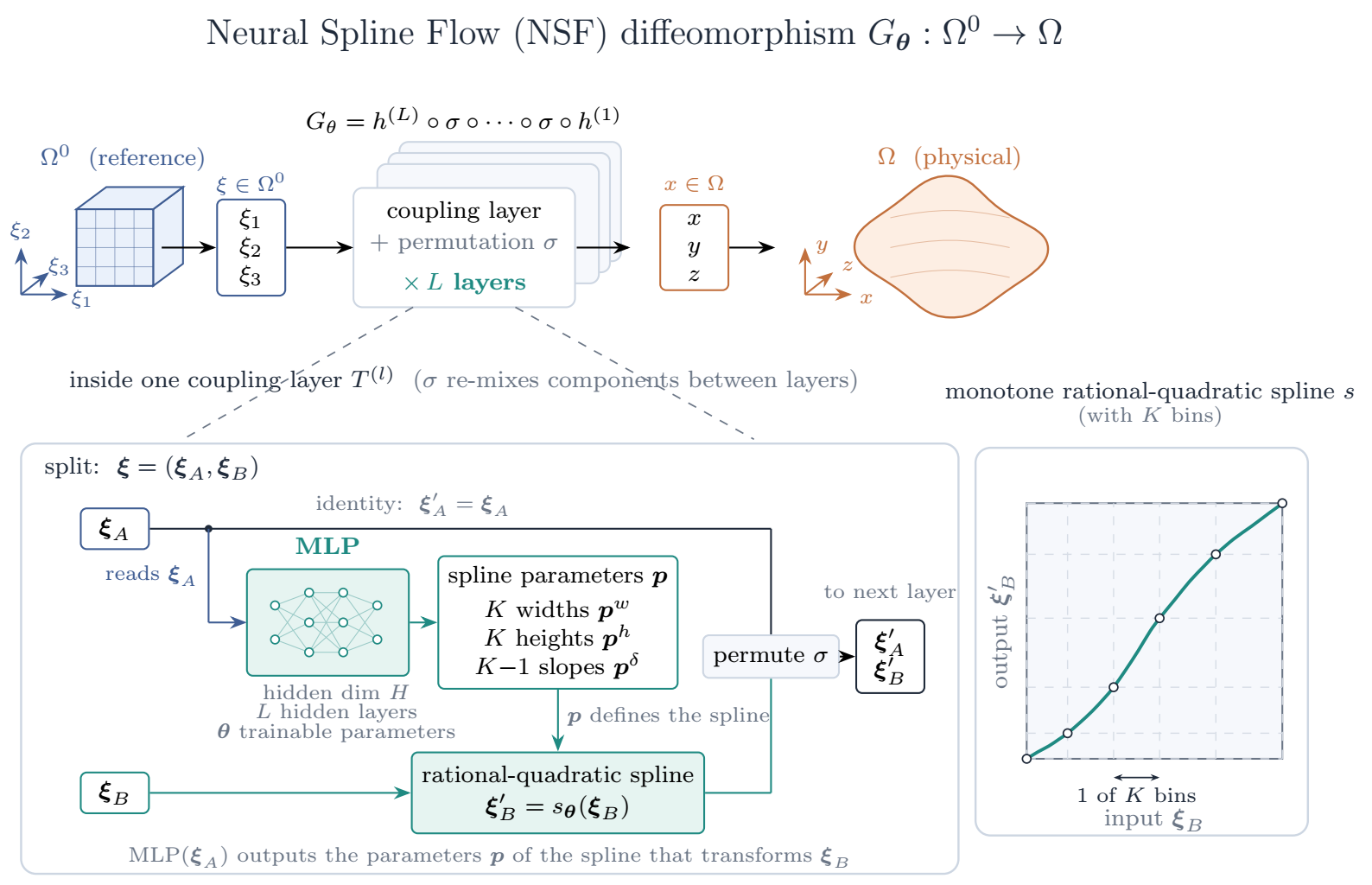}
  \caption{\textbf{NSF diffeomorphism architecture.} $G_{\boldsymbol\theta}:\Omega^0\!\to\!\Omega$
is built from $L$ rational-quadratic coupling layers $h^{(l)}$ and permutations
$\sigma$. Each layer splits $\boldsymbol\xi=(\boldsymbol\xi_A,\boldsymbol\xi_B)$, leaves $\boldsymbol\xi_A$ fixed, and uses $\mathrm{MLP}(\boldsymbol\xi_A)$ to predict the parameters $\boldsymbol p$ of a monotone spline that maps $\boldsymbol\xi_B\mapsto\boldsymbol\xi_B'=s_{\boldsymbol\theta}(\boldsymbol\xi_B)$; $\sigma$ then permutes the channels. Spline monotonicity guarantees invertibility and a strictly
positive Jacobian.}
  \label{fig:NSF_architecture}
\end{figure}

Based on the above, the MLP in each coupling layer outputs $3K-1$ parameters $\boldsymbol{p}\in\mathbb{R}^{3K-1}$, which consist of the $K$ bin widths $\boldsymbol{p}^w\in\mathbb{R}^{K}$ (passed through a softmax to ensure positivity and summation to the domain width), the $K$ bin heights $\boldsymbol{p}^h\in\mathbb{R}^{K}$ (similarly
constrained), and the $K-1$ derivatives $\boldsymbol{p}^{\delta}\in\mathbb{R}^{K-1}$ at the interior nodes (passed through a softplus to ensure positivity). The positivity of the knot derivatives is the mechanism that guarantees monotonicity of $s$. A schematic representation of the proposed NSF architecture and the associated operations is shown in Fig.~\ref{fig:NSF_architecture}.

The coupling-layer structure introduced above yields three properties of
$G_{\boldsymbol{\theta}}$ that are essential for the reference-domain FEM: (i) a guaranteed
positive Jacobian determinant, (ii) an analytic inverse, and (iii) efficient access to the full
Jacobian matrix.

Regarding the first property, the Jacobian of the coupling layer in
Eq.~\eqref{eq:coupling_layer} is lower triangular. For example, in a
2D domain with $\xi_1$ denoting the passive coordinate and $\xi_2$ the
active coordinate, the Jacobian reads
\begin{equation}\label{eq:coupling_jacobian}
\frac{\partial(\xi_1',\xi_2')}{\partial(\xi_1,\xi_2)}
=
\begin{pmatrix}
1 & 0 \\[3pt]
\displaystyle\frac{\partial s}{\partial \xi_1}
&
\displaystyle\frac{\partial s}{\partial \xi_2}
\end{pmatrix}
\end{equation}

\noindent Its determinant is therefore simply $\partial s/\partial \xi_2$, i.e., the derivative of the monotone spline, which is strictly positive by construction.

The complete NSF additionally contains fixed coordinate permutations between successive coupling layers. Denoting the corresponding permutation matrices by $P^{(\ell)}$, the determinant of the full
Jacobian can be written as
\begin{equation}
\det DG_{\boldsymbol{\theta}}(\boldsymbol{\xi})
=
\left(\prod_{\ell} \det P^{(\ell)}\right)
\prod_{\ell=1}^{L}
\frac{\partial s^{(\ell)}}
{\partial \xi_{\mathrm{active}}^{(\ell-1)}}
\end{equation}
In the architecture employed here, the sequence of permutations is chosen such that its total composition has positive parity, namely
\begin{equation}
\prod_{\ell} \det P^{(\ell)} = +1
\end{equation}
Since the derivative of every spline transformation is strictly positive,
it follows that
\[
\det DG_{\boldsymbol{\theta}}(\boldsymbol{\xi}) > 0
\qquad
\forall\,\boldsymbol{\xi}\in\Omega^0
\]
Consequently, taking the logarithm of the determinant and applying the
chain rule gives
\begin{equation}\label{eq:full_logdet}
\log \det DG_{\boldsymbol{\theta}}(\boldsymbol{\xi})
=
\sum_{\ell=1}^{L}
\log
\left(
\frac{\partial s^{(\ell)}}
{\partial \xi_{\mathrm{active}}^{(\ell-1)}}
\right)
\end{equation}
Thus, the full Jacobian determinant
$\mathcal{J}=\det DG_{\boldsymbol{\theta}}$ is guaranteed to remain strictly positive for any choice of network weights, without requiring Jacobian-sign regularisation or projection. Together with the exact invertibility of the coupling transformations, this ensures that the NSF map is invertible and orientation preserving by construction, so that the diffeomorphism assumptions employed in
Proposition~\ref{prop:equivalence} are satisfied by the architecture rather than enforced through optimisation.

In addition, the same structure makes each coupling layer analytically invertible: given $(\xi_1', \xi_2')$, one recovers $\xi_1 = \xi_1'$ immediately and then solves the monotone scalar equation
$s(\xi_2; \mathrm{MLP}(\xi_1, \mathbf{z})) = \xi_2'$ for $\xi_2$. Since $s$ is piecewise
rational-quadratic, this reduces to locating the correct bin and solving a quadratic, so no
iterative root-finding is needed. The full inverse $G_{\boldsymbol{\theta}}^{-1}$ follows by inverting the layers in reverse order.

Finally, the log-determinant in Eq.~\eqref{eq:full_logdet} is available analytically as a
byproduct of the forward pass, since it accumulates the univariate spline derivatives of the
individual layers. The FEM pullback, however, also requires the full $d_s \times d_s$
Jacobian matrix $\boldsymbol{J} = DG_{\boldsymbol{\theta}}(\boldsymbol{\xi})$, including its
off-diagonal entries, to form the geometric coefficients of the transformed weak form. This
is obtained via reverse-mode automatic differentiation through the forward pass. For example, for
$d_s = 2$ it requires two backward passes, one per output component. These evaluations are
batched over all quadrature points simultaneously, as described in
Section~\ref{sec:implementation}.

\begin{remark}[Choice of mapping direction]\label{rem:direction_choice}
The direction of the NSF forward pass is a design choice dictated by where the finite element method evaluates the map. From the density-estimation viewpoint the natural option is $F_{\boldsymbol{\theta}}: \Omega \to
\Omega^0$ (physical $\to$ reference), with the FEM pipeline relying on the inverse
$G_{\boldsymbol{\theta}} = F_{\boldsymbol{\theta}}^{-1}$. Both directions are learned from the same paired data and are exact bijections by construction, so the choice cannot rest on invertibility. The relevant distinction is which points the network is asked to evaluate at, and whether training supervises them. The FEM evaluates the map at the nodes and quadrature points of the fixed reference mesh, that is, at a prescribed set of points in $\Omega^0$, and it does so identically for every geometry. Training in the direction $G_{\boldsymbol{\theta}}: \Omega^0 \to
\Omega$ supervises the network on inputs drawn from $\Omega^0$ itself, so the evaluation set lies within the region covered by the training inputs, and the boundary points of $\partial\Omega^0$ at which Dirichlet conditions are imposed are supervised directly. In the opposite direction the training inputs are points of $\Omega$, while the FEM would require $F_{\boldsymbol{\theta}}^{-1}$ at prescribed points of $\Omega^0$. Those pre-images are not controlled by the training data and may fall outside the region of $\Omega$ where the flow was supervised, so the accuracy of the map there would be a matter of extrapolation rather than of fit. We therefore adopt
\begin{equation}\label{eq:forward_direction}
G_{\boldsymbol{\theta}}: \Omega^0 \to \Omega, \qquad
\boldsymbol{\xi} \mapsto \boldsymbol{x} = G_{\boldsymbol{\theta}}(\boldsymbol{\xi})
\end{equation}
so that the FEM queries only the forward pass and its Jacobian, while no inverse is required during
assembly or solution, and training supervises exactly the points at which the map is later evaluated.
\end{remark}

\subsection{Training data generation using harmonic maps}
\label{subsec:harmonic_maps}

Training the NSF requires paired data $\{(\boldsymbol{\xi}_k, \boldsymbol{x}_k)\}_{k=1}^{N_{\mathrm{train}}}$, where $\boldsymbol{x}_k \in \Omega$ is a physical point, $\boldsymbol{\xi}_k \in \Omega^0$ is its pre-image in the reference domain under a known diffeomorphism $G$, and $N_{\mathrm{train}}$ is the number of training points. A variety of established parameterization and registration techniques can be used to construct such point correspondences, including harmonic mappings \cite{Remacle2010}, optimal-transport-based methods \cite{Bonneel2023}, Teichm\"{u}ller mappings \cite{Lui2014,NIAN2016}, Ricci-flow-based parametrizations \cite{Zeng2013}, and mappings based on large elastic or hyperelastic deformations \cite{SHAMANSKIY2020}. In the present work, harmonic mappings are employed for this purpose, as described below.

The harmonic map $\Phi:\Omega \to \Omega^0$ assigns to each physical point $\boldsymbol{x}\in\Omega$
its coordinates in the reference square. Writing its two components (for the 2D case) as
\begin{equation*}
\Phi(\boldsymbol{x}) = \big(\mu(\boldsymbol{x}),\, \nu(\boldsymbol{x})\big) = \boldsymbol{\xi}
\end{equation*}
the scalar fields $\mu,\nu:\Omega\to\mathbb{R}$ return, respectively, the first and second
reference coordinates $\xi_1,\xi_2$ of the image of $\boldsymbol{x}$. Each component is taken to be
harmonic, and the geometry enters solely through the boundary data, thus, for a simply connected
domain $\Omega\subset\mathbb{R}^2$ with smooth boundary, the problem of obtaining $\Phi(\boldsymbol{x}) = \big(\mu(\boldsymbol{x}),\, \nu(\boldsymbol{x})\big)$ can be written as 

\begin{equation}\label{eq:harmonic_map}
\Delta \mu = 0,\quad \Delta \nu = 0 \ \ \text{in } \Omega,
\qquad
\mu = b_1,\quad \nu = b_2 \ \ \text{on } \partial\Omega
\end{equation}
where $\boldsymbol{b}=(b_1,b_2):\partial\Omega\to\partial\Omega^0$ is a prescribed boundary
correspondence. We construct $\boldsymbol{b}$ by parameterising $\partial\Omega$ by arclength and
mapping it monotonically onto the perimeter of $\Omega^0=[-1,1]^2$, with the four corners of the
square pinned to four boundary points of $\Omega$. The two Laplace problems
in Eq.~\eqref{eq:harmonic_map} are therefore decoupled Dirichlet problems that share only this
boundary correspondence, and together they define the map $\Phi=(\mu,\nu)$.

Equivalently, $\Phi$ is the minimiser of the Dirichlet energy
\begin{equation}\label{eq:dirichlet_energy}
E[\Phi] = \int_\Omega \big(\|\nabla \mu\|^2 + \|\nabla \nu\|^2\big)\, \mathrm{d}\boldsymbol{x}
\end{equation}
among all maps satisfying the prescribed boundary correspondence. Because the target $\Omega^0$ is convex, the Rad\'{o}--Kneser--Choquet theorem~\cite{Duren2004} guarantees that the harmonic extension of the monotone boundary homeomorphism $\boldsymbol{b}$ is univalent, and by Lewy's theorem \cite{Duren2004} its Jacobian does not vanish. Hence $\Phi$ is a diffeomorphism of the open sets for any Jordan domain $\Omega$, regardless of its own convexity. In particular, only the target domain $\Omega^0$ is required to be convex.

As noted in Remark~\ref{rem:corners}, smoothness cannot extend to the closures, since the four corners of $\partial\Omega^0$ are preimages of points at which $D\Phi^{-1}$ degenerates. This does not pose a problem in the proposed formulation, however, a practical issue is that strongly non-convex domains produce harmonic maps whose Jacobian $\det D\Phi$ becomes nearly degenerate near concavities. Consequently, the inverse $\Phi^{-1}$, which is the map the NSF must ultimately represent, develops large gradients and is poorly conditioned. Such cases would have an impact on the accuracy of the method, and even though these are not examined in this work, possible remedies are discussed as future work in the Conclusions section.

Closing this subsection, each Laplace problem in Eq.~\eqref{eq:harmonic_map} is numerically solved by finite differences on a regular grid covering a bounding box of $\Omega$, with the Dirichlet data $\boldsymbol{b}$ imposed on $\partial\Omega$. $\Phi$ is then evaluated at arbitrary points by interpolation from the grid. To build the training set, we sample interior and boundary physical points $\boldsymbol{x}_k$ (with counts $N_{\mathrm{int}}$ and $N_b$, so that $N_{\mathrm{train}}=N_{\mathrm{int}}+N_b$; see Section~\ref{subsec:single_shape}) and push them through the harmonic map to obtain $\boldsymbol{\xi}_k = \Phi(\boldsymbol{x}_k)$. The resulting pairs $\{(\boldsymbol{\xi}_k,
\boldsymbol{x}_k)\}$ are used to train $G_{\boldsymbol{\theta}}$ in the direction $\Omega^0\rightarrow \Omega$, so that $G_{\boldsymbol{\theta}}(\boldsymbol{\xi}_k)\approx\boldsymbol{x}_k$, thus, the network learns to approximate the inverse harmonic map $\Phi^{-1}:\Omega^0\to\Omega$.

\subsection{Input representation via Signed Distance Functions}
\label{subsec:input_rep}

\begin{figure}[H]
  \centering
  \includegraphics[width=\textwidth]{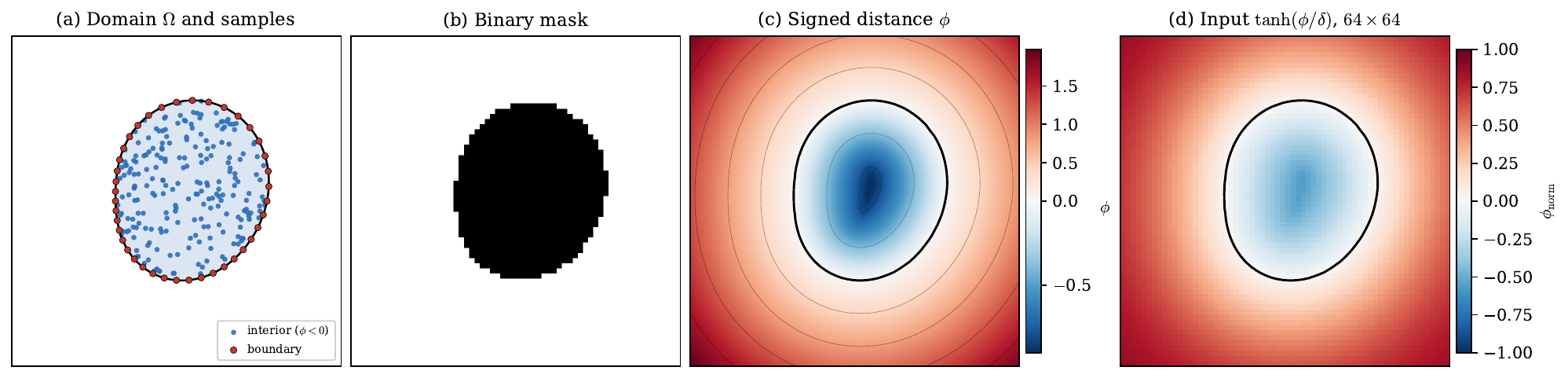}
  \caption{Signed-distance representation of a sample domain: (a) sampled interior
  ($\phi<0$) and boundary points used for training; (b) binary occupancy mask; (c) signed
  distance field $\phi$ (negative inside, positive outside) with the zero level set
  $\{\phi=0\}=\partial\Omega$ highlighted and several distance contours; (d) the normalised
  encoder input $I=\tanh(\phi/\delta)$ on the $64\times64$ grid.}
  \label{fig:sdf_representation}
\end{figure}

Both the point-sampling and the shape-conditioning stages of the pipeline require a numerical
representation of each domain $\Omega$. We use its \emph{signed distance function} (SDF)
$\phi:\mathbb{R}^2 \to \mathbb{R}$, which returns the distance from a point to $\partial\Omega$
carrying a sign: $\phi(\boldsymbol{x}) < 0$ inside $\Omega$, $\phi(\boldsymbol{x}) = 0$ on
$\partial\Omega$, and $\phi(\boldsymbol{x}) > 0$ outside (Figure~\ref{fig:sdf_representation}(c)).
The SDF serves two roles: (i) its sign provides a fast inside/outside test for sampling interior
points, and (ii) its values on a regular grid form the image descriptor that will later condition the encoder of the multi-shape model (Section~\ref{subsec:two_stage}).

Compared with a binary occupancy mask, which records only inside/outside membership, the SDF
carries richer geometric information. The mask is piecewise constant, so its gradient vanishes
everywhere except at a one-pixel jump on $\partial\Omega$ and conveys nothing about distance to
the boundary, curvature, or scale that are precisely the features the encoder must infer. The SDF,
by contrast, varies smoothly and its magnitude encodes distance to $\partial\Omega$, providing
an informative gradient throughout the domain, from which a convolutional encoder can extract
such cues readily. The difference is visible between the mask and the field in
Figure~\ref{fig:sdf_representation}(b,c). For network input we compress the unbounded SDF to
$[-1,1]$ with a smooth squashing function,
\begin{equation}\label{eq:sdf_normalize}
\phi_{\mathrm{norm}}(\boldsymbol{x}) = \tanh\!\left(\frac{\phi(\boldsymbol{x})}{\delta}\right)
\end{equation}
where $\delta > 0$ sets the length scale of the transition. With this approach the map is near-linear for $|\phi| \lesssim \delta$ (near the boundary, where geometric detail matters most) and
saturates far from it, where the precise distance is less informative. We take $\delta$ as the
$95$th percentile of $|\phi|$ over the grid, a robust scale that is insensitive to the few extreme values at the deep interior and at the corners of the bounding box. The normalised field is sampled on a regular $H \times W$ grid (e.g.\ $64 \times 64$) to produce the single-channel image $I \in \mathbb{R}^{H \times W \times 1}$ fed to the encoder (Figure~\ref{fig:sdf_representation}(d)). Training points are drawn from the same representation (Figure~\ref{fig:sdf_representation}(a)).
Specifically, boundary points are obtained by sampling the curve parameter $t \in [0, 2\pi)$ of
$\partial\Omega$ uniformly and evaluating its parameterisation, while interior points are
generated by rejection sampling within a bounding box of $\Omega$, retaining those with
$\phi(\boldsymbol{x}) < 0$.

\subsection{Single-shape training}
\label{subsec:single_shape}

Before presenting the multi-shape framework, we describe training for a single fixed domain, which serves as a building block for the general case.

For a given shape $\Omega$, the (unconditional) NSF $G_{\boldsymbol{\theta}}:\Omega^0 \to \Omega$ is trained by minimising a supervised loss on the paired data from Section~\ref{subsec:harmonic_maps}. The training objective is
\begin{equation}\label{eq:single_shape_loss}
\mathcal{L} = \mathcal{L}_{\mathrm{int}} + \lambda_\partial \, \mathcal{L}_\partial + \lambda_{\log\det} \, \mathcal{L}_{\log\det}
\end{equation}
where the three terms are:

\begin{itemize}
\item \textit{Interior loss}: Mean squared error between predicted and ground-truth physical coordinates:
\begin{equation}\label{eq:loss_interior}
\mathcal{L}_{\mathrm{int}} = \frac{1}{N_{\mathrm{int}}} \sum_{k=1}^{N_{\mathrm{int}}} \left\| G_{\boldsymbol{\theta}}(\boldsymbol{\xi}_k) - \boldsymbol{x}_k \right\|^2
\end{equation}
\item \textit{Boundary loss}: Mean squared error on boundary correspondences, weighted more heavily to enforce geometric fidelity:
\begin{equation}\label{eq:loss_boundary}
\mathcal{L}_\partial = \frac{1}{N_b} \sum_{k=1}^{N_b} \left\| G_{\boldsymbol{\theta}}(\boldsymbol{\xi}_k^b) - \boldsymbol{x}_k^b \right\|^2
\end{equation}
\item \textit{Log-determinant regularisation}: A mild penalty on the magnitude of the log-determinant to prevent extreme Jacobian values:
\begin{equation}\label{eq:loss_logdet}
\mathcal{L}_{\log\det} = \frac{1}{N_{\mathrm{int}}} \sum_{k=1}^{N_{\mathrm{int}}} \left(\log \left|\det DG_{\boldsymbol{\theta}}(\boldsymbol{\xi}_k)\right|\right)^2
\end{equation}
\end{itemize}

It is worth mentioning the absence of any Jacobian-sign penalty, since $\det DG_{\boldsymbol{\theta}} > 0$ by the NSF architecture (Eq.~\ref{eq:full_logdet}) and thus orientation preservation is guaranteed without an explicit constraint. Similarly, no cycle-consistency loss is needed, as the NSF inverse is exact (up to floating-point precision). The training objective is therefore considerably simpler than those required by Neural ODE or other unconstrained parameterisations.

\subsection{Multi-shape conditioning via a two-stage training strategy}
\label{subsec:two_stage}
 
Extending the single-shape NSF to a family of domains requires a mechanism that conditions the spline parameters on the target shape. We adopt a two-stage strategy: in Stage~1, an auto-decoder framework~\cite{Park2019} jointly learns the flow parameters and a per-shape latent code for each training domain, while in Stage~2, a lightweight CNN encoder is trained to predict the latent code directly from the SDF image, supervised by the codes obtained in Stage~1. At inference time, the encoder maps a new SDF image to a latent code in a single forward pass, and the conditioned flow produces the diffeomorphism.
 
\subsubsection{Stage 1: auto-decoder training}
 
Each training shape $\Omega_i$, $i=1,\ldots,N_{shape}$, is assigned a learnable latent code $\boldsymbol{z}_i \in \mathbb{R}^{d_z}$, stored as a row in an embedding table. No encoder network is involved at this stage. The conditioned NSF defines a family of diffeomorphisms
\begin{equation}\label{eq:conditioned_nsf}
G_{\boldsymbol{\theta}}(\,\cdot\,;\boldsymbol{z}): \Omega^0 \to \Omega_{\boldsymbol{z}}
\end{equation}
where the spline parameters in each coupling layer are predicted by MLPs that receive both the passive coordinate and $\boldsymbol{z}$ as input. Different values of $\boldsymbol{z}$ produce different spline shapes and therefore different diffeomorphisms.
 
The latent codes $\{\boldsymbol{z}_i\}_{i=1}^{N_{shape}}$ are optimised jointly with the flow parameters $\boldsymbol{\theta}$, using separate learning rates. In practice, the flow parameters train at a lower rate $\eta_{\mathrm{flow}}$, gradually acquiring shared deformation structure, while the latent codes train at a higher rate $\eta_{\boldsymbol{z}} \gg \eta_{\mathrm{flow}}$, allowing each $\boldsymbol{z}_i$ to rapidly locate the correct region of latent space for its shape. The training loss is
\begin{equation}\label{eq:autodecoder_loss}
\mathcal{L}_{\mathrm{S1}} = \frac{1}{N_{shape}} \sum_{i = 1}^{N_{shape}} \left(
    \mathcal{L}_{\mathrm{int}}^{(i)} + \lambda_\partial \, \mathcal{L}_\partial^{(i)} + \lambda_{\log\det} \, \mathcal{L}_{\log\det}^{(i)}
\right)
+ \lambda_{\boldsymbol{z}} \frac{1}{N_{shape}} \sum_{i=1}^{N_{shape}} \|\boldsymbol{z}_i\|^2
\end{equation}
where $\mathcal{L}_{\mathrm{int}}^{(i)}$, $\mathcal{L}_\partial^{(i)}$, and $\mathcal{L}_{\log\det}^{(i)}$ are the per-shape losses from Eqs.~\eqref{eq:loss_interior}--\eqref{eq:loss_logdet}, evaluated using the conditioned flow $G_{\boldsymbol{\theta}}(\cdot\,;\boldsymbol{z}_i)$. The final term is an $L^2$ penalty on the latent codes that prevents them from drifting far from the origin and keeps the latent space compact. At the end of Stage~1, the flow has learned a smooth family of deformations parameterised by $\boldsymbol{z}$, and each training shape possesses a high-quality latent code $\boldsymbol{z}_i^*$ that produces an accurate diffeomorphism when passed to the flow.
 
\subsubsection{Stage 2: encoder training}
 
The auto-decoder provides accurate latent codes for training shapes, but it cannot produce codes for unseen shapes without ground-truth harmonic map data. To enable instant inference on new geometries, we train a CNN encoder $\mathrm{Enc}_{\boldsymbol{\theta}_{E}}: I \mapsto \boldsymbol{z}$ that maps the SDF image directly to the latent space.
 
In Stage~2, the flow parameters $\boldsymbol{\theta}$ are frozen and the stored codes $\{\boldsymbol{z}_i^*\}_{i=1}^{N_{shape}}$ from Stage~1 serve as regression targets. The encoder is trained to minimise
\begin{equation}\label{eq:encoder_loss}
\mathcal{L}_{\mathrm{S2}} = \frac{1}{N} \sum_{i=1}^N \Big(
    \lambda_{\mathrm{code}} \, \|\mathrm{Enc}_{\boldsymbol{\theta}_{E}}(I_i) - \boldsymbol{z}_i^*\|^2
    + \lambda_{\mathrm{map}} \, \mathcal{L}_{\mathrm{int}}^{(i)}
    + \lambda_{\partial} \, \mathcal{L}_\partial^{(i)}
\Big)
\end{equation}
where $\mathcal{L}_{\mathrm{int}}^{(i)}$ and $\mathcal{L}_\partial^{(i)}$ are the interior and boundary mapping errors evaluated using the predicted code $\hat{\boldsymbol{z}}_i = \mathrm{Enc}_{\boldsymbol{\theta}_{E}}(I_i)$ and the frozen flow $G_{\boldsymbol{\theta}}(\cdot\,;\hat{\boldsymbol{z}}_i)$. The first term drives the encoder toward the known-good codes, while the second and third terms ensure that the predicted codes actually produce accurate mappings, guarding against codes that are numerically close to $\boldsymbol{z}_i^*$ but fall in a region of latent space where the flow behaves differently.
 
The two-stage decomposition was adopted to improve training stability by avoiding the moving-target problem that arises when the encoder and latent representation are learned simultaneously. By contrast, Stage~2 is a standard supervised regression: the targets $\boldsymbol{z}_i^*$ are fixed, the flow is frozen, and the encoder simply learns a mapping from images to a well-defined latent space. 

\begin{figure}[H]
  \centering
  \includegraphics[width=\textwidth]{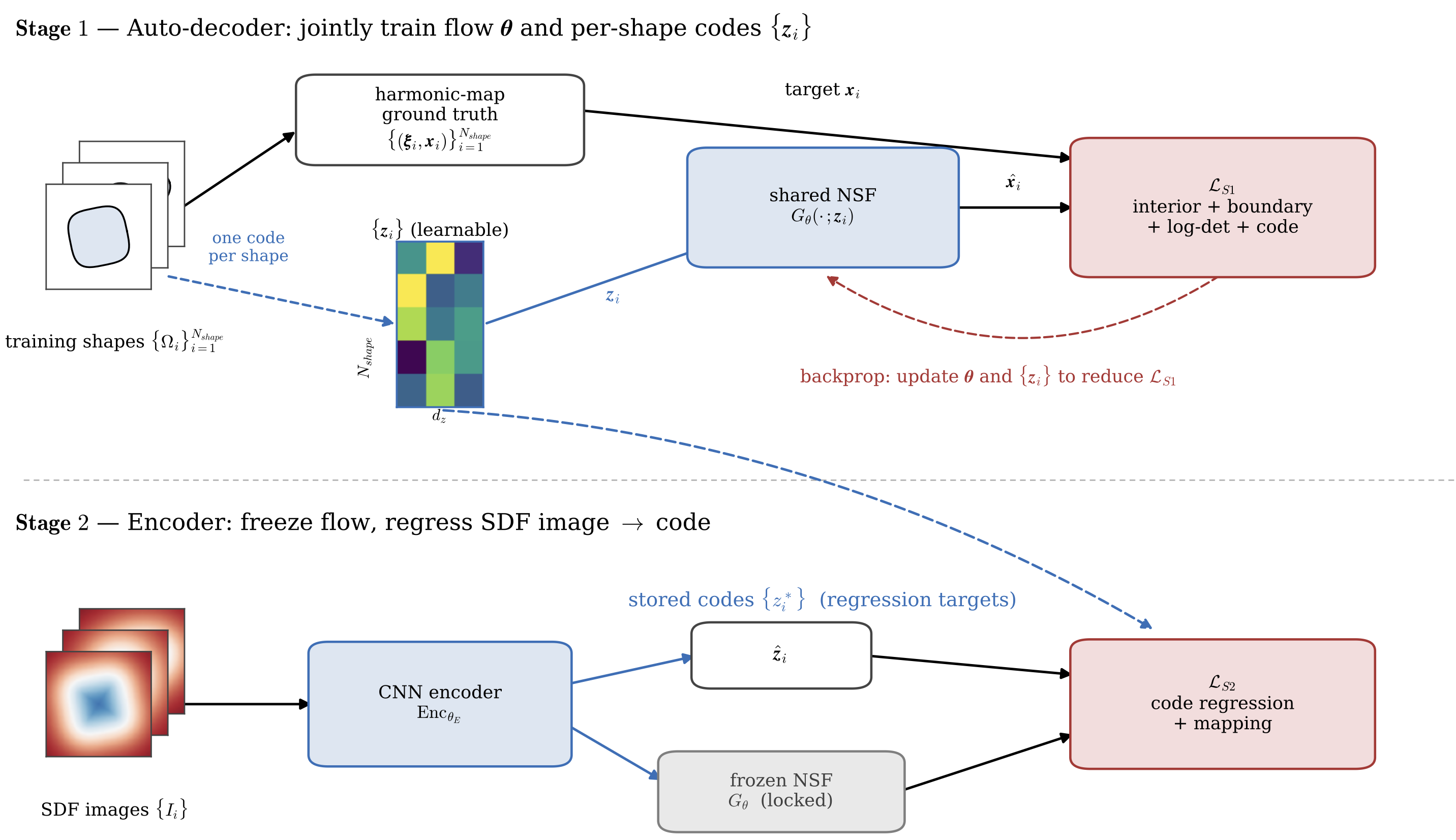}
  \caption{Two-stage training of the multi-shape model. \textbf{Stage 1} (auto-decoder):
  training shapes and a learnable table of per-shape latent codes $\{\boldsymbol{z}_i\}$ condition the shared NSF $G_{\theta}(\cdot\,;\boldsymbol{z}_i)$; the flow parameters $\boldsymbol{\theta}$ and the codes are optimised jointly against the harmonic-map ground truth.
  \textbf{Stage 2} (encoder): with the flow frozen, a CNN encoder $\mathrm{Enc}_{\theta_E}$
  maps each SDF image $I_i$ to a predicted code $\hat{\boldsymbol{z}}_i$, supervised by the stored codes $\{\boldsymbol{z}_i^*\}$ from Stage~1 (regression) together with the resulting mapping quality. The code table is the hinge between the two stages.}
  \label{fig:two_stage_training}
\end{figure}

The proposed two-stage training approach is illustrated in Fig.~\ref{fig:two_stage_training}
 
\subsubsection{Inference on new shapes}
\label{subsubsec:inference}
At test time, the full inference pipeline for a new domain $\Omega^*$ is:
\begin{equation}\label{eq:inference_pipeline}
I^* \;\xrightarrow{\;\mathrm{Enc}_{\boldsymbol{\theta}_{E}}\;}\; \boldsymbol{z}^* \;\xrightarrow{\;NSF\;}\; G_{\boldsymbol{\theta}}(\,\cdot\,;\boldsymbol{z}^*): \Omega^0 \to \Omega^*
\end{equation}
where $I^*$ is the SDF image of the new domain, computed on a regular grid without any meshing. The encoder produces the latent code in a single forward pass (a few milliseconds), and the conditioned flow provides the diffeomorphism $G_{\boldsymbol{\theta}}(\,\cdot\,;\boldsymbol{z}^*)$ that can be evaluated at arbitrary points. The mapping $G_{\boldsymbol{\theta}}$ is then used in the FEM pipeline exactly as described in Section~\ref{sec:ref_fem}, with the Jacobian $\boldsymbol{J} = DG_{\boldsymbol{\theta}}$ and its determinant evaluated at quadrature points via automatic differentiation.
 
The outcome of the proposed two-stage learning procedure is twofold. First, the shared flow  acquires the common structure of how reference squares deform into the training shapes, with the latent code $\boldsymbol{z}$ controlling shape-specific parameters such as aspect ratio and curvature. Next, the encoder learns to extract the geometric features of the SDF image (e.g. boundary curvature, aspect ratio, area) and map them to the corresponding region of latent space.

\begin{remark}[Diffeomorphism quality metrics]\label{rem:quality_metrics}
The following quality metrics can be used to assess the accuracy of the learned diffeomorphism:
\begin{enumerate}
\item \textit{Boundary fidelity.} We evaluate the mapping at boundary points of the reference
square and verify that the images lie on the physical boundary, reporting
\begin{equation}\label{eq:boundary_metric}
E_\partial = \max_{\boldsymbol{\xi} \in \partial \Omega^0}
\left| c\!\left(G_{\boldsymbol{\theta}}(\boldsymbol{\xi};\boldsymbol{z})\right) \right|
\end{equation}
where $c(\boldsymbol{x}) = 0$ defines the boundary curve (e.g.\ $(x_1/a)^2 + (x_2/b)^2 - 1 = 0$
for an ellipse). Values near zero indicate accurate boundary mapping.

\item \textit{Distortion.} Let $\sigma_1(\boldsymbol{\xi}) \ge \sigma_2(\boldsymbol{\xi}) > 0$
denote the singular values of $\boldsymbol{J}(\boldsymbol{\xi}) = DG_{\boldsymbol{\theta}}(\boldsymbol{\xi};\boldsymbol{z})$. The condition number $\kappa_J = \sigma_1 / \sigma_2$ quantifies the local anisotropy of the mapping, since large values correlate with anisotropic geometric coefficients $\boldsymbol{A}_{\mathrm{geo}}$ in the pulled-back weak form in Eq.~\eqref{eq:diffusion_pullback} and with a worse-conditioned discrete system.

\item \textit{Mapping RMSE.} On the training or test data we report the root-mean-square error
between predicted and ground-truth physical coordinates,
\begin{equation}\label{eq:mse_metric}
\mathrm{RMSE} = \left( \frac{1}{N} \sum_{k=1}^{N}
\left\| G_{\boldsymbol{\theta}}(\boldsymbol{\xi}_k;\boldsymbol{z}) - \boldsymbol{x}_k \right\|^2
\right)^{1/2}
\end{equation}

\end{enumerate}
\end{remark}

\section{Algorithm and Computational Considerations}
\label{sec:implementation}

This section consolidates the mathematical framework of Sections~\ref{sec:ref_fem}--\ref{sec:nsf} into a concrete algorithm and discusses the computational structure that arises from solving PDEs on a fixed reference mesh. We distinguish between an offline \emph{training phase} (performed once per shape family) and an online \emph{inference phase} (performed for each new shape and PDE instance), and highlight the structural advantages that the reference-square formulation provides for efficient implementation.

\subsection{End-to-end algorithm}
\label{subsec:algorithm}

The complete pipeline consists of an offline training phase (Algorithms~\ref{alg:stage1}--\ref{alg:stage2}) and an online inference phase (Algorithm~\ref{alg:inference}).

\begin{algorithm}[H]
\caption{Offline Stage 1: Auto-Decoder Training}\label{alg:stage1}
\begin{algorithmic}[1]
\Require Dataset of shapes $\{\Omega_i\}_{i=1}^{N_{shape}}$ with harmonic map ground truth $\left\{\left(\boldsymbol{\xi}_{ij}, \boldsymbol{x}_{ij}\right)_{j=1}^{N_{\mathrm{train}}}\right\}_{i=1}^{N_{\mathrm{shape}}}$ and SDF images $\{I_i\}_{i=1}^{N_{shape}}$
\Ensure Trained flow parameters $\boldsymbol{\theta}$ and latent codes $\{\boldsymbol{z}_i^*\}_{i=1}^{N_{shape}}$
\State Initialise flow parameters $\boldsymbol{\theta}$ and latent codes $\boldsymbol{z}_i \sim \mathcal{N}(\mathbf{0}, \sigma^2 \mathbf{I})$ for $i=1,\ldots,N_{shape}$
\For{epoch $= 1, \ldots, E_1$}
    \For{each mini-batch $\mathcal{B} \subset \{1,\ldots,N_{shape}\}$}
        \For{each shape $i \in \mathcal{B}$}
            \State Evaluate $G_{\boldsymbol{\theta}}(\boldsymbol{\xi}_i; \boldsymbol{z}_i)$ for interior and boundary points \Comment{NSF forward}
            \State Compute per-shape losses $\mathcal{L}_{\mathrm{int}}^{(i)}$, $\mathcal{L}_\partial^{(i)}$, $\mathcal{L}_{\log\det}^{(i)}$ via Eqs.~\eqref{eq:loss_interior}--\eqref{eq:loss_logdet}
        \EndFor
        \State Compute batch loss $\mathcal{L}_{\mathrm{S1}}$ via Eq.~\eqref{eq:autodecoder_loss} 
        \State Update $\boldsymbol{\theta}$ with learning rate $\eta_{\mathrm{flow}}$ and $\{\boldsymbol{z}_i\}_{i\in\mathcal{B}}$ with learning rate $\eta_{\boldsymbol{z}}$
    \EndFor
\EndFor
\State Store converged codes $\{\boldsymbol{z}_i^*\}_{i=1}^{N_{shape}}$ as regression targets for Stage~2
\end{algorithmic}
\end{algorithm}

\begin{algorithm}[H]
\caption{Offline Stage 2: Encoder Training}\label{alg:stage2}
\begin{algorithmic}[1]
\Require Frozen flow $G_{\boldsymbol{\theta}}$, stored codes $\{\boldsymbol{z}_i^*\}_{i=1}^{N_{shape}}$, SDF images $\{I_i\}_{i=1}^{N_{shape}}$, training data $\{(\boldsymbol{\xi}_{ij}, \boldsymbol{x}_{ij})\}$
\Ensure Trained encoder parameters $\boldsymbol{\theta}_E$
\State Initialise encoder $\mathrm{Enc}_{\boldsymbol{\theta}_E}$
\For{epoch $= 1, \ldots, E_2$}
    \For{each mini-batch $\mathcal{B} \subset \{1,\ldots,N_{shape}\}$}
        \State Predict latent codes: $\hat{\boldsymbol{z}}_i = \mathrm{Enc}_{\boldsymbol{\theta}_E}(I_i)$ for $i \in \mathcal{B}$ \Comment{CNN forward}
        \State Evaluate mappings $G_{\boldsymbol{\theta}}(\boldsymbol{\xi}_i; \hat{\boldsymbol{z}}_i)$ using predicted codes \Comment{Frozen flow}
        \State Compute $\mathcal{L}_{\mathrm{S2}}$ via Eq.~\eqref{eq:encoder_loss}: code regression + mapping quality
        \State Update $\boldsymbol{\theta}_E$ only \Comment{Flow $\boldsymbol{\theta}$ remains frozen}
    \EndFor
\EndFor
\end{algorithmic}
\end{algorithm}

\begin{algorithm}[H]
\caption{Online: PDE Solution on a New Domain $\Omega^*$}\label{alg:inference}
\begin{algorithmic}[1]
\Require Trained encoder $\mathrm{Enc}_{\boldsymbol{\theta}_E}$, trained flow $G_{\boldsymbol{\theta}}$; SDF image $I^*$ of $\Omega^*$; PDE data $(\boldsymbol{f}, \boldsymbol{u}_D, \boldsymbol{t}_N)$; reference mesh $\mathcal{T}^0_h$
\Ensure Discrete solution $\hat{\boldsymbol{u}}$ on $\Omega^0$ and physical-domain representation $\{(\boldsymbol{x}_j, \boldsymbol{u}_j)\}$
\Statex \hrulefill
\Statex \textbf{Stage 1: Encode shape}
\State Compute SDF image $I^*$ on a regular grid (no meshing required)
\State Predict latent code: $\boldsymbol{z}^* = \mathrm{Enc}_{\boldsymbol{\theta}_E}(I^*)$ \Comment{Single forward pass}
\Statex \hrulefill
\Statex \textbf{Stage 2: Precompute geometric quantities}
\State Collect all quadrature points $\{\boldsymbol{\xi}_q\}$ from $\mathcal{T}^0_h$
\State Evaluate $\boldsymbol{x}_q = G_{\boldsymbol{\theta}}(\boldsymbol{\xi}_q; \boldsymbol{z}^*)$ in batch \Comment{NSF forward pass}
\State Evaluate $\boldsymbol{J}_q = DG_{\boldsymbol{\theta}}(\boldsymbol{\xi}_q; \boldsymbol{z}^*)$ in batch \Comment{Automatic differentiation}
\State Compute $\mathcal{J}_q = \det \boldsymbol{J}_q$, $\;\boldsymbol{J}_q^{-1}$, $\;\boldsymbol{A}_{\mathrm{geo},q} = \mathcal{J}_q \, \boldsymbol{J}_q^{-1} \boldsymbol{A}(\boldsymbol{x}_q) \, \boldsymbol{J}_q^{-T}$ for each $q$
\Statex \hrulefill
\Statex \textbf{Stage 3: Assemble and solve}
\State Assemble $\hat{\mathbf{K}}$ and $\hat{\mathbf{f}}$ on $\mathcal{T}^0_h$ using precomputed $\{\boldsymbol{A}_{\mathrm{geo},q}, \mathcal{J}_q, \boldsymbol{x}_q\}$
\State Apply boundary conditions on $\partial\Omega^0$ (with values $\boldsymbol{u}_D(G_{\boldsymbol{\theta}}(\boldsymbol{\xi}; \boldsymbol{z}^*))$ on $\Gamma^0_D$)
\State Solve $\hat{\mathbf{K}}\,\hat{\mathbf{u}} = \hat{\mathbf{f}}$
\Statex \hrulefill
\Statex \textbf{Stage 4: Map solution to physical domain}
\State For each reference node $\boldsymbol{\xi}_j$, compute $\boldsymbol{x}_j = G_{\boldsymbol{\theta}}(\boldsymbol{\xi}_j; \boldsymbol{z}^*)$
\State Return $\{(\boldsymbol{x}_j, \hat{\boldsymbol{u}}_j)\}$ for post-processing 
\end{algorithmic}
\end{algorithm}

A key feature of this decomposition is that the neural network is invoked only in Stages~1 and~2 of Algorithm~\ref{alg:inference}. The encoder evaluation (Stage~1) consists of a single CNN forward pass requiring milliseconds, while the geometric precomputation (Stage~2) involves a batched evaluation of the flow and its Jacobian. Once these quantities are cached, the finite element assembly and solve in Stage~3 proceed with standard numerical linear algebra and involve no neural network evaluations. The neural network therefore handles the geometric representation, while the finite element method retains responsibility for the discretisation and numerical solution of the governing PDE on the learned domain. Notably, the online algorithm requires only the SDF image of the new domain as input, avoiding the need to generate a new boundary-conforming mesh for $\Omega^*$ or compute a harmonic map.

\subsection{Efficient evaluation of geometric quantities}
\label{subsec:geom_quantities}

The reference-domain FEM requires three geometric quantities at every quadrature point $\boldsymbol{\xi}_q$:
\begin{enumerate}[label=(\roman*),leftmargin=2em]
    \item the physical coordinate $\boldsymbol{x}_q = G_{\boldsymbol{\theta}}(\boldsymbol{\xi}_q; \boldsymbol{z})$,
    \item the Jacobian matrix $\boldsymbol{J}_q = DG_{\boldsymbol{\theta}}(\boldsymbol{\xi}_q; \boldsymbol{z}) \in \mathbb{R}^{2\times 2}$,
    \item the Jacobian determinant $\mathcal{J}_q = \det \boldsymbol{J}_q$.
\end{enumerate}

\noindent These can be computed very efficiently as follows:

\begin{enumerate}
\item \textit{Batched forward evaluation}:
The physical coordinates are obtained by a single batched forward pass of the NSF through all $N_q$ quadrature points simultaneously. On GPU or accelerator hardware, this is a dense matrix computation with cost $O(N_q \cdot L \cdot C_{\mathrm{MLP}})$, where $L$ is the number of coupling layers and $C_{\mathrm{MLP}}$ is the cost of one MLP evaluation. For the mesh sizes considered in this work ($N_q \sim 10^3$--$10^4$), this evaluation completes in milliseconds.

\item \textit{Log-determinant from the forward pass}:
The NSF coupling-layer structure provides $\log\mathcal{J}_q$ analytically as a byproduct of the forward pass (Eq.~\ref{eq:full_logdet}), at negligible additional cost. The determinant is recovered as $\mathcal{J}_q = \exp(\log\mathcal{J}_q)$.

\item \textit{Full Jacobian via automatic differentiation}:
The $2\times 2$ Jacobian matrix $\boldsymbol{J}_q$ is obtained via reverse-mode automatic differentiation through the forward pass. For $d_s = 2$ output dimensions, this requires two backward passes, each computing one row of $\boldsymbol{J}_q$. Crucially, these backward passes are batched over all $N_q$ quadrature points, yielding the full set $\{\boldsymbol{J}_q\}_{q=1}^{N_q}$ in two vectorised operations. The Jacobian inverse $J_q^{-1}$ and the geometric coefficient tensor $\boldsymbol{A}_{\mathrm{geo},q}$ are then computed via standard $2\times 2$ matrix algebra.

\item \textit{Precomputation and caching}:
Since the geometric quantities depend only on the mesh and the mapping (not on the PDE data), they can be precomputed once and reused across multiple right-hand sides, boundary conditions, or parameter sweeps on the same geometry. This is particularly advantageous in parametric studies where the domain is fixed but the loading or material properties vary.
\end{enumerate}

\subsection{Discussion on computational opportunities enabled by the fixed reference mesh}
\label{subsec:fixed_mesh_assembly}
\label{subsec:solvers}

The primary computational advantage of the proposed formulation is the elimination of
geometry-specific meshing. Once the diffeomorphic map has been obtained, every PDE solve is
performed on the same fixed, structured mesh $\mathcal{T}^0_h$ defined on
$\Omega^0=(-1,1)^2$, rather than on a new mesh fitted to each physical domain.
This removes mesh generation and mesh-quality control from the online analysis of a new
geometry.

Beyond this main advantage, the fixed reference-domain formulation also exposes several
additional computational opportunities that could be exploited in more specialised or
high-performance implementations. These possibilities follow from the regularity and
geometry-independent topology of the reference mesh, but are not all implemented or
quantitatively benchmarked in the present work. We summarise the most relevant ones below.

\begin{enumerate}

\item \textbf{Fixed sparsity pattern and reusable algebraic structure.}
Because the topology of the reference mesh is unchanged across geometries, the sparsity
pattern of the discrete operator $\hat{\mathbf{K}}$ is also fixed. Consequently, data structures
associated with sparse storage, graph-based orderings, and other topology-dependent solver
operations can in principle be constructed once and reused for different physical domains.
Only the numerical entries of the operator change through the geometry-dependent coefficient
field $\boldsymbol{A}_{\mathrm{geo}}$.

This property may be particularly useful in many-query settings. Across different geometries,
for example, a sparse direct solver could reuse its fill-reducing ordering and symbolic
factorisation, requiring only numerical refactorisation for each new coefficient field.
For repeated solves on the same geometry with different right-hand sides, the numerical
factorisation itself may also be reused.

\item \textbf{Regular memory access and vectorisation.}
A regular structured mesh permits predictable and contiguous layouts of nodal, elemental,
and quadrature-point data. This structure is favourable for cache-efficient assembly and
matrix--vector products and is also well suited to SIMD vectorisation and accelerator-based
execution.

In particular, all reference elements have identical topology and geometry, so the same
element kernel can be applied throughout the mesh. Combined with the batched evaluation of
the NSF map and its Jacobian described in Section~\ref{subsec:geom_quantities}, this structure
provides a natural basis for GPU-oriented or otherwise vectorised implementations.

\item \textbf{Separation of geometry-dependent and geometry-independent quantities.}
The discrete operator can be written schematically as

\begin{equation}
\label{eq:assembly_decomposition}
\hat{K}_{ij}
=
\sum_e
\int_{\Omega^0_e}
(\nabla \hat{\varphi}_i)^T
\boldsymbol{A}_{\mathrm{geo}}(\boldsymbol{\xi})
\nabla \hat{\varphi}_j
\,\mathrm{d}\boldsymbol{\xi}.
\end{equation}

The basis-function values and gradients depend only on the reference element, whereas all
information about the physical geometry is contained in
$\boldsymbol{A}_{\mathrm{geo}}(\boldsymbol{\xi})$. Since the reference elements are congruent,
shape-function values, gradients, quadrature weights, and other reference-element quantities
can be precomputed once and reused throughout the mesh and across geometries. This separation also makes the formulation naturally compatible with matrix-free approaches, in which $\hat{\mathbf{K}}$ need not be assembled explicitly and the action $\hat{\mathbf{K}}\mathbf{v}$ can instead be evaluated element by element using the precomputed reference quantities together with the geometry-dependent coefficients.

\end{enumerate}

These additional opportunities should therefore be viewed as consequences of the computational structure created by the diffeomorphic pullback, rather than as optimisations required for the method itself. In the present implementation, the central benefit is that the complex physical geometry does not need to be discretised by a new boundary-conforming mesh for each analysis. The price of working on the fixed reference mesh is that the geometry-dependent quantities, including $\boldsymbol{J}$, $\mathcal{J}$, $\boldsymbol{J}^{-1}$, and $\boldsymbol{A}_{\mathrm{geo}}$, must be evaluated for each new shape, and the resulting non-polynomial integrands may benefit from higher-order quadrature when the geometric coefficients vary strongly within an element. Moreover, the learned diffeomorphism will inevitably exhibit some discrepancies  between the physical domain and the mapped domain. Lastly, strongly distorted mappings may lead to anisotropic transformed coefficients and poorer conditioning. These considerations define the trade-off between geometric simplification and algebraic complexity in the proposed formulation.

\section{Numerical Experiments} 
\label{sec:experiments}
 
\noindent This section presents numerical experiments that validate the proposed reference-domain FEM framework, and demonstrate its applicability and performance on various cases in computational mechanics.

\subsection{Error metrics}
\label{subsec:error_metrics}
An important thing that needs to be established before proceeding to the experiments is the error metric used to assess the quality of pullback solutions, since in all cases there are three finite element solutions involved and this might cause ambiguities. In particular, let $\Omega$ be the physical domain and $G:\Omega^0\to\Omega$ the learned diffeomorphism. The three finite element solutions involved are:

\begin{itemize}
  \item $u_{\mathrm{ref}}$, either the analytical solution (Sec.~\ref{subsec:poisson_ellipse}) or the direct FE solution computed on a highly refined body-fitted mesh $\mathcal{T}_{\mathrm{ref}}$ ($\approx 4\times10^4$ elements) and taken as the reference solution (Sec.~\ref{subsec:poisson_superelliptical}, \ref{subsec:elasticity_blobs});
  \item $u_{\mathrm{dir}}$, the direct solve on a coarse body-fitted mesh $\mathcal{T}_{\mathrm{dir}}$;
  \item $\hat{u}_{\mathrm{pb}}$, the pullback solve on the structured reference
        mesh $\mathcal{T}^0$, whose physical counterpart is $u_{\mathrm{pb}} = \hat{u}_{\mathrm{pb}}\circ G^{-1}$.
\end{itemize}

The three meshes share no nodes, so both approximations are evaluated at the nodes $\boldsymbol{x}_i\in \mathcal{T}_{\mathrm{ref}}$ through their own interpolants:

\begin{equation}
u_{\mathrm{dir}}^{\,h}(\boldsymbol{x}_i)=\sum_{a} N_a^{\mathrm{dir}}(\boldsymbol{x}_i)\,
u_{\mathrm{dir},a},
\qquad
u_{\mathrm{pb}}^{\,h}(\boldsymbol{x}_i)=\sum_{a} N_a^{\mathrm{pb}}(\boldsymbol{x}_i)\,
\hat{u}_{\mathrm{pb},a}
\end{equation} 

\noindent where $N_a^{\mathrm{pb}}$ are the shape functions of the mapped mesh. 

The errors of both solutions are integrated on the same mesh $\mathcal{T}_{\mathrm{ref}}$ by quadrature 

\begin{equation}\label{eq:l2norm}
\|v\|_{L^2(\Omega)}
=\Bigl(\sum_{e\in\mathcal{T}_{\mathrm{ref}}}\ \sum_{q}
 w_q\,\bigl|\det \boldsymbol{J}_{e,q}\bigr|\;
 \bigl|v_h(\boldsymbol{x}_q)\bigr|^2\Bigr)^{1/2}
\end{equation}

\noindent where $\boldsymbol{J}_{e,q}$ denotes the isoparametric Jacobian of element $e$ at quadrature point $q$ of $\mathcal{T}_{\mathrm{ref}}$, and $|\cdot|$ the absolute value for the scalar problems (Sec.~\ref{subsec:poisson_ellipse},\ref{subsec:poisson_superelliptical}). Also, for vector-valued fields (Sec.~\ref{subsec:elasticity_blobs}) we compute:

\begin{equation}
\|\boldsymbol{v}\|_{L^2(\Omega)}^2=\int_\Omega\bigl(v_x^2+v_y^2\bigr)\,\mathrm{d}\boldsymbol{x}
\end{equation}

\noindent The reported quantities are relative:
\begin{equation}\label{eq:relerr}
\varepsilon_{\mathrm{dir}}=
\frac{\|u_{\mathrm{dir}}^{\,h}-u_{\mathrm{ref}}\|_{L^2}}{\|u_{\mathrm{ref}}\|_{L^2}}
\qquad
\varepsilon_{\mathrm{pb}}=
\frac{\|u_{\mathrm{pb}}^{\,h}-u_{\mathrm{ref}}\|_{L^2}}{\|u_{\mathrm{ref}}\|_{L^2}} 
\end{equation}

It is important to note that neither coarse discretisation covers $\Omega$ exactly. The mesh
$\mathcal{T}_{\mathrm{dir}}$ is a polygon inscribed in $\Omega$, omitting a boundary band of width $O(h^2/R)$ with $R$ the local radius of curvature, while $G(\Omega^0)$ reproduces $\partial\Omega$ only to within the training error of the map. Reference nodes falling outside either region are therefore not in the domain of the corresponding interpolant. We evaluate both errors on the common
support

\begin{equation}
\mathcal{K}=\bigl\{\,i:\ \boldsymbol{x}_i\in
\mathcal{T}_{\mathrm{dir}}\cap G(\Omega^0)\,\bigr\}
\end{equation}

\noindent setting $v_h(\boldsymbol{x}_i)=0$ for $i\notin\mathcal{K}$ in both the numerator and the denominator of \eqref{eq:relerr}, so that the two errors are integrated over the same region and remain directly comparable. In all cases reported here $|\mathcal{K}^{c}|$ amounts to $1$--$2\%$ of the reference nodes and consists entirely of nodes on $\partial\Omega$.

\subsection{Poisson problem on an elliptical domain}
\label{subsec:poisson_ellipse}
 
We begin with an illustrative example of a Poisson problem on an elliptical domain with a known analytical solution,  which allows a detailed assessment of the mapping quality, the transformed weak form, and the accuracy of the resulting finite element solution. For this example, the NSF model was trained to map a square to a single elliptical domain, without employing the multi-shape auto-decoder framework, in order to focus exclusively on the representational capacity of the NSF. 

Specifically, let us consider the Poisson equation with homogeneous Dirichlet boundary conditions on an elliptical domain $\Omega = \{(x_1, x_2) \in \mathbb{R}^2 : x_1^2/a^2 + x_2^2/b^2 < 1\}$ with semi-axes $a = 2$ and $b = 1$:
\begin{equation}\label{eq:poisson_ellipse}
-\Delta u = f \quad \text{in } \Omega, \qquad u = 0 \quad \text{on } \partial\Omega
\end{equation}
To enable exact error quantification, we employ the \emph{method of manufactured solutions}. The function
\begin{equation}\label{eq:manufactured_u}
u_{\mathrm{exact}}(x_1, x_2) = 1 - \frac{x_1^2}{a^2} - \frac{x_2^2}{b^2}
\end{equation}
satisfies $u_{\mathrm{exact}} = 0$ on $\partial\Omega$ by construction and yields the constant source term
\begin{equation}\label{eq:manufactured_f}
f = -\Delta u_{\mathrm{exact}} = \frac{2}{a^2} + \frac{2}{b^2} = 2.5
\end{equation}
This manufactured solution has a smooth, bowl-shaped profile that peaks at the centre of the ellipse ($u_{\mathrm{exact}}(0,0) = 1$) and vanishes on the boundary.
 
\begin{figure}[H]
    \centering

    \begin{subfigure}[b]{0.48\textwidth}
        \centering
        \includegraphics[width=\textwidth]{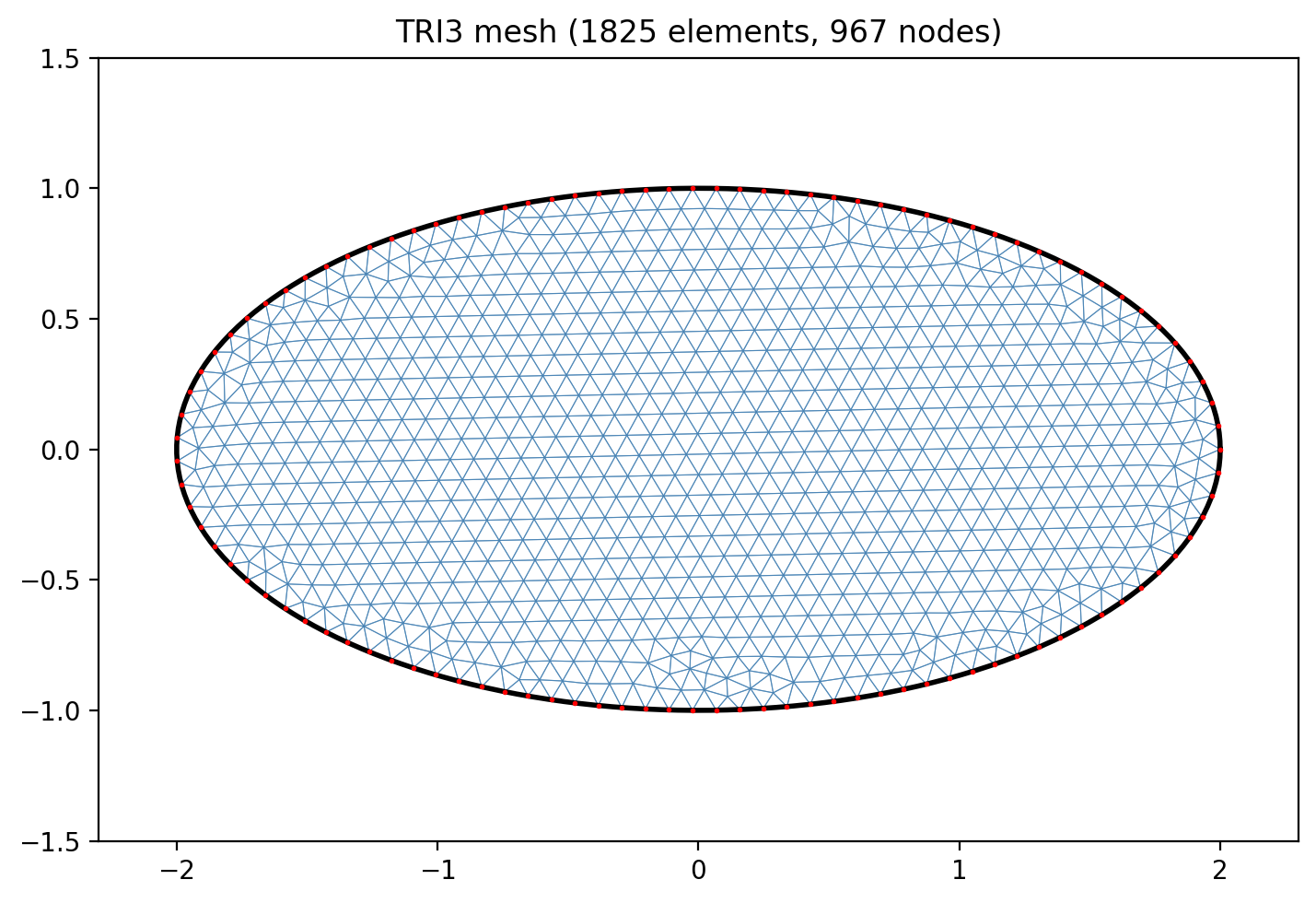}
        \caption{Triangular baseline mesh}
        \label{fig:example1_tri3_baselineMesh}
    \end{subfigure}
    \hfill
    \begin{subfigure}[b]{0.48\textwidth}
        \centering
        \includegraphics[width=\textwidth]{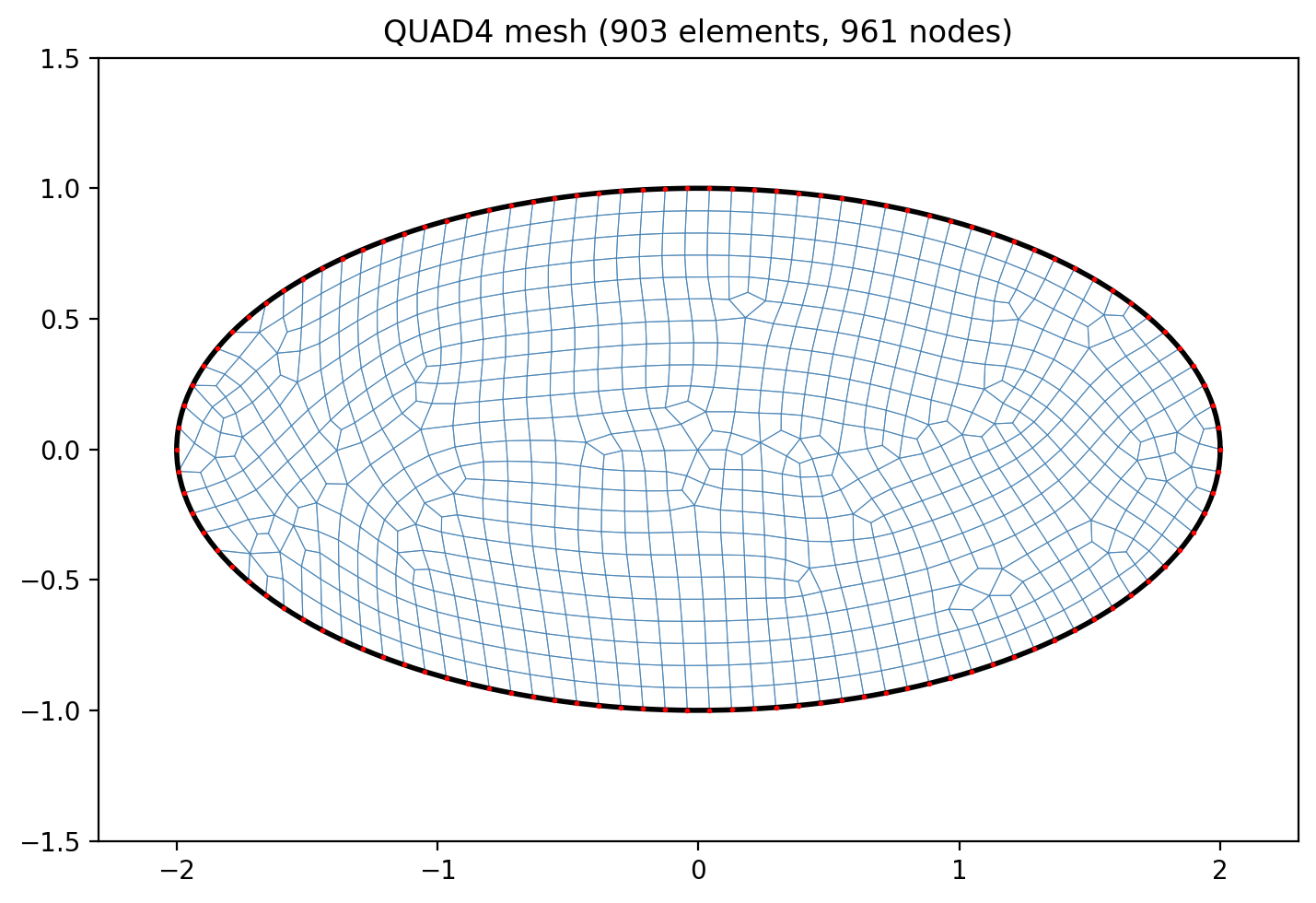}
        \caption{Quadrilateral baseline mesh}
        \label{fig:example1_quad4_baselineMesh}
    \end{subfigure}

    \begin{subfigure}[b]{0.40\textwidth}
        \centering
        \includegraphics[width=\textwidth]{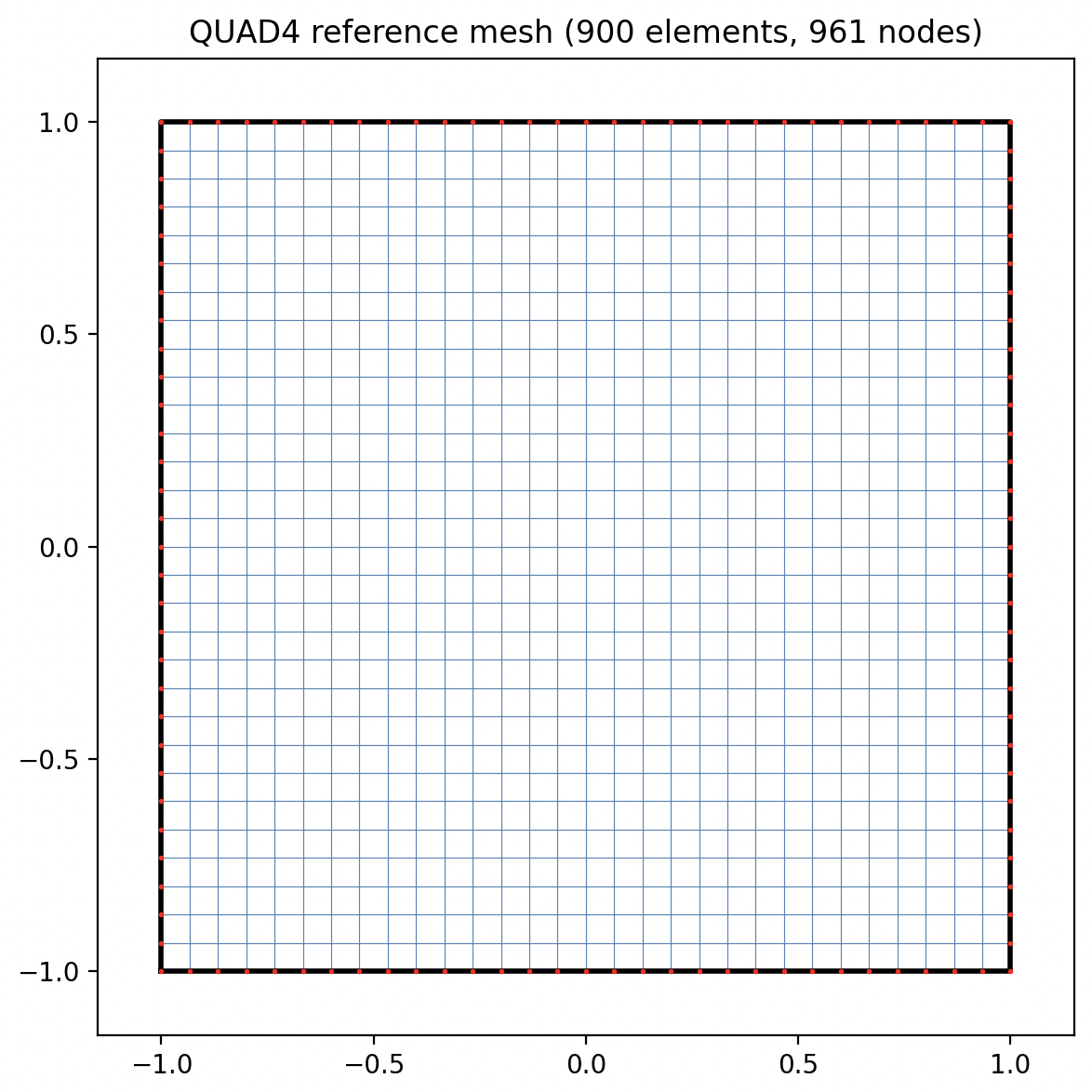}
        \caption{Quadrilateral mesh on the reference domain}
        \label{fig:example1_quad4_referenceMesh}
    \end{subfigure}
    
    \caption{Baseline meshes and mesh on the reference (square) domain}
    \label{fig:example1_baselineMeshes}
\end{figure}

As a baseline, we solve Eq.~\eqref{eq:poisson_ellipse} directly on $\Omega$ using both an unstructured triangular and an unstructured quadrilateral mesh generated by Gmsh~\cite{Geuzaine2009}. The triangular mesh consists of linear elements (Tri3) with a characteristic element size $h = 0.091$, yielding 967 nodes and 1825 elements. The stiffness matrix and load vector are assembled using the standard Galerkin FEM with one-point quadrature on each triangle. Similarly, the quadrilateral mesh consists of linear elements (Quad4) with a characteristic element size $h = 0.085$, yielding 961 nodes and 903 elements, using a $2\times 2$ Gauss integration to compute the element matrices and vectors. Figure \ref{fig:example1_baselineMeshes} displays both baseline meshes (\ref{fig:example1_tri3_baselineMesh},\ref{fig:example1_quad4_baselineMesh}), as well as the $30\times30$ quad element mesh on the reference domain (\ref{fig:example1_quad4_referenceMesh}). In all cases, the resulting linear systems are solved by a sparse direct solver.

The NSF is trained in the reversed direction $G_{\boldsymbol{\theta}}: \Omega^0 \to \Omega$ as described in Section~\ref{subsec:single_shape}. For this single-shape example, no auto-decoder conditioning is required, and the unconditional NSF simply maps the reference square $\Omega^0 = (-1,1)^2$ directly to the ellipse. The ground truth data  for training is generated by solving the harmonic map problem of Eq.~\eqref{eq:harmonic_map} on a $512 \times 512$ finite-difference grid. Interior points ($N_{\mathrm{in}} = 5000$) are sampled inside the ellipse, and boundary points ($N_b = 1000$) are sampled uniformly along the ellipse perimeter. For each interior point $\boldsymbol{x}_k \in \Omega$, the harmonic map solver provides the corresponding reference-square coordinate $\boldsymbol{\xi}_k = G^{-1}(\boldsymbol{x}_k )$. Thus, the training data for the NSF is the pairs $(\boldsymbol{\xi}_k, \boldsymbol{x}_k)$.
 
\begin{figure}[H]
  \centering
    \includegraphics[width=0.8\linewidth]{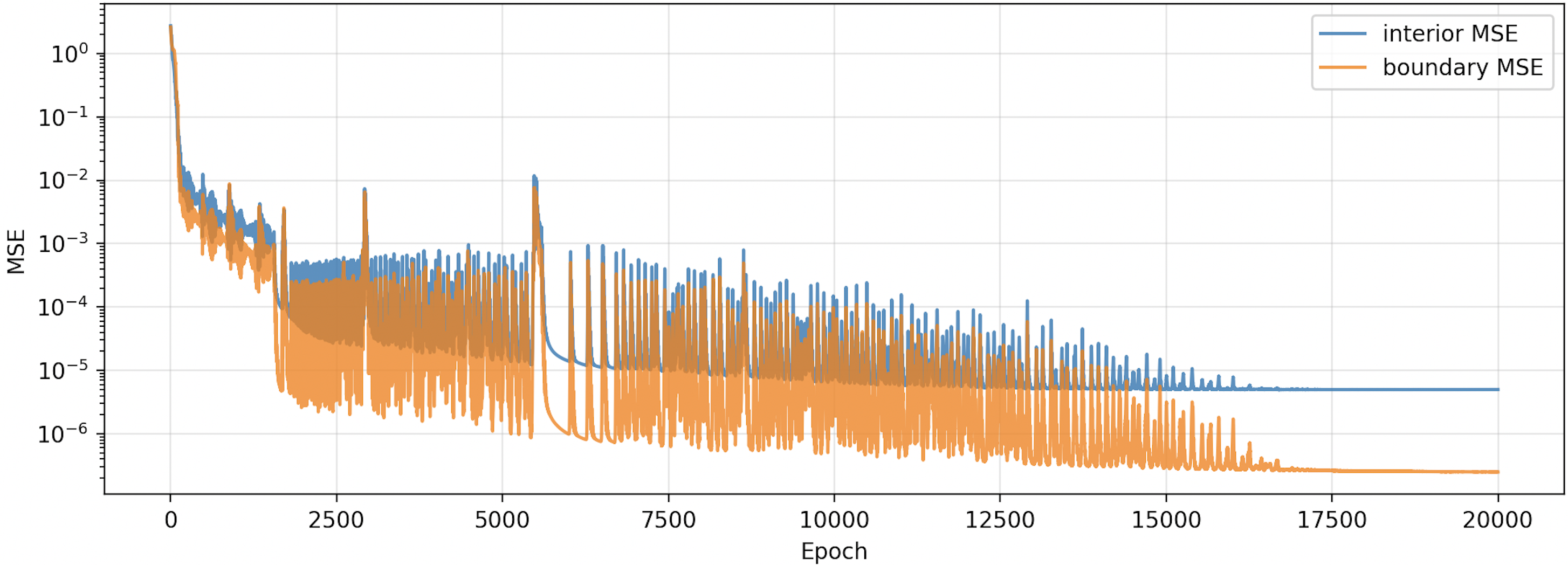}
  \caption{Evolution of the training loss for the NSF architecture.}
  \label{fig:ex1_trainingCurve}
\end{figure}

The NSF architecture consists of $L = 8$ coupling layers with $K = 20$ spline bins per layer, hidden dimension 128, and 3 hidden layers in each parameter-prediction MLP. The domain bounds are set to $(-1.10, 1.10)$ (input, slightly larger than the unit square) and $(-2.1, 2.1)$ (output, encompassing the ellipse). Training proceeds for 20000 epochs with the Adam optimiser, initial learning rate $1 \times 10^{-3}$, and cosine annealing to $10^{-5}$. After training, the model achieves an interior MSE of $5 \times 10^{-5}$ and a boundary MSE of $2 \times 10^{-6}$, as shown in the training curve depicted in Fig.~\ref{fig:ex1_trainingCurve}. The Jacobian determinant $\mathcal{J} = \det DG_{\boldsymbol{\theta}}$ remains strictly positive across all evaluation points, as expected, confirming the diffeomorphism property. 

\begin{figure}[H]
  \centering
    \includegraphics[width=0.8\linewidth]{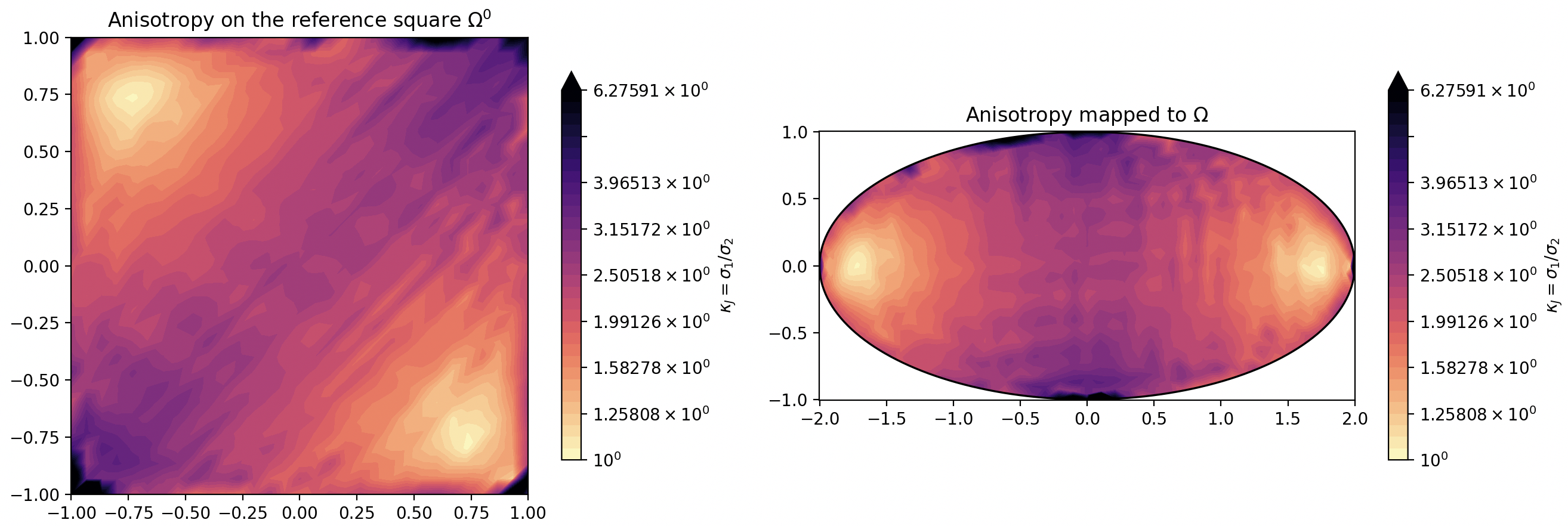}
  \caption{Distribution of anisotropy in the reference square and its mapping back to the physical domain.}
  \label{fig:ex1_anisotropy}
\end{figure}

Before using the trained NSF in the FEM pipeline, we verify its suitability through the diagnostics reported in Table \ref{tab:nsf_diagnostics} (also described in Section~\ref{rem:quality_metrics}). Also, figure~\ref{fig:ex1_anisotropy} depicts the anisotropy, computed as the ratio $\kappa_J=\sigma_1 / \sigma_2$, in the reference square, as well as its mapping back to the physical domain. 

\begin{table}[htbp]
\centering
\caption{Diagnostic checks used to assess the suitability of the trained NSF before its incorporation into the FEM pipeline.}
\label{tab:nsf_diagnostics}
\begin{tabular}{p{3.2cm} p{5.2cm} p{6.0cm}}
\hline
\textbf{Metric} & \textbf{Result} & {} \\
\hline

Boundary fidelity &
$\dfrac{x_1^2}{a^2} + \dfrac{x_2^2}{b^2} \in [0.997,\ 1.008]$ for all $\boldsymbol{\xi} \in \partial\Omega^0$ &
The mapped boundary deviates from the exact ellipse by less than $0.8\%$, which is considered sufficient. \\[0.3cm]

Round-trip consistency &
For 200 test points, the composition
$G_{\boldsymbol{\theta}}^{-1}(G_{\boldsymbol{\theta}}(\boldsymbol{\xi}))$
yields a mean absolute error of $O(10^{-6})$ &
This confirms the numerical stability and accuracy of the exact inverse provided by the NSF. \\

\hline
\end{tabular}
\end{table}
 
Finally, following Algorithm~\ref{alg:inference} (with the latent-code fitting stage omitted since this is a single-shape example), we solve the transformed Poisson problem on the reference square. The reference domain $\Omega^0= (-1,1)^2$ is discretised with a uniform grid of $30 \times 30 = 900$ bilinear quadrilateral elements (Quad4), yielding 961 nodes, as shown in Fig. \ref{fig:example1_quad4_referenceMesh}. Tensor-product $2 \times 2$ Gauss--Legendre quadrature is used for element integration. At each of the $900 \times 4 = 3600$ quadrature points, the NSF provides:
\begin{itemize}[leftmargin=2em]
    \item the physical coordinate $\boldsymbol{x}_q = G_{\boldsymbol{\theta}}(\boldsymbol{\xi}_q)$ via a batched forward pass,
    \item the full $2\times 2$ Jacobian $J_q = DG_{\boldsymbol{\theta}}(\boldsymbol{\xi}_q)$ via two batched backward passes,
    \item the determinant $\mathcal{J}_q = \det J_q$ from the log-determinant output of the forward pass.
\end{itemize}
The geometric coefficient tensor is then computed as $\boldsymbol{A}_{\mathrm{geo},q} = \mathcal{J}_q \, \boldsymbol{J}_q^{-1} \boldsymbol{A} \boldsymbol{J}_q^{-T}$ and the transformed source as $\hat{f}_q = f \cdot \mathcal{J}_q = 2.5 \, \mathcal{J}_q$. All geometric quantities are evaluated in a single batch and cached before the element assembly loop begins. The stiffness matrix and load vector are assembled using the pulled-back weak form from Eqs.~\eqref{eq:diffusion_pullback}--\eqref{eq:rhs_pullback}:
\begin{equation*}
\hat{K}_{ij} = \sum_e \sum_q w_q \, (\nabla\hat{\varphi}_i)^T \boldsymbol{A}_{\mathrm{geo},q} \, \nabla\hat{\varphi}_j \, |\det \boldsymbol{J}_{\mathrm{iso},q}|,
\qquad
\hat{f}_i = \sum_e \sum_q w_q \, \hat{\varphi}_i \, f \, \mathcal{J}_q \, |\det \boldsymbol{J}_{\mathrm{iso},q}|
\end{equation*}
where $\boldsymbol{J}_{\mathrm{iso}}$ is the isoparametric Jacobian of the reference-element mapping (which, for uniform Quad4 elements on $(-1,1)^2$, is simply a constant scaling). The resulting $961 \times 961$ sparse linear system is solved by a direct solver. Since the discrete solution $\hat{u}_h$ lives on the reference-square nodes, to visualise and compare with the exact solution, we compute the physical coordinates $\boldsymbol{x}_j = G_{\boldsymbol{\theta}}(\boldsymbol{\xi}_j)$ at each reference node via a batched forward pass and evaluate the error $|u_h(\boldsymbol{x}_j) - u_{\mathrm{exact}}(\boldsymbol{x}_j)|$ pointwise.
 
Table~\ref{tab:poisson_results} compares the direct and mapped FEM solutions for roughly the same number of nodes. The pullback formulation is less accurate than the direct one: its relative error is 0.72\%, against 0.23\% and 0.28\% for the Tri3 and Quad4 discretisations. The gap is mainly attributable to the learned mapping. Since the pullback is an exact reformulation (Proposition~\ref{prop:equivalence}), the finite element error on the reference square retains its classical behaviour, and the geometry enters only through the coefficients $\boldsymbol{A}_{geo}$. A positional error in $G_{\boldsymbol{\theta}}$ therefore propagates to the solution, which the reported figures reflect. The mapping accuracy, rather than the discretisation, is thus the limiting factor here, and it is a quantity that can be measured and improved independently of the solver.

\begin{figure}[H]
    \centering

    \begin{subfigure}[b]{0.48\textwidth}
        \centering
        \includegraphics[width=\linewidth]{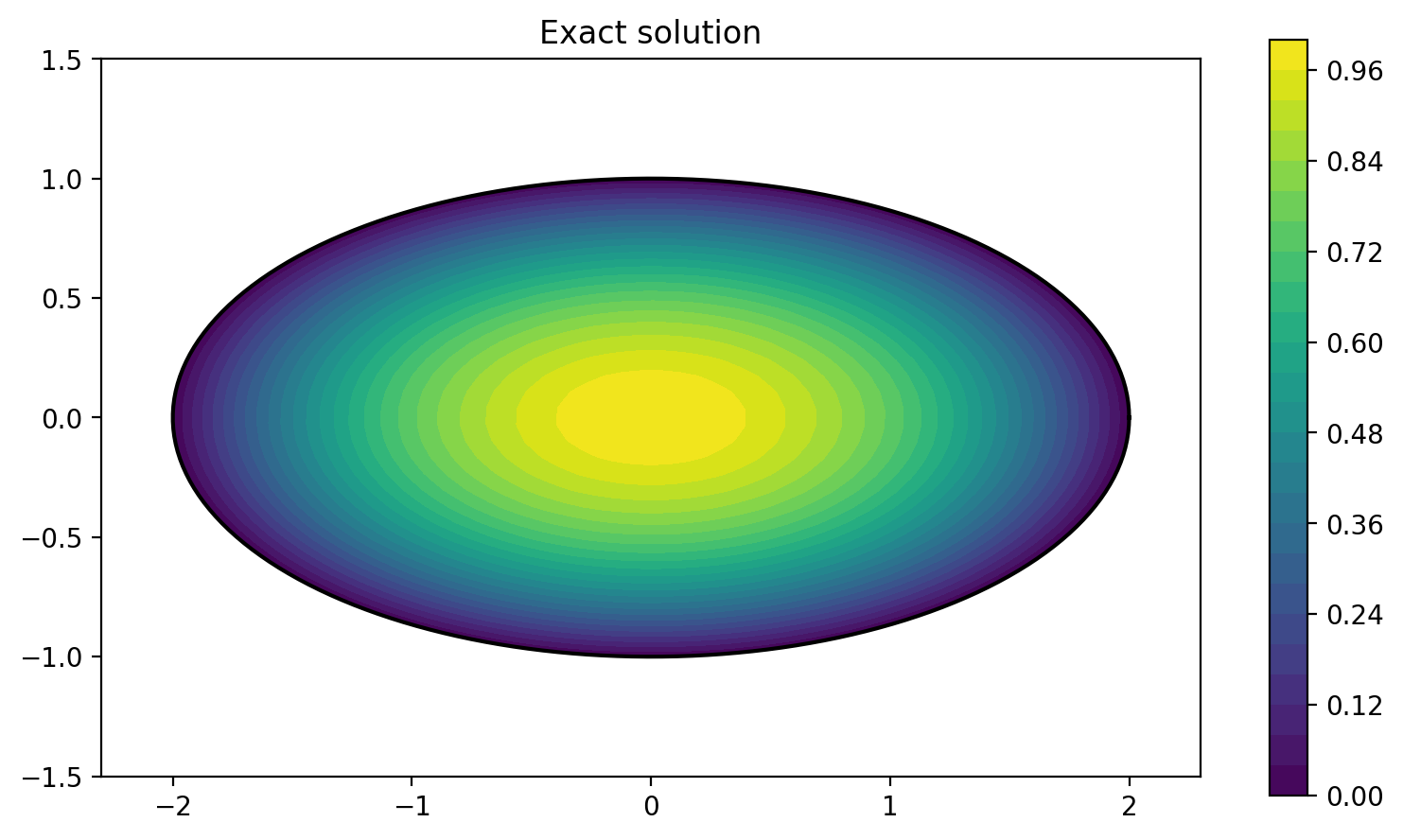}
        \caption{Exact solution}
        \label{fig:ex1_exact}
    \end{subfigure}
    \hfill
    \begin{subfigure}[b]{0.48\textwidth}
        \centering
        \includegraphics[width=\linewidth]{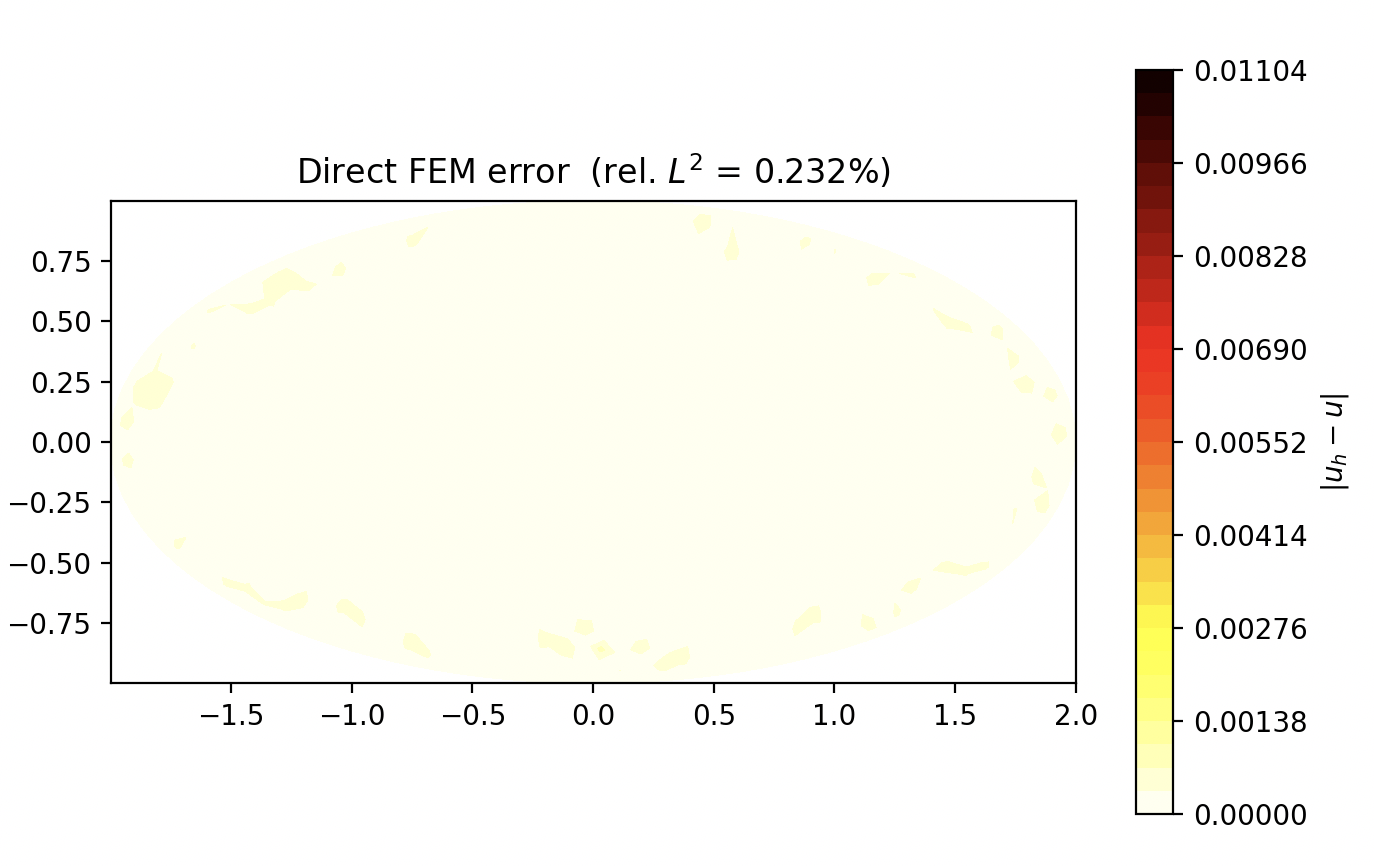}
        \caption{Pointwise absolute error map for the direct FEM with Tri3 elements}
        \label{fig:ex1_tri3}
    \end{subfigure}

    \vspace{0.5cm}

    \begin{subfigure}[b]{0.48\textwidth}
        \centering
        \includegraphics[width=\linewidth]{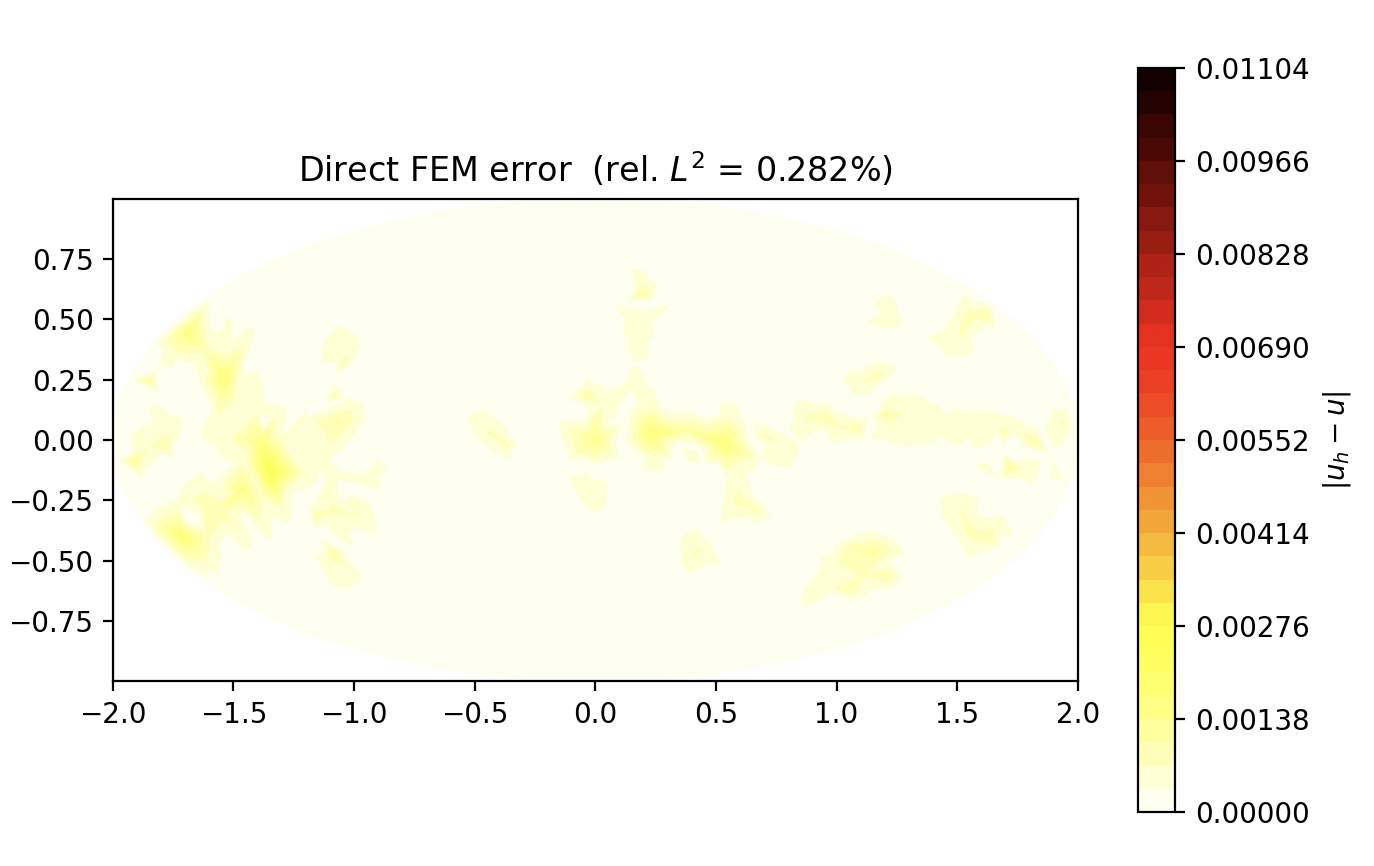}
        \caption{Pointwise absolute error map for the direct FEM with Quad4 elements}
        \label{fig:ex1_quad4}
    \end{subfigure}
    \hfill
    \begin{subfigure}[b]{0.48\textwidth}
        \centering
        \includegraphics[width=\linewidth]{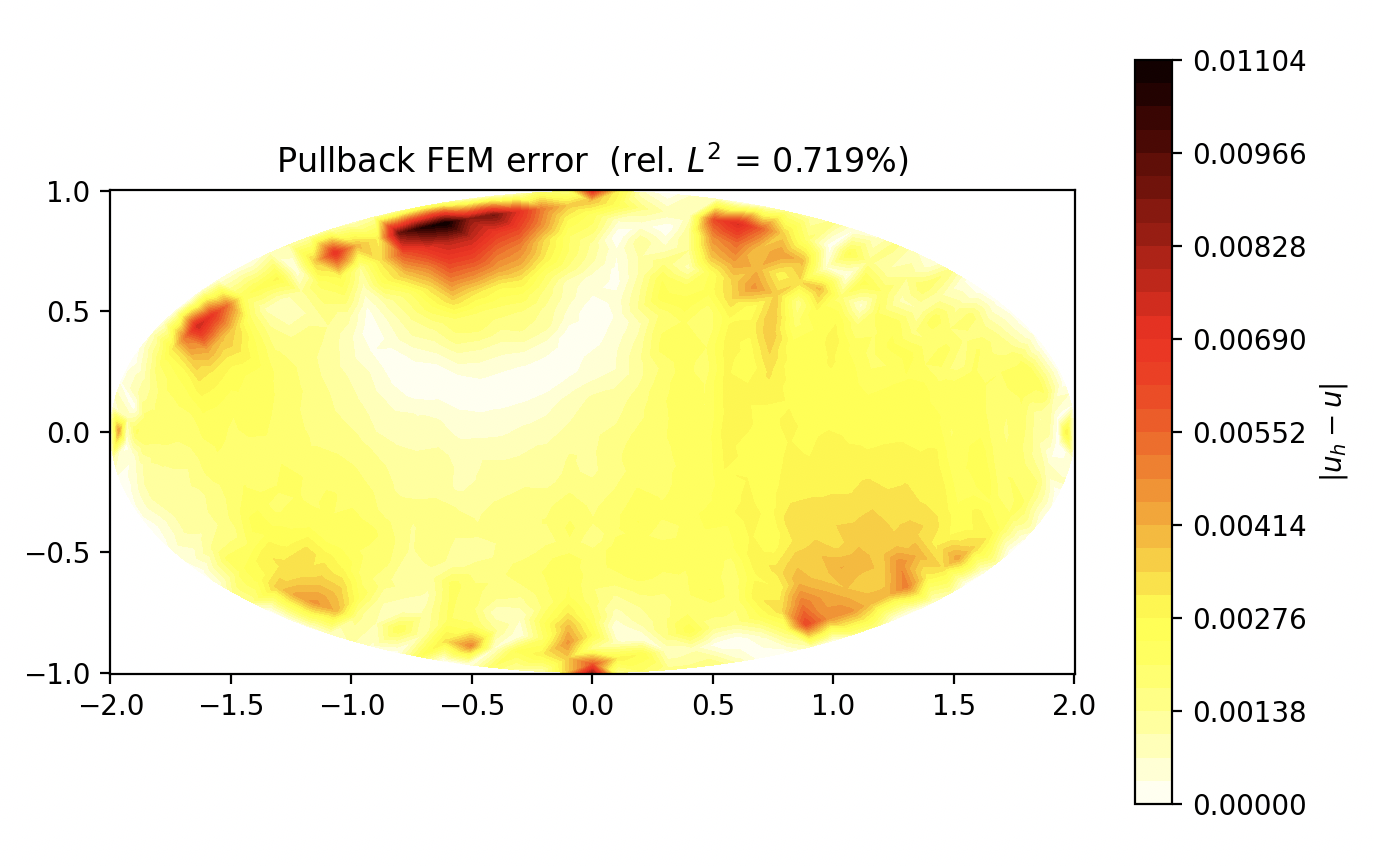}
        \caption{Pointwise absolute error map for the pullback FEM with Quad4 elements}
        \label{fig:ex1_mapped}
    \end{subfigure}

    \caption{Visual comparison of direct and pullback FEM solutions}
    \label{fig:Example1_visualComparison}
\end{figure}

Figure~\ref{fig:Example1_visualComparison} provides a visual comparison. The first row shows the exact solution on the ellipse (left) and the pointwise absolute error map for the direct FEM with Tri3 elements (right), while the second row depicts the same error map for the direct FEM with Quad4 (left) and for the pullback FEM on the reference square (right). 

\begin{table}[ht]
\centering
\caption{Comparison of direct and mapped FEM solutions for the Poisson problem on the ellipse ($a=2$, $b=1$). Errors are measured against the manufactured solution $u_{\mathrm{exact}} = 1 - x_1^2/a^2 - x_2^2/b^2$.}
\label{tab:poisson_results}
\smallskip
\begin{tabular}{lcccc}
\toprule
\textbf{Method} & \textbf{Elements} & \textbf{Nodes} & rel. $\boldsymbol{L^2}$ \textbf{error}\\
\midrule
Direct FEM (Tri3, Gmsh) & 1825 & 967 & $2.32 \times 10^{-3}$ \\
Direct FEM (Quad4, Gmsh) & 903 & 961 & $2.82 \times 10^{-3}$ \\
Mapped FEM (Quad4, $30\times 30$) & 900 & 961 & $7.19 \times 10^{-3}$ \\
\bottomrule
\end{tabular}
\end{table}

Overall, this example provided a demonstration of the validity of the proposed mathematical framework, as described in Sections~\ref{sec:ref_fem}--\ref{sec:nsf}, on a problem with a known analytical solution. We conclude that the NSF mapping is sufficiently accurate for FEM purposes, and, despite the slight boundary deviation ($<1\%$), the mapped FEM solution can be highly accurate. 

\subsection{Poisson problem on multiple elliptical and superelliptical domains}
\label{subsec:poisson_superelliptical}

The second numerical example extends the framework to the Poisson equation but now it is solved on a family of randomly generated domains. Unlike the previous example of Section~\ref{subsec:poisson_ellipse}, which validated the method on a single fixed geometry with a known analytical solution, this example tests the multi-shape conditioning pipeline of Section~\ref{subsec:two_stage}.

In this example the test domains are generated randomly as ellipses or superellipses. Writing
$\boldsymbol{c}=(c_1,c_2)$ for the centre, $a,b>0$ for the semi-axes, $\vartheta$ for the rotation
angle, $R(\vartheta)$ for the rotation matrix
\begin{equation}
R(\vartheta) = \begin{pmatrix} \cos\vartheta & -\sin\vartheta \\ \sin\vartheta & \cos\vartheta \end{pmatrix}
\end{equation}
then, the boundary of an ellipse is the curve
\begin{equation}\label{eq:ellipse}
\partial\Omega = \left\{\, \boldsymbol{c} + R(\vartheta)
\begin{pmatrix} a\cos t \\ b\sin t \end{pmatrix} \;:\; t\in[0,2\pi) \,\right\}
\end{equation}
A superellipse (Lam\'{e} curve) generalises this through an exponent $n$, with boundary
\begin{equation}\label{eq:superellipse}
\partial\Omega = \left\{\, \boldsymbol{c} + R(\vartheta)
\begin{pmatrix} a\,\operatorname{sgn}(\cos t)\,\lvert\cos t\rvert^{2/n} \\[2pt]
b\,\operatorname{sgn}(\sin t)\,\lvert\sin t\rvert^{2/n} \end{pmatrix} \;:\; t\in[0,2\pi) \,\right\}
\end{equation}
In the local, unrotated frame the superellipse satisfies the implicit relation
$\lvert x_1/a\rvert^{n} + \lvert x_2/b\rvert^{n} = 1$. The case $n=2$ recovers the ellipse
of Eq.~\eqref{eq:ellipse}, while larger $n$ yields increasingly rectangular shapes.

Each domain is obtained by sampling its parameters from the following distributions: for the
ellipses, $c_1,c_2\sim\mathcal{U}(-0.5,0.5)$, $a\sim\mathcal{U}(0.8,1.5)$,
$b\sim\mathcal{U}(0.4,a)$, and $\vartheta\sim\mathcal{U}(0,\pi)$; for the superellipses,
$c_1,c_2\sim\mathcal{U}(-0.3,0.3)$, $a\sim\mathcal{U}(0.8,1.4)$, $b\sim\mathcal{U}(0.6,1.4)$,
$n\sim\mathcal{U}(1.5,5)$, and $\vartheta\sim\mathcal{U}(0,\pi)$. 

\begin{figure}[H]
  \centering
  \begin{subfigure}[t]{0.49\textwidth}
    \includegraphics[width=\linewidth]{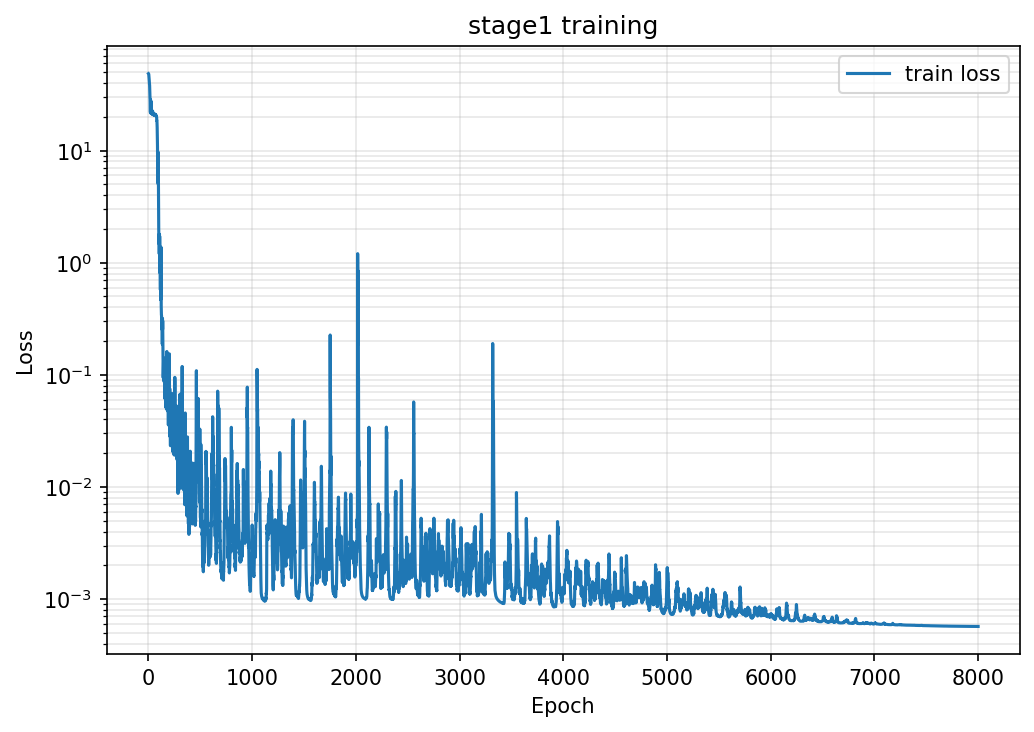}
    \caption{Stage~1: NSF auto-decoder training}
    \label{fig:ex2_loss1}
  \end{subfigure}\hfill
  \begin{subfigure}[t]{0.49\textwidth}
    \includegraphics[width=\linewidth]{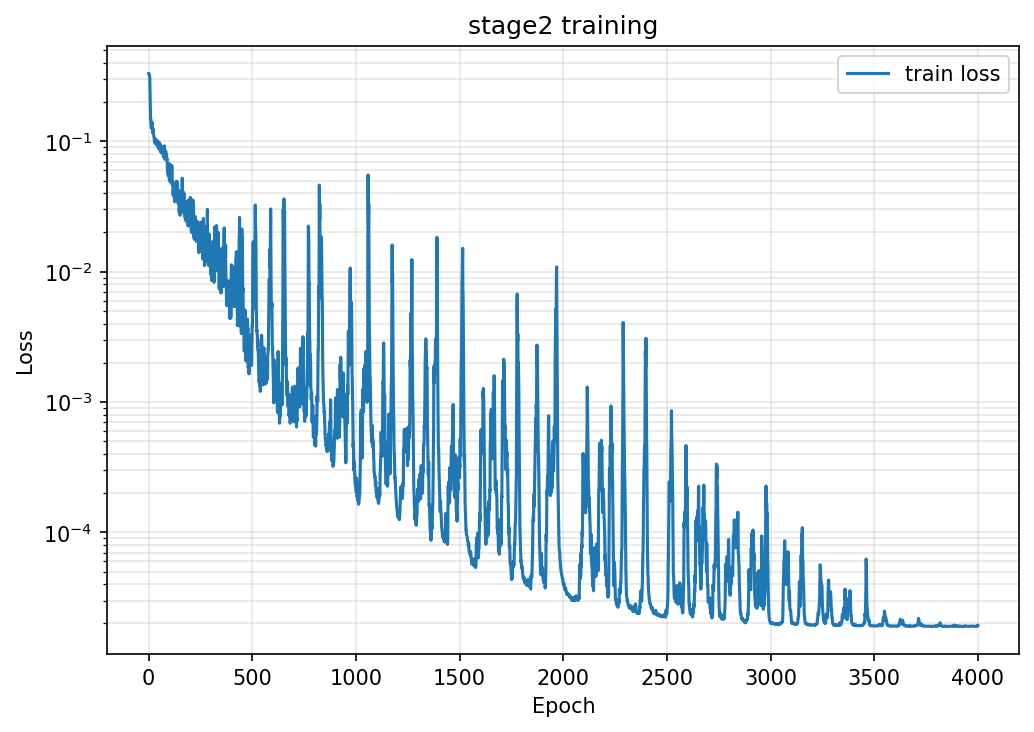}
    \caption{Stage~2: CNN encoder training}
    \label{fig:ex2_loss2}
  \end{subfigure}

  \caption{Evolution of the training loss for the two-stage geometry-learning framework:
(a) Stage~1, joint training of the NSF auto-decoder and latent representations;
(b) Stage~2, training of the CNN encoder to predict the learned latent representations
from the input SDF images.}
  \label{fig:training_losses}
\end{figure}

Based on these parameter choices, we generated 2000 random ellipses and superellipses to train the proposed 2-stage architecture. The hyperparameters chosen for the networks are reported in Tables 
\ref{tab:ex2_stageA} and \ref{tab:ex2_stageB}. Also, Figures~\ref{fig:ex2_loss1} and \ref{fig:ex2_loss2} show the minimization of the training loss for the first and second stage of the proposed architecture.

\begin{table}[htbp]
\centering

\begin{minipage}{0.45\textwidth}
\centering
\caption{Hyperparameters of the NSF model in Stage 1}
\label{tab:ex2_stageA}
\begin{tabular}{cc}
\hline
{} & NSF parameters \\
\hline
number of layers & 18 \\
number of bins & 20 \\
hidden dimension & 128 \\
number of hidden layers & 3 \\
training epochs & 8000 \\
optimizer & ADAM\\
learning rate & 5e-4 \\
\hline
\end{tabular}
\end{minipage}
\hfill
\begin{minipage}{0.45\textwidth}
\centering
\caption{Hyperparameters of the Encoder in Stage 2}
\label{tab:ex2_stageB}
\begin{tabular}{cc}
\hline
{} & CNN parameters \\
\hline
input dimension & $64\times 64$ \\
latent code/dimension & 8 \\
convolutional layers & $1\to 32\to 64\to 128$ \\
convolutional kernel sizes & $4\times 4$, $4\times 4$, $4\times 4$ \\
activation functions & GelU \\
training epochs & 4000 \\
optimizer & ADAM\\
learning rate & 1e-3 \\
\hline
\end{tabular}
\end{minipage}

\end{table}

\begin{figure}[H]
    \centering
    \includegraphics[width=\textwidth]{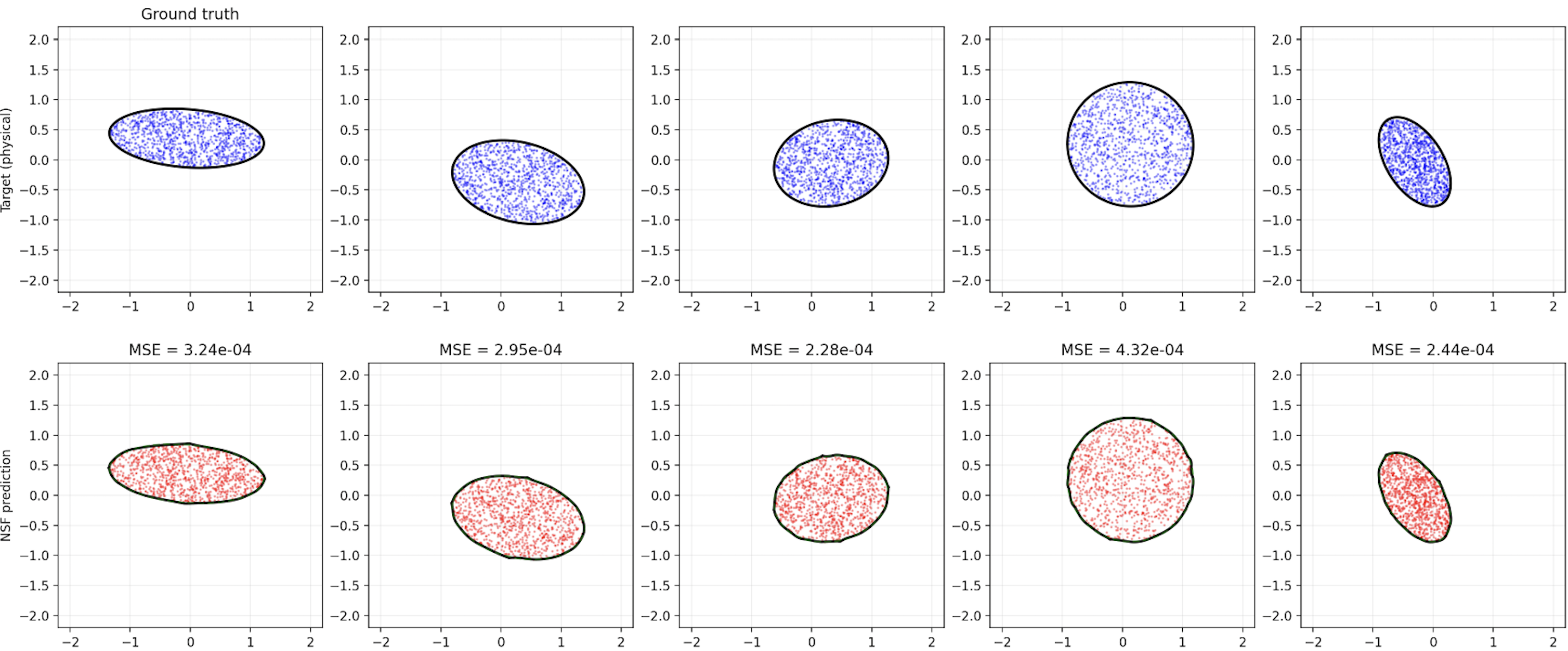}
    \caption{The first row depicts five representative elliptical domains from the test set, generated via Eq.~\eqref{eq:ellipse}. The second row represents their mapping from the reference square via the NSF diffeomorphism}
    \label{fig:ellipse_gallery}
\end{figure}

\begin{figure}[H]
    \centering
    \includegraphics[width=\textwidth]{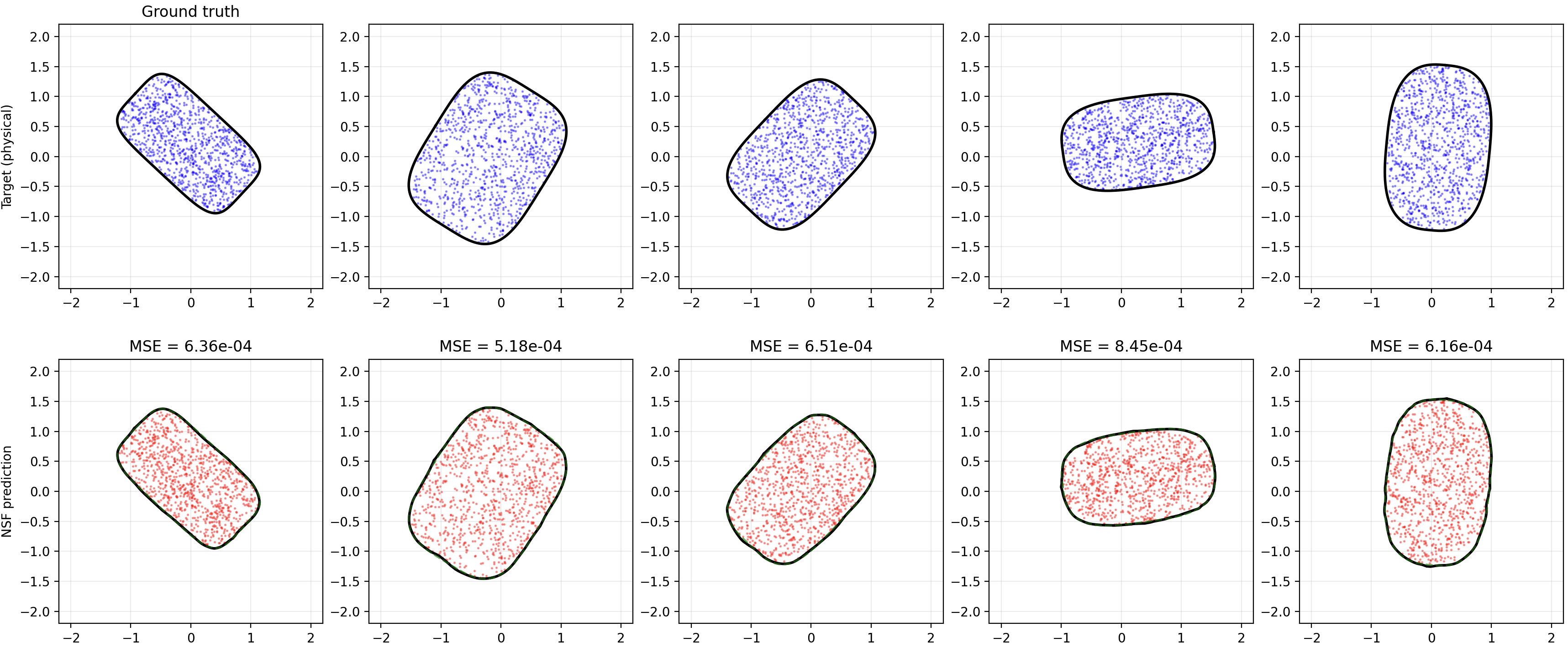}
    \caption{The first row depicts five representative superelliptical domains from the test set, generated via Eq.~\eqref{eq:superellipse}. The second row represents their mapping from the reference square via the NSF diffeomorphism}
    \label{fig:superellipse_gallery}
\end{figure}

Figures~\ref{fig:ellipse_gallery} and \ref{fig:superellipse_gallery} show five representative ellipses and superellipses from the test dataset, respectively, and their corresponding mappings from the reference square. 

\begin{figure}[H]
  \centering
  \begin{subfigure}[t]{0.49\textwidth}
    \includegraphics[width=\linewidth]{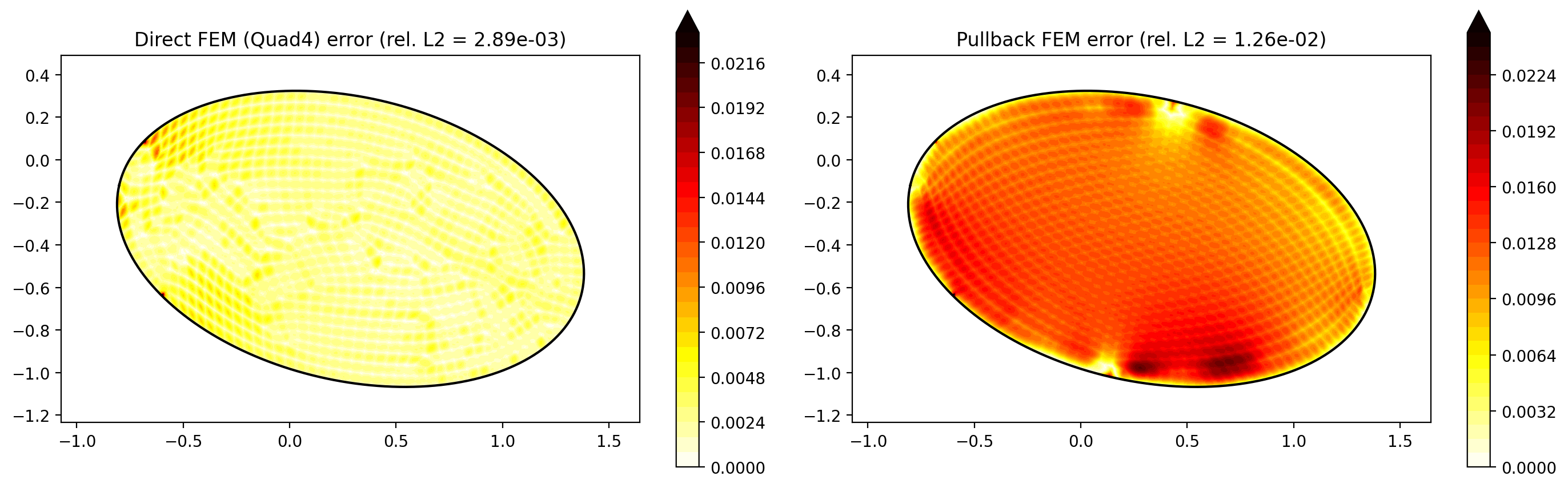}
    \caption{test case 1}
    \label{fig:err0}
  \end{subfigure}\hfill
  \begin{subfigure}[t]{0.49\textwidth}
    \includegraphics[width=\linewidth]{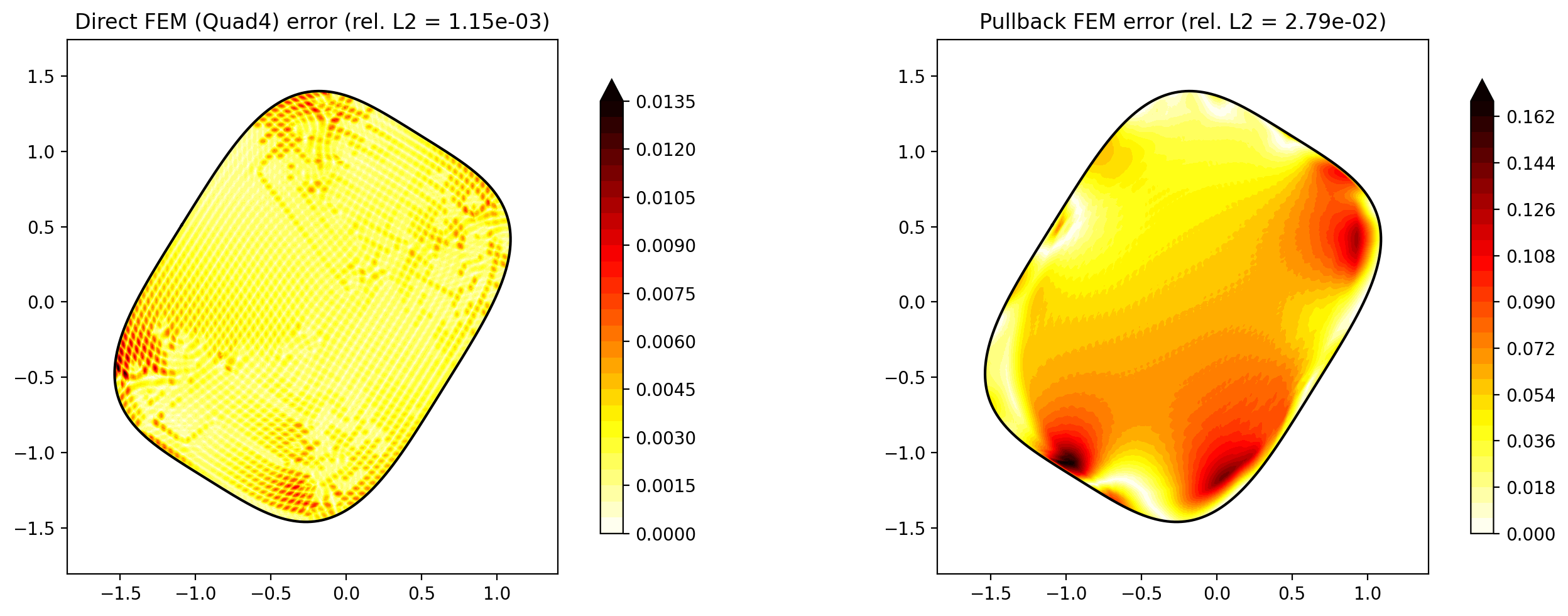}
    \caption{test case 2}
    \label{fig:err1}
  \end{subfigure}

  \vspace{1ex}

  \begin{subfigure}[t]{0.49\textwidth}
    \includegraphics[width=\linewidth]{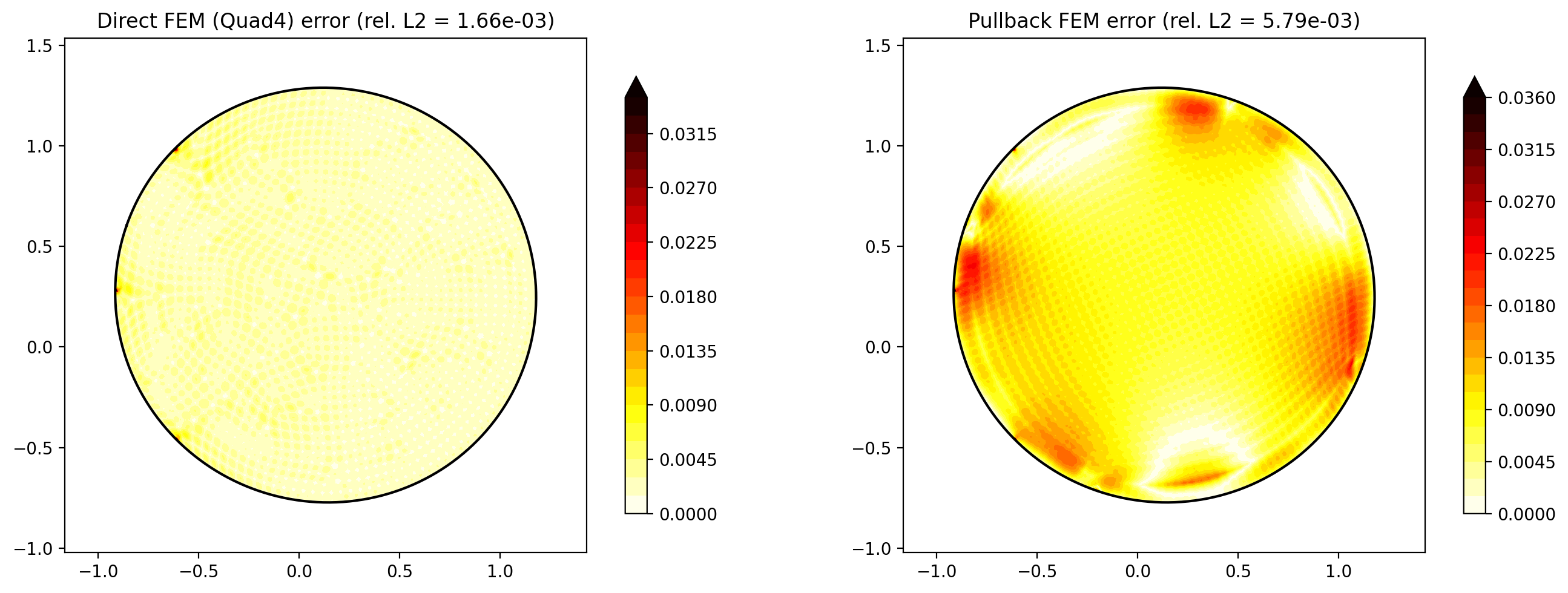}
    \caption{test case 3}
    \label{fig:err2}
  \end{subfigure}\hfill
  \begin{subfigure}[t]{0.49\textwidth}
    \includegraphics[width=\linewidth]{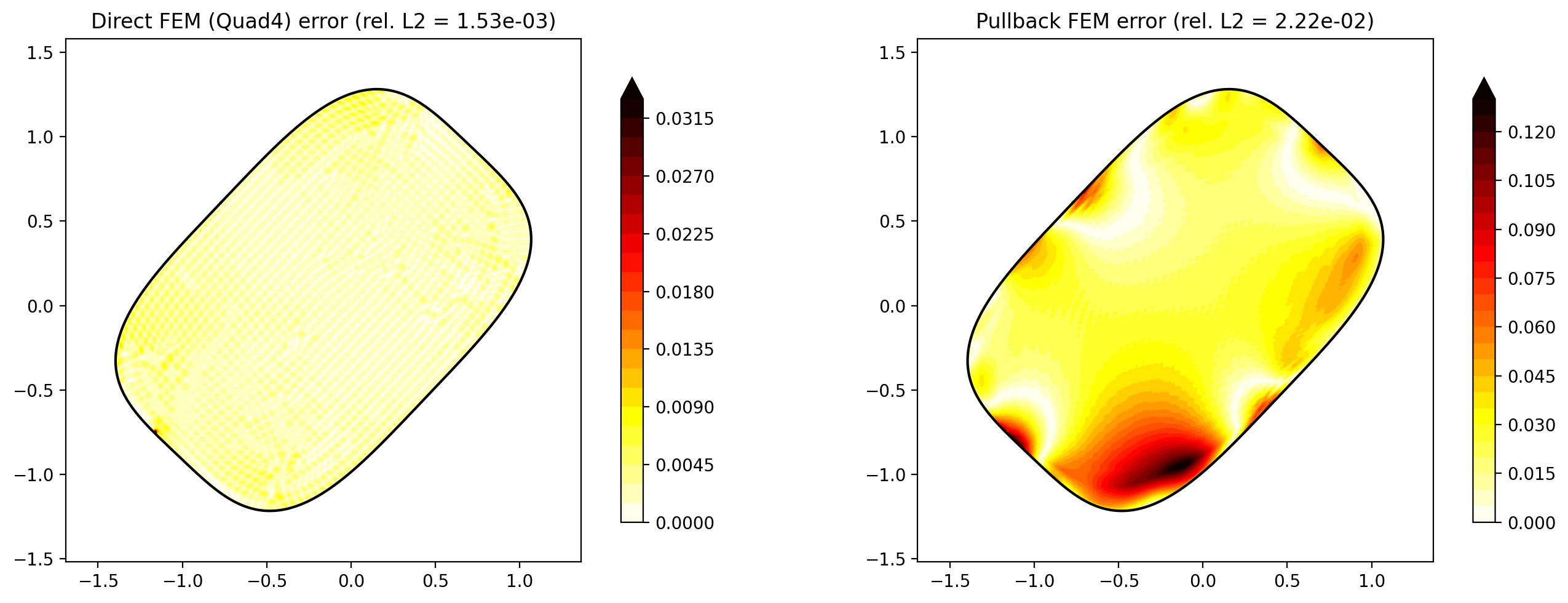}
    \caption{test case 4}
    \label{fig:err3}
  \end{subfigure}
  
    \vspace{1ex}

  \begin{subfigure}[t]{0.49\textwidth}
    \includegraphics[width=\linewidth]{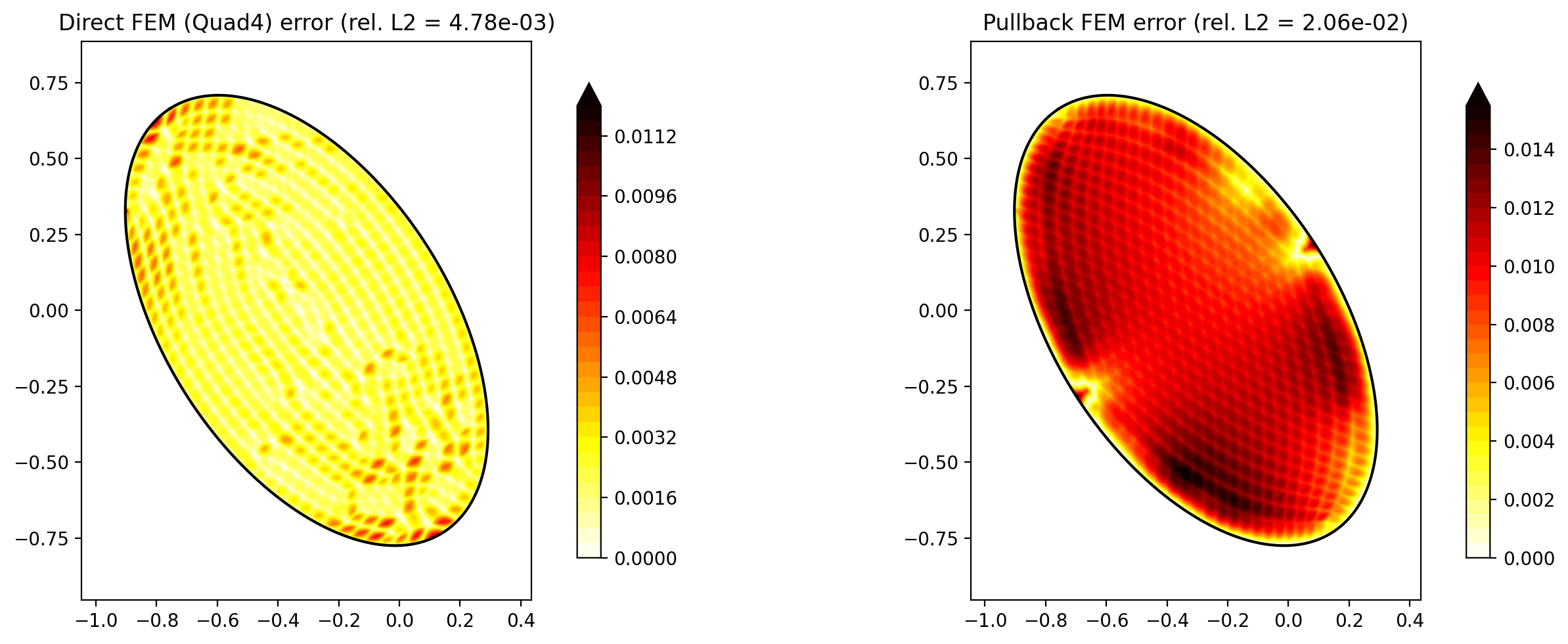}
    \caption{test case 5}
    \label{fig:err4}
  \end{subfigure}\hfill
  \begin{subfigure}[t]{0.49\textwidth}
    \includegraphics[width=\linewidth]{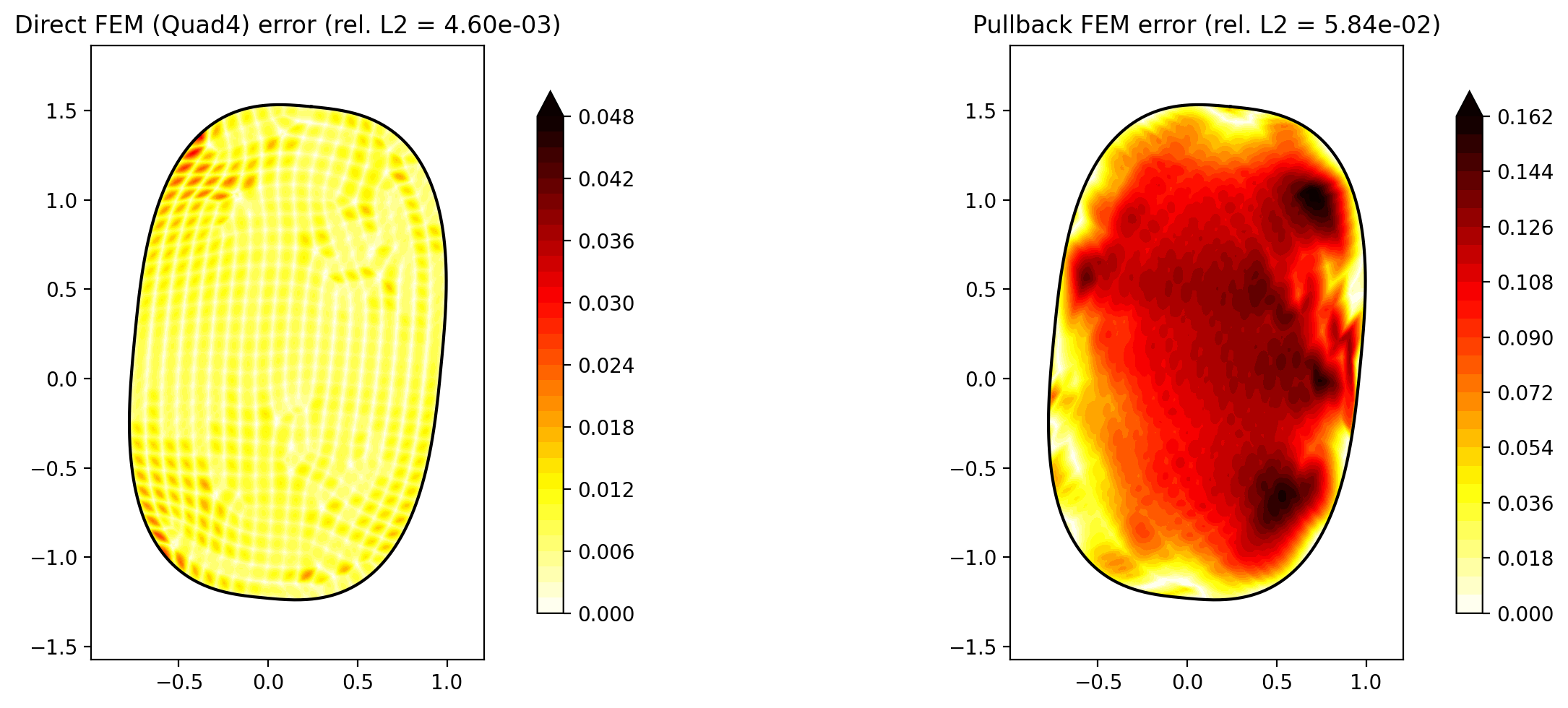}
    \caption{test case 6}
    \label{fig:err5}
  \end{subfigure}

  \caption{Pointwise absolute error maps for six test domains: direct FEM (left of each pair)
  and pullback FEM (right), each measured against the high-fidelity reference.}
  \label{fig:multishape_errors}
\end{figure}

Having established the accuracy of the learned geometric representation, we next assess the performance of the proposed approach in the FE simulations.  Fig.~\ref{fig:multishape_errors} presents the pointwise error fields for six test geometries randomly selected from the test dataset, comparing the direct FEM and the pullback FEM solutions against a high-fidelity reference solution obtained using a highly refined finite element mesh of approximately 40,000 elements in each case. The corresponding relative $L^2$ errors are reported in Table~\ref{tab:multishapePoisson_results}. For a fair comparison, the direct and pullback discretisations were selected to have approximately comparable numbers of elements and degrees of freedom.

Upon inspection of the results, we observe that in this multi-shape setting the accuracy of the pullback FEM approach is almost one order of magnitude below the direct FEM approach. This behaviour differs from the single-shape case in the previous example, where the learned mapping reproduced the target geometry with very high accuracy. In the present example, the diffeomorphism for each unseen geometry is inferred through the two-stage, multi-shape architecture from a model trained across an entire family of domains. This constitutes a considerably more challenging approximation problem, and the resulting mapping errors propagate through the transformed geometric coefficients of the pullback formulation, contributing to the observed increase in the finite element error.

Nevertheless, the relative error of the pullback solutions, being less than $6\%$ in all cases, is still considered to be acceptable, demonstrating the feasibility of the proposed multi-shape framework, while also indicating that the accuracy of the learned diffeomorphism becomes an important factor in the overall solution accuracy. Further improvements may be expected from larger and more diverse training datasets, and improved network and training hyperparameters.

\begin{table}[ht]
\centering
\caption{Comparison of direct and mapped FEM solutions for the Poisson problem on six domains}
\label{tab:multishapePoisson_results}
\smallskip
\begin{tabular}{lcccccc}
\toprule
& \multicolumn{3}{c}{\textbf{Direct FEM (Quad4, Gmsh)}} & \multicolumn{3}{c}{\textbf{Mapped FEM (Quad4)}} \\
\cmidrule(lr){2-4} \cmidrule(lr){5-7}
\textbf{Shape} & Elements & Nodes & rel. L2 error & Elements & Nodes & rel. L2 error \\
\midrule
1 & 989 & 1047 & $2.89\times 10^{-3}$ & 900 (30$\times$30) & 961 & $1.26\times 10^{-2}$ \\
2 & 1274 & 1344 & $1.15\times 10^{-3}$ & 1225 (35$\times$35) & 1296 & $2.79\times 10^{-2}$ \\
3 & 1442 & 1509 & $1.66\times 10^{-3}$ & 1444 (38$\times$38) & 1521 & $5.79\times 10^{-3}$ \\
4 & 1597 & 1672 & $1.53\times 10^{-3}$ & 1600 (40$\times$40) & 1681 & $2.22\times 10^{-2}$ \\
5 & 518 & 561 & $4.78\times 10^{-3}$ & 529 (23$\times$23) & 576 & $2.06\times 10^{-2}$ \\
6 & 1718 & 1795 & $4.60\times 10^{-3}$ & 1681 (41$\times$41) & 1764 & $5.84\times 10^{-2}$ \\
\bottomrule
\end{tabular}
\end{table}

\subsection{Plane-stress elasticity on Fourier blob domains}
\label{subsec:elasticity_blobs}

The third numerical example applies the framework to a vector-valued PDE, that of linear elasticity under plane-stress conditions, which is solved on a family of randomly generated domains. Here, the train and test domains are generated as random smooth perturbations of circles using a
truncated Fourier series in polar coordinates, also called Fourier blobs. Specifically, the radial function is defined as
\begin{equation}\label{eq:fourier_blob}
r(\vartheta) = r_0 + \sum_{m=1}^{M} \left[ a_m \cos(m\vartheta) + b_m \sin(m\vartheta) \right], \qquad \vartheta \in [0, 2\pi)
\end{equation}

\noindent where $r_0 > 0$ is the base radius and $M$ is the number of Fourier modes. The corresponding boundary is parameterized as
\begin{equation}
\partial\Omega
=
\left\{
\boldsymbol{c}
+
\begin{pmatrix}
r(\vartheta)\cos\vartheta\\
r(\vartheta)\sin\vartheta
\end{pmatrix}
:
\vartheta\in[0,2\pi)
\right\}
\label{eq:fourier_blob_curve}
\end{equation}
where $\boldsymbol{c}=(c_1,c_2)$ denotes the centre of the domain.

In the present study, $M=5$ modes are retained and the geometric parameters are sampled independently according to

\begin{equation}
c_1,c_2 \sim \mathcal{U}(-0.3,0.3),
\qquad
r_0 \sim \mathcal{U}(0.8,1.2)
\end{equation}

\noindent while the Fourier coefficients are sampled as

\begin{equation}
a_m,b_m
\sim
\mathcal{U}
\left(
-\frac{0.3}{m},
\frac{0.3}{m}
\right),
\qquad
m=1,\ldots,M
\end{equation}

The auto-decoder (Stage~1) is trained on 4000 randomly generated Fourier blob shapes using the configuration described in Section~\ref{subsec:two_stage}, with the harmonic map computed on a $256 \times 256$ finite-difference grid for each training shape. The CNN encoder (Stage~2) is subsequently trained on the SDF images and stored latent codes. The hyperparameters chosen for the networks of this example are reported in Tables  \ref{tab:ex3_stageA} and \ref{tab:ex3_stageB}. 

\begin{table}[htbp]
\centering

\begin{minipage}{0.45\textwidth}
\centering
\caption{Hyperparameters of the NSF model in Stage 1}
\label{tab:ex3_stageA}
\begin{tabular}{cc}
\hline
{} & NSF parameters \\
\hline
number of layers & 20 \\
number of bins & 20 \\
hidden dimension & 128 \\
number of hidden layers & 3 \\
training epochs & 8000 \\
optimizer & ADAM\\
learning rate & 5e-4 \\
\hline
\end{tabular}
\end{minipage}
\hfill
\begin{minipage}{0.45\textwidth}
\centering
\caption{Hyperparameters of the Encoder in Stage 2}
\label{tab:ex3_stageB}
\begin{tabular}{cc}
\hline
{} & CNN parameters \\
\hline
input dimension & $64\times 64$ \\
latent code/dimension & 8 \\
convolutional layers & $1\to 32\to 64\to 128$ \\
convolutional kernel sizes & $4\times 4$, $4\times 4$, $4\times 4$ \\
activation functions & GelU \\
training epochs & 4000 \\
optimizer & ADAM\\
learning rate & 1e-3 \\
\hline
\end{tabular}
\end{minipage}

\end{table}

\begin{figure}[H]
    \centering
    \includegraphics[width=\textwidth]{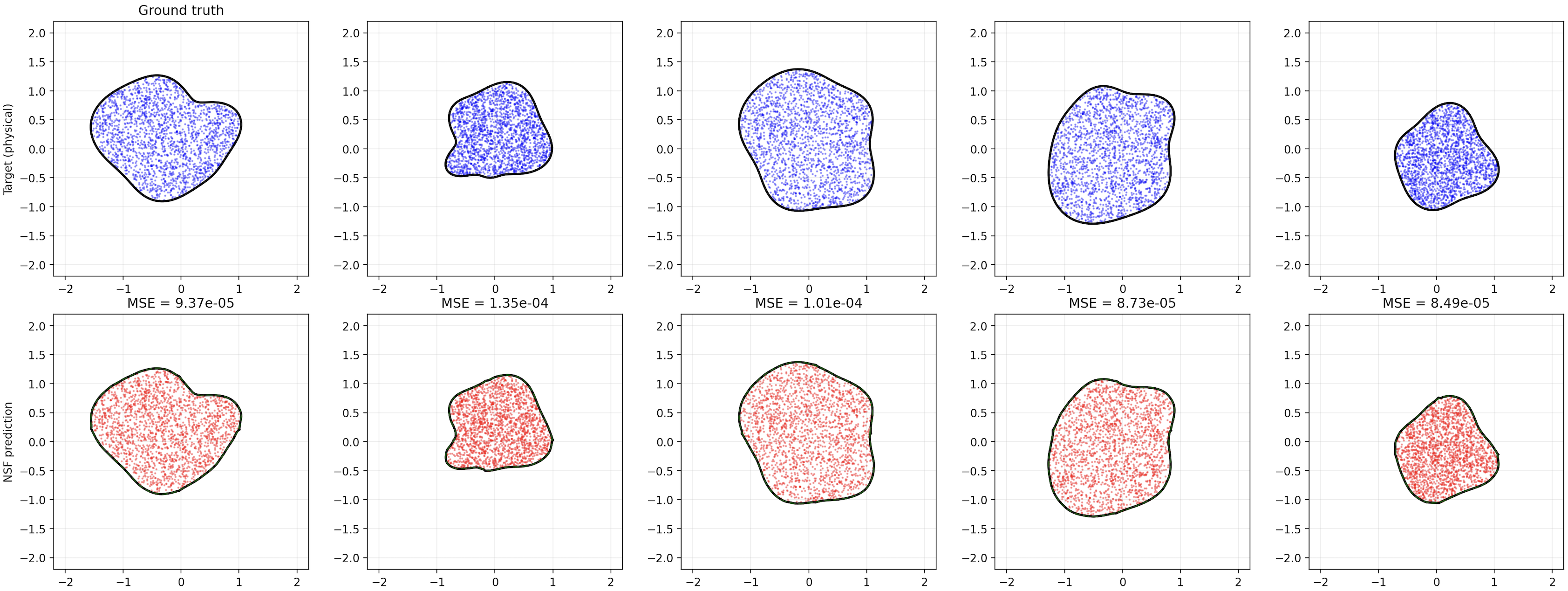}
    \caption{The first row depicts five representative Fourier blob domains from the test set, generated via Eq.~\eqref{eq:fourier_blob}. The second row represents their mapping from the reference square via the NSF diffeomorphism}
    \label{fig:blob_gallery}
\end{figure}

Figure~\ref{fig:blob_gallery} shows five representative Fourier blobs from the test set, illustrating the geometric diversity of the shape family. The top row depicts the target physical domains together with the sampled points, and the bottom row the images of the corresponding reference-square points under the learned map $G_{\boldsymbol{\theta}}(\boldsymbol{\xi}_k;\boldsymbol{z}^\ast)$.

In this setting, we consider the equilibrium of a linear elastic body under plane-stress conditions. The strong form reads: find the displacement field $\mathbf{u}: \Omega \to \mathbb{R}^2$ such that
\begin{equation}\label{eq:elasticity_strong}
-\nabla \cdot \boldsymbol{\sigma}(\mathbf{u}) = \mathbf{f} \quad \text{in } \Omega, \qquad
\mathbf{u} = \mathbf{0} \quad \text{on } \Gamma_D, \qquad
\boldsymbol{\sigma}(\mathbf{u}) \cdot \mathbf{n} = \mathbf{t} \quad \text{on } \Gamma_N
\end{equation}
where $\boldsymbol{\sigma} = \mathbf{C} : \boldsymbol{\varepsilon}(\mathbf{u})$ is the Cauchy stress tensor, $\boldsymbol{\varepsilon}(\mathbf{u}) = \tfrac{1}{2}(\nabla\mathbf{u} + \nabla\mathbf{u}^T)$ is the symmetric strain tensor, $\mathbf{C}$ is the plane-stress elasticity tensor with Young's modulus $E = 210 \times 10^3 \, \mathrm{MPa}$ and Poisson's ratio $\nu = 0.3$, $\mathbf{f}=(0, -10)^T$ is the body force, and $\mathbf{t}=(0, 0)^T$ is the prescribed traction on the Neumann boundary $\Gamma_N$. The Dirichlet boundary is defined geometrically as the set of boundary nodes whose $x_1$-coordinate lies within the leftmost $20\%$ of the shape's horizontal extent:
\begin{equation}\label{eq:clamp_criterion}
\Gamma_D = \left\{ \boldsymbol{x} \in \partial\Omega \,:\, x_1 < x_1^{\min} + 0.2 \, (x_1^{\max} - x_1^{\min}) \right\}
\end{equation}
where $x_1^{\min}$ and $x_1^{\max}$ are the extremal $x_1$-coordinates of the boundary. 

Under the change of variables $\boldsymbol{x}=G_\theta(\boldsymbol{\xi};\boldsymbol{z})$, with Jacobian $\boldsymbol{J}=DG_\theta$ and $\mathcal{J}=\det\boldsymbol{J}>0$, the physical-domain gradients of the finite element shape functions follow from the composition of the isoparametric map $\boldsymbol{\eta}\mapsto \boldsymbol{\xi}$, with Jacobian $\boldsymbol{J}_{iso}=\frac{\partial\boldsymbol{\xi}}{\partial\boldsymbol{\eta}}$, and the learned map $\boldsymbol{\xi}\mapsto \boldsymbol{x}$:

\begin{equation}
\nabla_{\boldsymbol{x}} N_i
=
\nabla_{\boldsymbol{\eta}} N_i\,\boldsymbol{J}_{iso}^{-1}\boldsymbol{J}^{-1}
\end{equation}

\noindent where $\boldsymbol{\eta}$ denotes the coordinates on the parent element, while for the uniform Quad4 elements of the structured reference mesh, $\boldsymbol{J}_{iso}$ is a constant diagonal scaling.

For elasticity, the strain--displacement matrix $\boldsymbol{B}_q$ is therefore constructed in the usual Voigt form using these transformed physical-domain derivatives, i.e. $\partial N_i/\partial x_1$ and $\partial N_i/\partial x_2$, evaluated at each quadrature point. Thus, the element stiffness matrix and load vector on the reference square are
\begin{equation}\label{eq:ke_elasticity}
\hat{\mathbf{k}}_e = \sum_q w_q \, \mathbf{B}_q^T \, \mathbf{C} \, \mathbf{B}_q \, \mathcal{J}_q \, |\det \boldsymbol{J}_{\mathrm{iso},q}|
\end{equation}
\begin{equation}\label{eq:fe_elasticity}
\hat{\mathbf{f}}_e = \sum_q w_q \, \boldsymbol{N}^T(\boldsymbol{\xi}_q) \, \mathbf{f}\!\left(G_{\boldsymbol{\theta}}(\boldsymbol{\xi}_q; \mathbf{z})\right) \, \mathcal{J}_q \, |\det \boldsymbol{J}_{\mathrm{iso},q}|
\end{equation}
where $\boldsymbol{J}_{\mathrm{iso}}$ is the isoparametric Jacobian of the reference-element mapping. 

\begin{figure}[H]
    \centering

    \begin{subfigure}[b]{0.48\textwidth}
        \centering
        \includegraphics[width=0.8\linewidth]{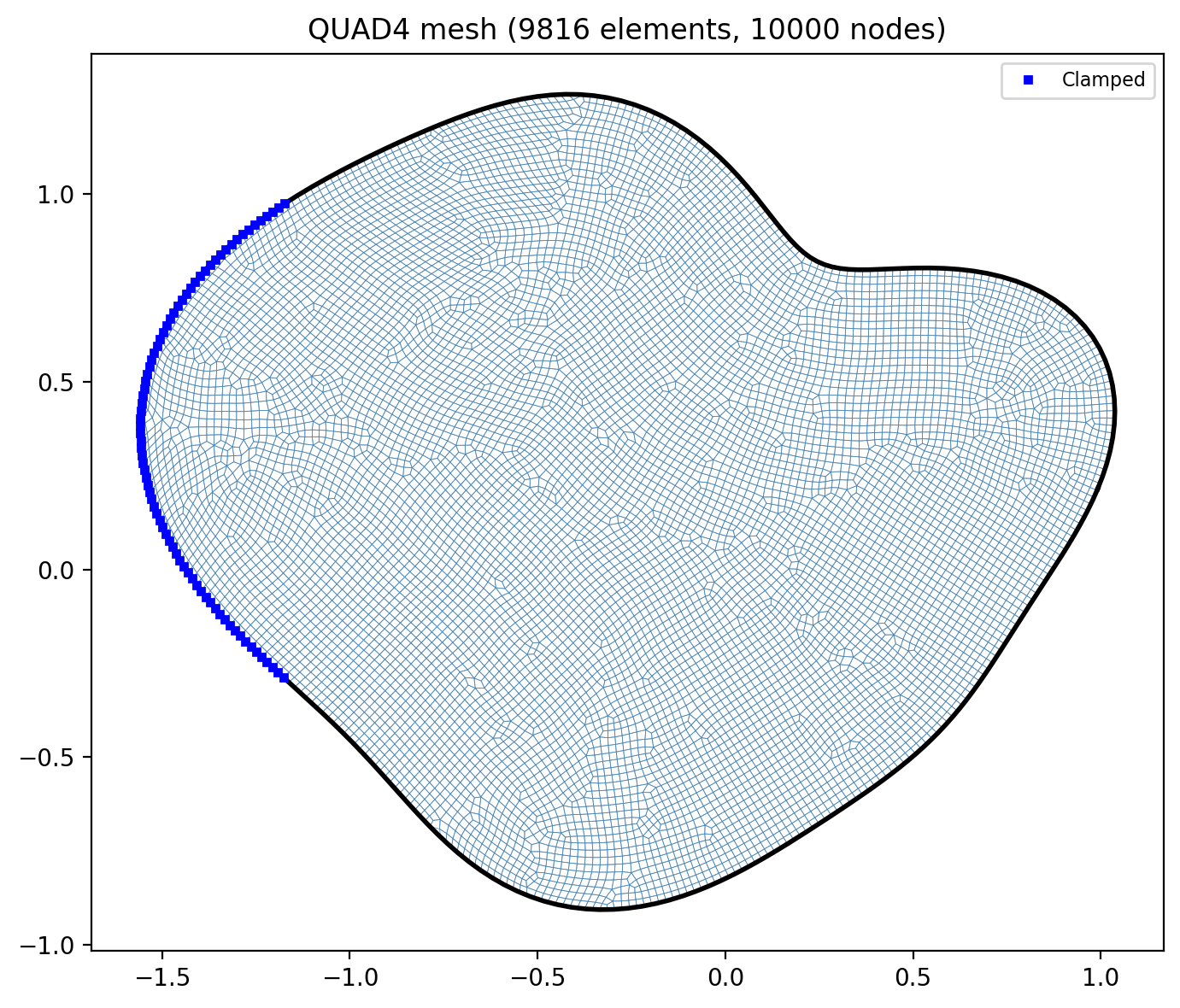}
        \caption{Mesh and boundary conditions on an arbitrary shape}
        \label{fig:ex3_shape1_discretization}
    \end{subfigure}
    \hfill
    \begin{subfigure}[b]{0.48\textwidth}
        \centering
        \includegraphics[width=\linewidth]{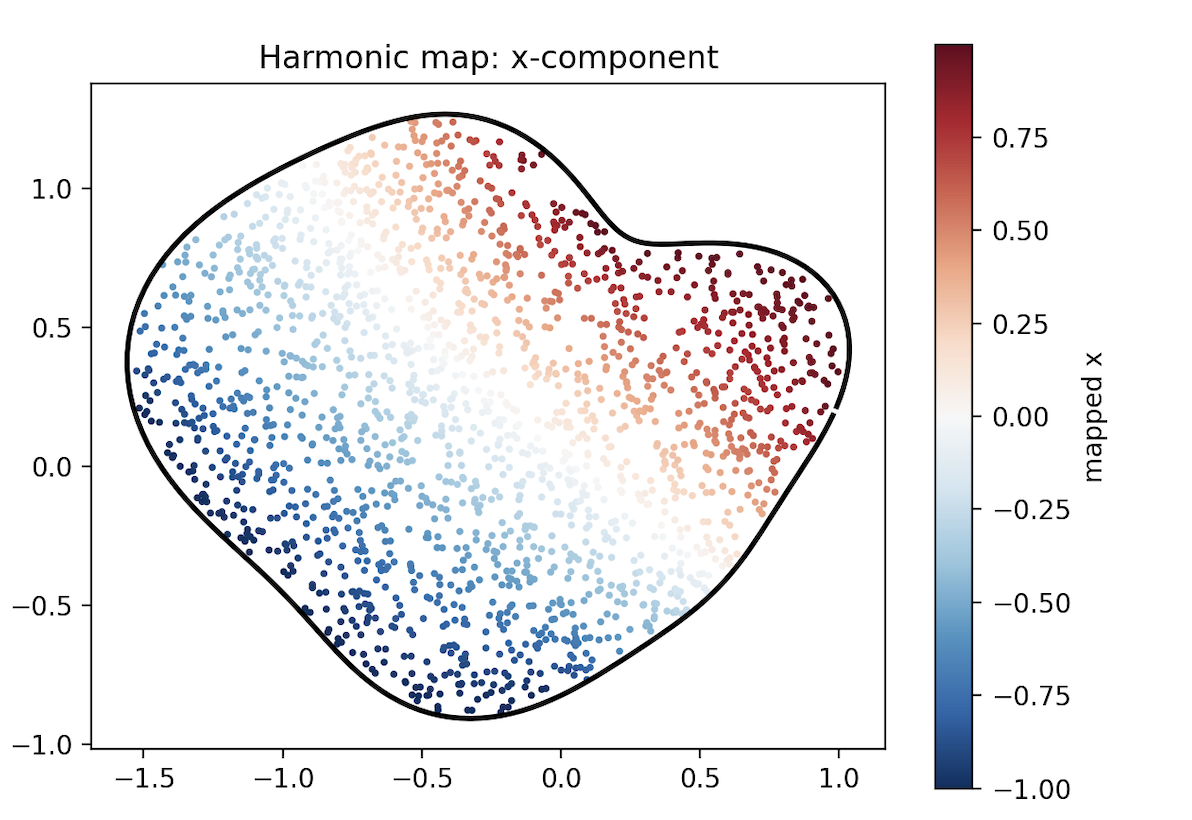}
        \caption{$x$-displacement of the internal points under the harmonic maps}
        \label{fig:ex3_HarmonicMaps_x}
    \end{subfigure}

    \vspace{0.5cm}

    \begin{subfigure}[b]{0.48\textwidth}
        \centering
        \includegraphics[width=\linewidth]{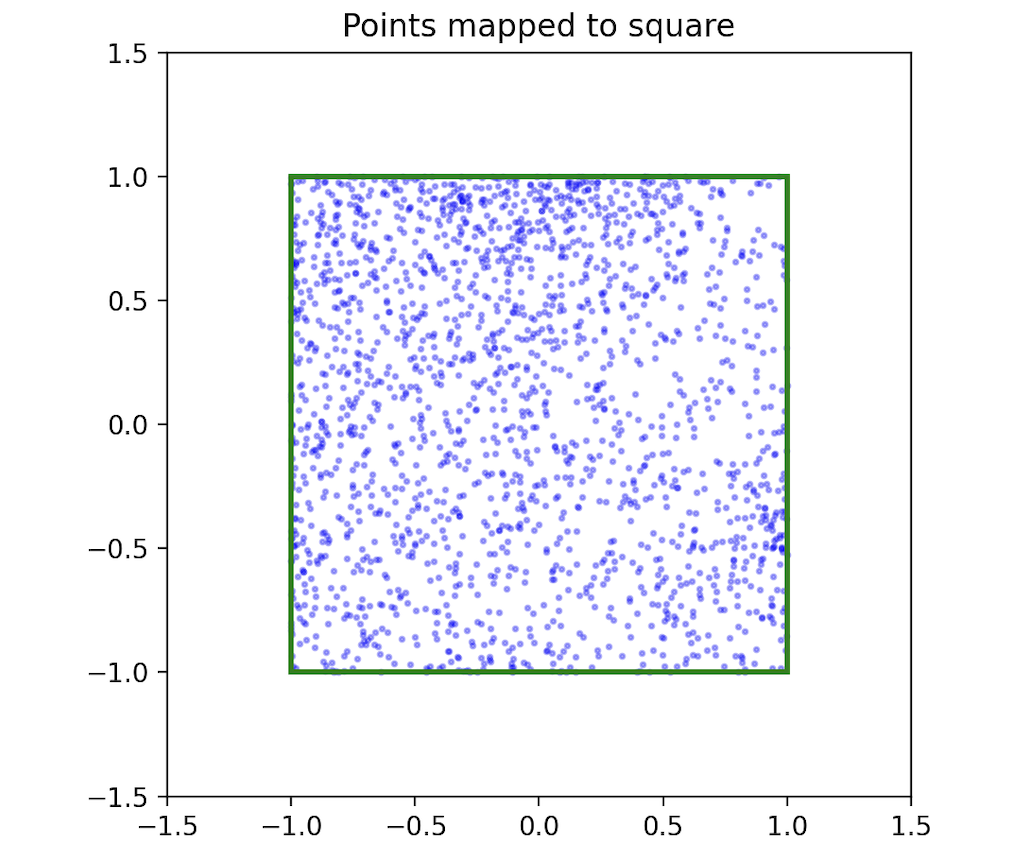}
        \caption{End locations of the internal points inside the reference square}
        \label{fig:ex3_PointsOnTheSquare}
    \end{subfigure}
    \hfill
    \begin{subfigure}[b]{0.48\textwidth}
        \centering
        \includegraphics[width=\linewidth]{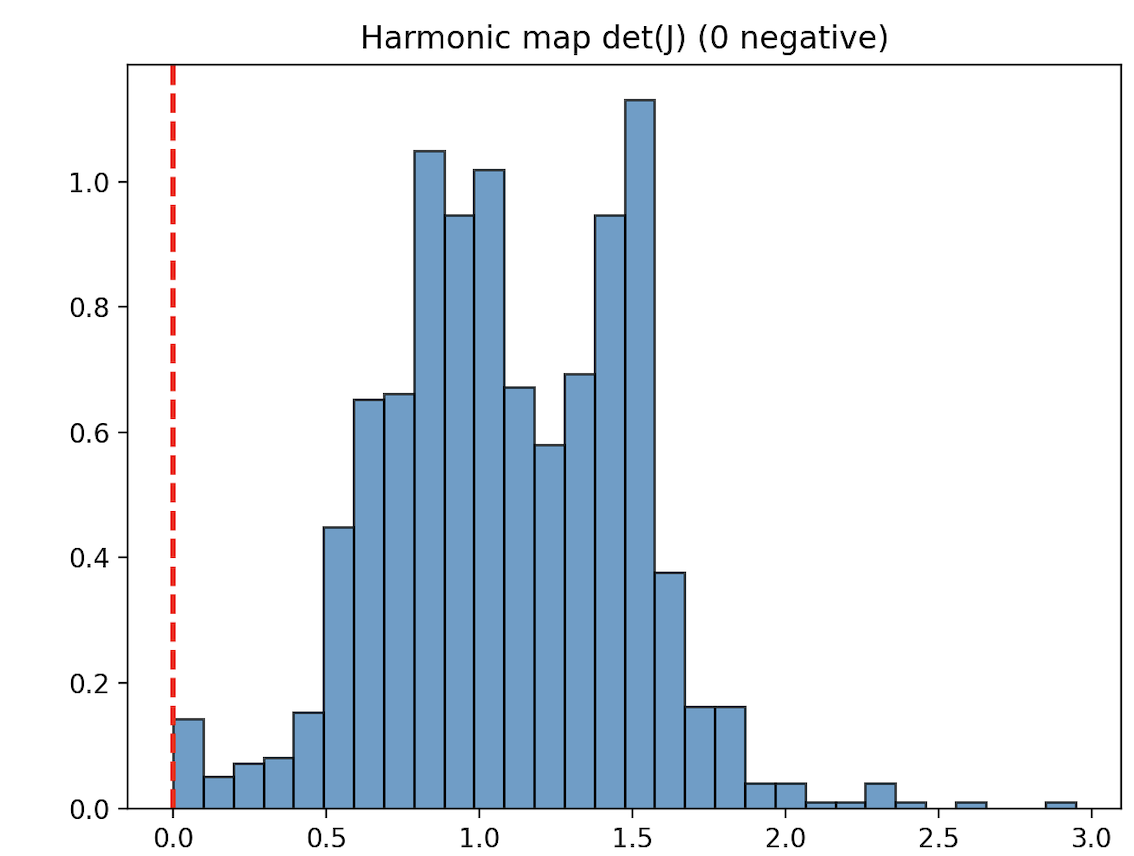}
        \caption{Distribution of Jacobian determinant values}
        \label{fig:ex3_JacobianDistribution}
    \end{subfigure}

    \caption{Mapping of a Fourier blob to the reference square}
    \label{fig:2x2}
\end{figure}

For the direct FEM baseline, each Fourier blob is meshed independently using Gmsh~\cite{Geuzaine2009} with linear quadrilateral elements (Quad4) and a characteristic element size $h = 0.02$. Figure~\ref{fig:ex3_shape1_discretization} depicts the mesh and boundary conditions for a representative Fourier-blob geometry. Figure~\ref{fig:ex3_HarmonicMaps_x} shows the $x$-displacement of the interior points under the harmonic mapping, while figure~\ref{fig:ex3_PointsOnTheSquare} shows their final locations in the reference square. The corresponding distribution of the Jacobian determinant is shown in figure~\ref{fig:ex3_JacobianDistribution}.

For each test shape $\Omega^*$, the inference pipeline (Algorithm~\ref{alg:inference}) proceeds as follows: the SDF image is computed on a $64 \times 64$ grid, the encoder predicts the latent code $\boldsymbol{z}^*$ in a single forward pass, and the geometric quantities ($\boldsymbol{J}_q^{-1}$, $\mathcal{J}_q$, $\boldsymbol{x}_q$) are precomputed at all quadrature points of the reference mesh. The stiffness matrix and load vector are then assembled using Eqs.~\eqref{eq:ke_elasticity}--\eqref{eq:fe_elasticity}, and the resulting linear system is solved by a sparse direct solver.

We present results for five test geometries not seen during training. For each geometry, the pullback FEM solution, computed on the fixed structured reference mesh, is compared with a direct FEM solution obtained on an independently generated Gmsh mesh. A high-fidelity reference solution is computed using a substantially refined discretisation of approximately 40,000 elements for each geometry. The accuracy of both approaches is assessed through the relative $L^2$ error of the displacement magnitude $|\boldsymbol{u}|$ with respect to this reference solution. The results are summarised in
Table~\ref{tab:elasticity_results}.

\begin{table}[ht]
\centering
\caption{Comparison of direct and mapped FEM solutions for plane-stress elasticity on five Fourier blob domains.}
\label{tab:elasticity_results}
\smallskip
\begin{tabular}{lcccccc}
\toprule
& \multicolumn{3}{c}{\textbf{Direct FEM (Quad4, Gmsh)}} & \multicolumn{3}{c}{\textbf{Mapped FEM (Quad4, $80\times 80$)}} \\
\cmidrule(lr){2-4} \cmidrule(lr){5-7}
\textbf{Shape} & Elements & Nodes & rel. L2 error & Elements & Nodes & rel. L2 error \\
\midrule
1 & 9816 & 10000 & 3.87$\times 10^{-3}$ & 6400 & 6561 & 2.07$\times 10^{-2}$ \\
2 & 6088 & 6234 & 1.74$\times 10^{-2}$ & 6400 & 6561 & 2.35$\times 10^{-2}$ \\
3 & 11533 & 11727 & 4.20$\times 10^{-3}$ & 6400 & 6561 & 5.05$\times 10^{-2}$ \\
4 & 10510 & 10697 & 3.14$\times 10^{-3}$ & 6400 & 6561 & 3.64$\times 10^{-3}$ \\
5 & 5852 & 5992 & 7.46$\times 10^{-3}$ & 6400 & 6561 & 7.99$\times 10^{-3}$ \\
\bottomrule
\end{tabular}
\end{table}

The pullback FEM errors remain below approximately 5.1\% for all five geometries. For several cases, the accuracy is comparable to that of the direct FEM, whereas Shapes 1 and 3 exhibit greater sensitivity to the learned geometric mapping. As in the previous examples, mapping inaccuracies or distortions propagate through the transformed gradients and integration measure, thereby contributing to the overall solution error.

\begin{figure}[H]
    \centering

    \begin{subfigure}[b]{0.68\textwidth}
        \centering
        \includegraphics[width=\textwidth]{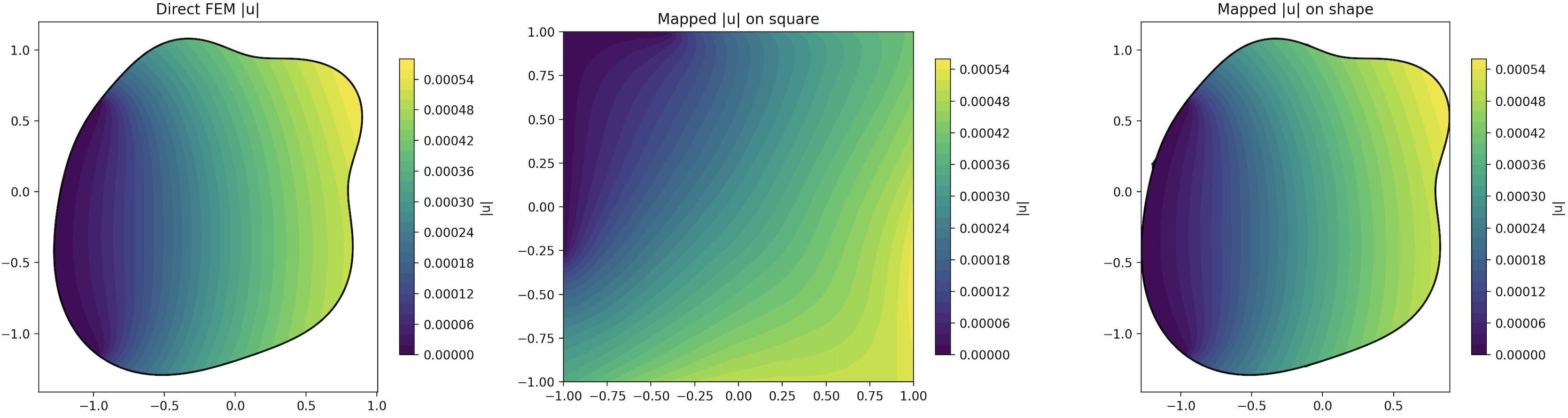}
        \caption{Displacement magnitude obtained with the direct FEM, pullback FEM on the reference square, and pullback FEM transferred to the physical geometry.}
        \label{fig:elasticity_bestResult}
    \end{subfigure}
    
    \begin{subfigure}[b]{0.68\textwidth}
        \centering
        \includegraphics[width=\textwidth]{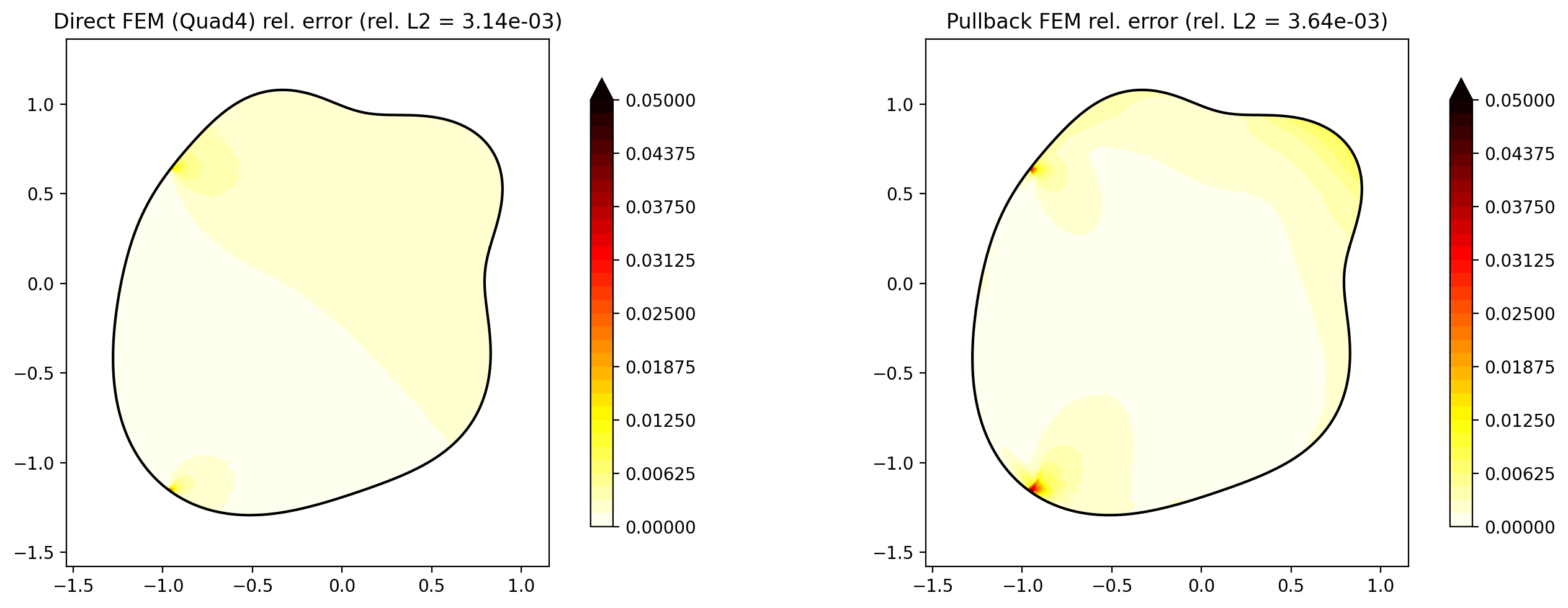}
        \caption{Corresponding pointwise error fields relative to the
high-fidelity reference solution.}
        \label{fig:elasticity_bestResultComparison}
    \end{subfigure}

    \caption{Best relative-performance case (Shape~4)}
    \label{fig:example3_resultsBestCase}
\end{figure}

\begin{figure}[H]
    \centering

    \begin{subfigure}[b]{0.68\textwidth}
        \centering
        \includegraphics[width=\textwidth]{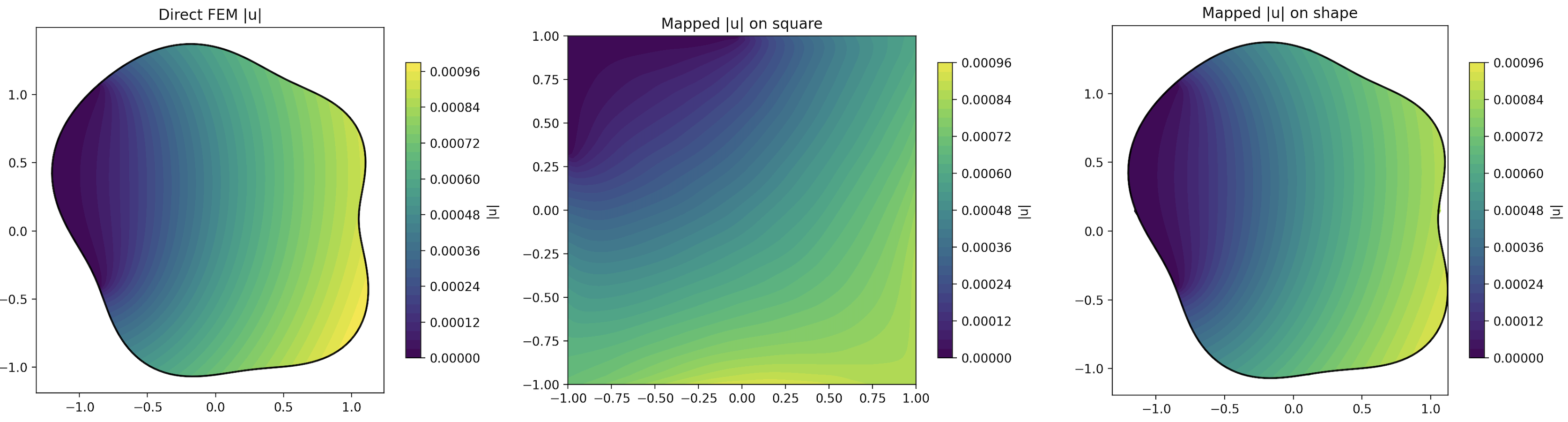}
        \caption{Displacement magnitude obtained with the direct FEM, pullback FEM on the reference square, and pullback FEM transferred to the physical geometry.}
        \label{fig:elasticity_worstResult}
    \end{subfigure}
    
    \begin{subfigure}[b]{0.68\textwidth}
        \centering
        \includegraphics[width=\textwidth]{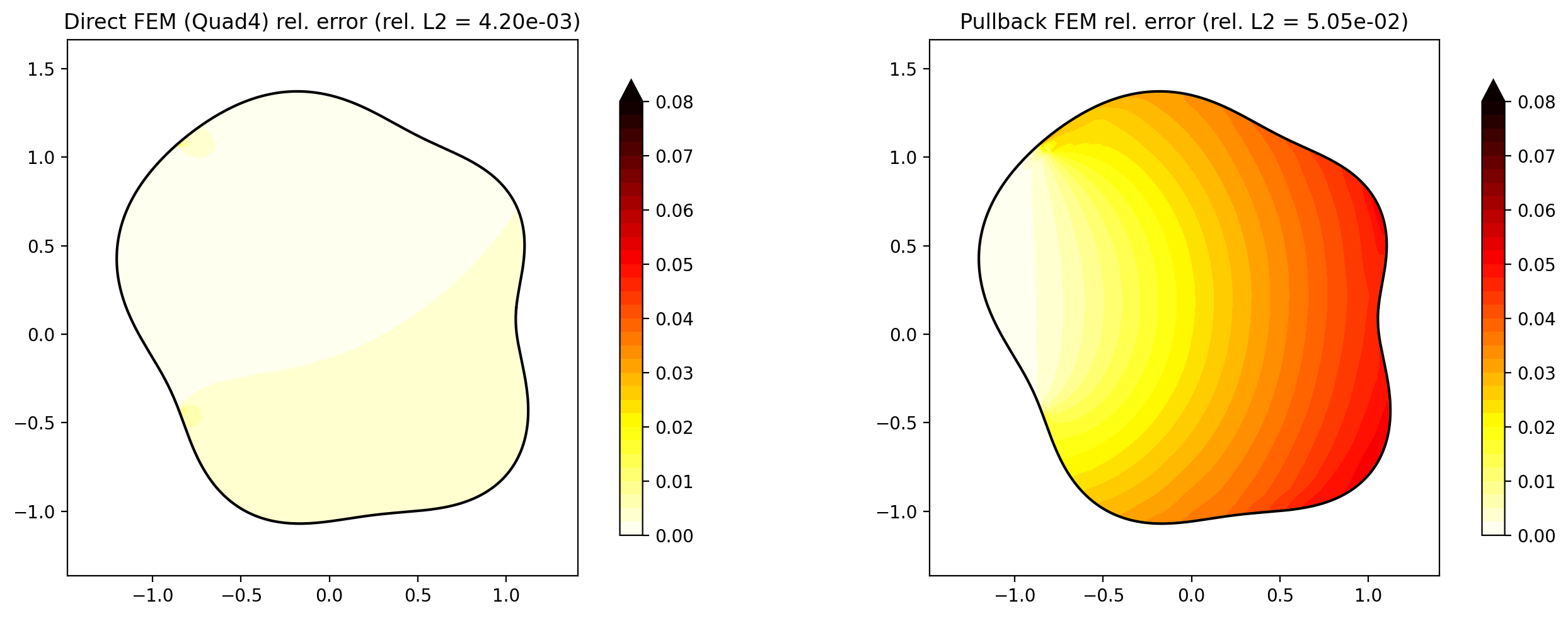}
        \caption{Corresponding pointwise error fields relative to the
high-fidelity reference solution.}
        \label{fig:elasticity_worstResultComparison}
    \end{subfigure}

    \caption{Largest relative discrepancy case (Shape~3)}
    \label{fig:example3_resultsWorstCase}
\end{figure}

Figures~\ref{fig:example3_resultsBestCase}--\ref{fig:example3_resultsWorstCase} illustrate the two representative cases corresponding to the best and the worst relative performance of the pullback formulation. These correspond to Shapes~4 and 3, respectively. For each case, the displacement magnitude obtained with the direct FEM is compared with the pullback solution both on the reference square and after transfer to the physical domain. The corresponding pointwise error fields with respect to the high-fidelity reference solution are also shown.

\begin{figure}[H]
    \centering
    \includegraphics[width=0.6\textwidth]{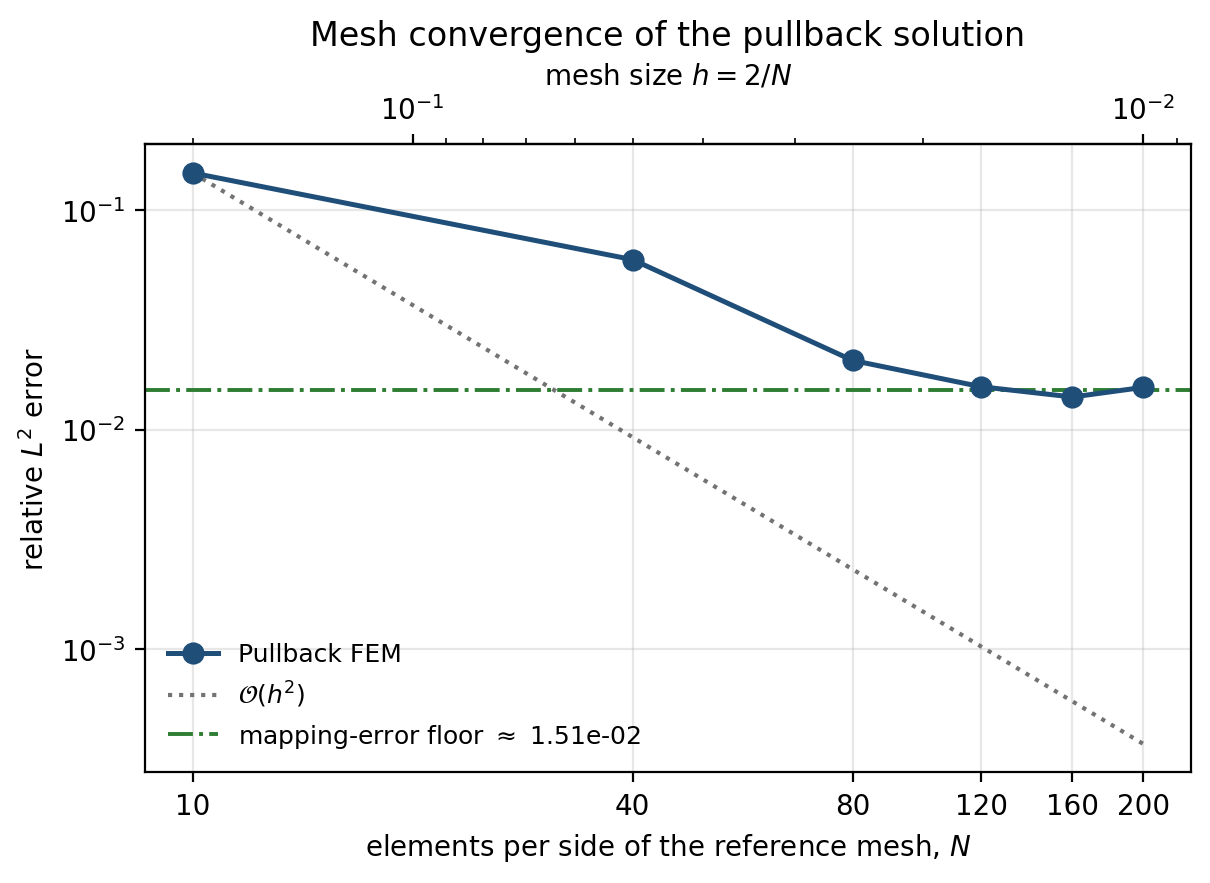}
    \caption{Convergence of the pullback finite element solution under refinement
of the reference mesh, for the first Fourier blob}
    \label{ex3_pullback_convergence}
\end{figure}

As a final assessment of the pullback formulation we examine its behaviour under refinement of the reference mesh. Fig.~\ref{ex3_pullback_convergence} reports the relative $L^2$ error of the pullback solution for the first Fourier blob as $\mathcal{T}^0$ is refined from $10\times10$ to $200\times200$ elements, with the learned mapping and the reference solution held fixed. The error decreases from
$1.48\times10^{-1}$ to $1.41\times10^{-2}$ and then ceases to improve. This behaviour is consistent with the error decomposition discussed in Sec.~\ref{subsec:solution_equivalence}. For a fixed diffeomorphism, the pullback is an exact reformulation of the boundary-value problem posed on $G(\Omega^0)$. However, this learned domain does not in general coincide exactly with the target domain $\Omega$. Refinement of $\mathcal{T}^0$ therefore reduces the finite element discretisation error, while the discrepancy associated with the learned geometric representation remains essentially independent of $h$. Once the discretisation contribution becomes sufficiently small, further refinement consequently produces little improvement in the total error.

The observed plateau can therefore be interpreted as an empirical estimate of the mapping-induced error floor for this problem, together with any remaining quadrature, interpolation, and reference-solution errors. Two observations follow. First, beyond a certain level of mesh refinement, the attainable accuracy becomes primarily limited by the quality of the learned geometric representation rather than by the finite element discretisation. Second, since the mapping is constructed independently of the PDE solve, its geometric quality can be assessed beforehand using the diagnostics introduced in Sec.~\ref{subsubsec:inference}, providing an a priori indication of mapping quality, although the resulting effect on the PDE solution remains problem-dependent.

\section{Conclusions}
\label{sec:conclusions}

In this work, we have presented a framework for solving partial differential equations on varying geometries by learning a diffeomorphic map from a fixed reference square to the physical domain using a Neural Spline Flow, and solving the corresponding pulled-back weak form on a structured reference mesh. The rational-quadratic coupling architecture guarantees a strictly positive Jacobian determinant, so the learned map remains invertible and orientation-preserving by construction. For a given learned diffeomorphism, the pullback constitutes an exact reformulation of the boundary-value problem on the domain represented by the learned map, preserving the mathematical structure of the original problem. In practical applications, the accuracy of this representation depends on how closely the learned map reproduces the target geometry.

An auto-decoder combined with a lightweight encoder extends the approach to families of geometries, enabling the mapping for an unseen shape to be inferred in a single forward pass. Unlike classical surrogate modelling, the neural network does not approximate the PDE solution itself, but only the geometric transformation, while the governing physics remains handled by the finite element solver. The overall accuracy of the approach therefore depends both on the finite element discretisation and on the quality of the learned mapping. The numerical results showed that residual mapping inaccuracies and distortion may propagate through the transformed coefficients and contribute to the overall solution error. However, for the geometries examined in this work the geometric approximation remained sufficiently accurate to yield reliable finite element solutions. 

The use of a fixed structured reference mesh removes the need to generate a new boundary-conforming mesh for each geometry and also creates opportunities for efficient structured-grid implementations, including reusable algebraic structures, tensor-product quadrature and/or matrix-free operator evaluation. The main limitations, however, arise for increasingly complex geometries, where a single reference-to-physical map may become strongly distorted and lead to highly anisotropic transformed coefficients. Moreover, multiply-connected domains cannot be represented by a diffeomorphism from a simply connected square, while in the multi-shape setting generalisation naturally depends on the geometric diversity represented in the training distribution.

Future work will therefore focus on improving mapping accuracy and robustness, quantifying and reducing the resulting geometric error, extending the framework through  multi-patch approaches or domain decomposition to more complex topologies, and investigating three-dimensional implementations. These developments could extend the proposed approach beyond the two-dimensional setting considered here toward more general computational mechanics applications.

\appendix

\bibliographystyle{elsarticle-num} 
\bibliography{myBib}

\end{document}